\documentclass[draft]{article}
\def\today{24.1.12} 
\usepackage{amsmath,amsfonts,amsthm,amssymb,amscd}

\renewcommand{\Im}{\mathop{\rm Im}\nolimits}

\newcommand{\beq}{\begin{equation}}
\newcommand{\ee}{\end{equation}}

\theoremstyle{plain} \newtheorem{theorem}{Theorem}[section]
\newtheorem{lemma}[theorem]{Lemma}
\newtheorem{proposition}[theorem]{Proposition}
 \theoremstyle{definition}
\newtheorem{definition}[theorem]{Definition} \theoremstyle{remark}
\newtheorem{remark}[theorem]{Remark}

\newcommand{\R}{{\mathbb R}} \newcommand{\U}{{\mathcal U}}

\newcommand{\Z}{{\mathbb Z}}

\newcommand{\N}{{\mathbb N}}

\newcommand{\Tr}{{\mathcal T}}

\newcommand{\Ph}{{\mathcal P}}

\newcommand{\resto}{{\mathcal R}}
\def\im{{\rm i}}

\newcommand{\C}{\mathbb{C}}

\def\uno{{\kern+.3em {\rm 1} \kern -.22em {\rm l}}}

\def\norma#1{\left\| #1\right\|}

\numberwithin{equation}{section}

\begin{document}

  \title{On asymptotic stability of moving  ground
 states    of the nonlinear Schr\"odinger equation}

 \author {Scipio Cuccagna}

\date{\today}
 \maketitle
\begin{abstract} We extend to the case of moving solitons, the result on asymptotic
stability of ground states of the NLS obtained by the author in \cite{Cu1}.
For technical reasons we consider only smooth solutions. The
 proof is similar to  the earlier paper. However now the flows
 required for the Darboux Theorem and the Birkhoff normal forms, instead of falling within the framework of standard  theory of ODE's, are   related to quasilinear hyperbolic symmetric systems.  It is also not obvious that   Darboux Theorem can be applied,
  since we need to compare   two symplectic forms   in a neighborhood of the ground states  not in
  $H^{1}(\R ^3)$, but rather in the space $\Sigma $ where
  also the variance is bounded.  But the NLS does not preserve small
  neighborhoods of the ground states in $\Sigma $.
\end{abstract}

\section{Introduction}

We consider the   nonlinear
Schr\"odinger equation (NLS)

\begin{equation}\label{NLS}
 \im u_{t }=-\Delta u   +\beta  (|u|^2) u  \quad , \quad u(0,x)=u_0(x), \, (t,x)\in\mathbb{ R}\times
 \mathbb{ R}^3
\end{equation}
assuming:  $\beta (|u|^2)u$ is  "short range"
  and smooth;  \eqref{NLS} has a smooth family of ground states.

  In \cite{Cu1} we proved that, if we break the translation invariance of the equations by either taking solutions with $ u_0(-x) \equiv u_0(x)$ or by
  introducing some spacial inhomogeneity in the equation, for instance by adding
  a term $V(x)u$, the ground states are asymptotically stable, under what looks
  a generic hypothesis. More precisely, we assumed  the sufficient conditions for
 orbital stability by M. Weinstein \cite{W1}. We assumed  some spectral hypotheses on the linearizations
 (absence of embedded eigenvalues, this probably   always true under Weinstein's
  hypotheses) and a number
 of other hypotheses which hold generically (nondegeneracy of the
 thresholds of the continuous spectrum;   some mild
 non resonance conditions on the eigenvalues) and which are stated
 in Subsect. \ref{subsec:statement}. We then proved a form of Fermi
 golden rule (FGR). Specifically, we proved   that certain coefficients of the system
 are square powers. This implies that they are non negative.  We then assumed
 that these coefficients are in fact positive, which is probably true
 generically, and used this to prove asymptotic stability of the ground states.
A   result similar  to
\cite{Cu1}, with some restrictions, is proved for the Dirac equation in  \cite{boussaidcuccagna}.

 In this paper we extend the proof of \cite{Cu1}  to equations
 like  \eqref{NLS}  without requiring any symmetry for the initial data $u_0$.
  Hence we prove that  a solution $u(t)$
of \eqref{NLS} starting sufficiently close to ground states, is
asymptotically as $t\nearrow +\infty$ of the form $e^{i\theta (t) +\frac{\im }{2} v_+\cdot x} \phi _{\omega _+} (x -D(t))+
e^{it\Delta }h_+$, for $\omega _+$  and $v_+$ fixed,  for $\theta\in C^1(\R , \R)$, $D\in C^1(\R , \R ^3)$ and
  for $h_+\in
H^1(\mathbb{R}^3)$ a small energy function. For technical reasons we need
a certain known a priority regularity  and decay of  $u_0$, although
in the estimates we use only the norm $\norma{u_0}_{H^1}$.

The problem of
stability of ground states of the NLS has a long history, partially reviewed in
\cite{Cu1}. Orbital stability  was  well understood in the 80's, see
\cite{CL,shatah,W1,GSS1,GSS2}, and there is a long literature.
For asymptotic  stability we highlight  \cite{SW1,SW2,BP1,BP2,zhousigal}, for more references see \cite{Cu1}.

    One of the crucial difficulties  in   asymptotic stability
 is in showing that the discrete
 modes appearing naturally in the problem,  which left on their own would   oscillate,
 are dumped by the nonlinear interaction with the continuous modes.
 By conservation of energy, this happens by
    passage of energy from the discrete
 to the continuous modes.   The proof
  requires nonlinear versions of the  FGR,   see \cite{sigal}. In our setup,
   versions of the FGR  of ever growing generality where proved in
 special cases in \cite{BP2,SW3,TY1,TY2,TY3,T,BS,Cu2,zhousigal,cuccagnamizumachi,zhouweinstein1}.
They involved search of appropriate coordinates through Poincar\'{e} Dulac normal forms.
For related linear problems, see \cite{KW} and the references therein.
 \cite{Cu3}
seems to be the first reference to have noticed the relevance of the hamiltonian structure
of the NLS. The intuition in \cite{Cu3} was implemented in the series
  \cite{bambusicuccagna,Cu1,Cu4,boussaidcuccagna}.

The FGR consists in proving that certain coefficients   are
square powers, and so are generally positive.  The coefficients
will typically have the structure $A\cdot \overline{B}$, with $A$ and $B$ coefficients
of the system in appropriate coordinate systems. The square power structure will follow from $B=A$.
Proving such identities among the coefficients in the NLS, is certainly
easier if we exploit the hamiltonian structure.  We search an appropriate  system of coordinates through the method of Birkhoff normal forms. In the cases considered in  \cite{bambusicuccagna,Cu4} this
  is   easier because the natural coordinates which appear linearizing
  the system around the 0 solution,  are canonical coordinates. So one can start the Birkhoff normal forms from the initial system of coordinates.

  In analysis of    the stability of solitons, the   natural coordinates of the linearization
  are not canonical.  Before starting the method of normal
  forms one has to find canonical coordinates, through an implementation
  of the Darboux theorem. This has to be done in the right non abstract way,
  in order not to lose the property that the NLS is a   semilinear
  system. This process is  done in \cite{Cu1,boussaidcuccagna}.
  However these papers, as well as most of  the papers quoted so far,
  avoid the case of moving solitons. Special cases,
  without discrete modes, were treated in \cite{BP1,Cu5}. For
   multisolitons with weak interaction and no discrete modes see
in \cite{P,RSS}.

Moving solitons present   three special difficulties. First of all, they yield terms in the equation of the continuous modes which are non linear but which is difficult
to see as perturbations of the  linear equation. It is not obvious how to eliminate them through an integrating factor.
Fortunately
work by Beceanu, such as \cite{beceanu}, has solved this problem. Early solutions
in particular cases are in \cite{BP1,BP2} (see \cite{BS,Cu5} for proofs).

The second difficulty involves  the Darboux theorem. The method followed in \cite{Cu1}
 becomes too complicated in the moving solitons setting. It   is useful to use charge and linear momenta  as  coordinates. In the case of the charge, had this been done in \cite{Cu1},
 it would have simplified the proof there.  One difficulty with the
  Darboux theorem is  the  determination of the vectorfield
$\mathcal{X}^t$ obtained as dual of an appropriate 1 form, in the Moser
version of Darboux Theorem used here.  In  \cite{Cu1} the existence of such $\mathcal{X}^t$ and some of its properties are rather elementary. In this paper,  we are comparing two symplectic forms which are not   both defined in
$H^1(\R ^3)$. Rather, they are  symplectic forms in the smaller space $\Sigma _1$
 formed by functions of bounded  $H^1$ norm and bounded variance,
see \eqref{eq:sigma} and Sect.\ref{sec:sigma}. The proof of the existence of
$\mathcal{X}^t$  would be easy
if we could assume that the variance of the solutions of the NLS,
assuming it is small at time $t=0$, remained small  for all times.
But this is not the case, so the discussion is rather complicated.

The third difficulty present here and not in    \cite{Cu1} is that, the vector fields whose flows are
 used to change coordinates in the implementation of
 Darboux theorem and of the method of Birkhoff normal forms, do not
 fall as   in  \cite{Cu1} within the framework
 of smooth vectorfields in Banach spaces.
  Here instead we have to deal with
 quasilinear hyperbolic symmetric systems. So well posedness and
 regularity  of the flows, which in   \cite{Cu1} are elementary,
 are here more delicate. Particular attention requires the issue of regularity
 of the flows with respect to the initial data.
  Fortunately our systems have quite simple structure.
  In a rather standard way, our flows are obtained as limits
  of flows of systems with viscosity, which fall within the classical
  framework of ODE's.  In the limit we lose some regularity.
  It is at this juncture that we use
the qualitative information on regularity and decay of the initial
datum $u_0$. The more we iterate, the more we lose regularity.
Fortunately we have as much regularity and decay of $u_0$ and of the
ground states as we want, to start with.

{\bf Acknowledgments} I wish to thank G. Tondo for
discussions about the reduction of variables and for pointing out reference \cite{Olver}.

\subsection{Statement of the main result}
\label{subsec:statement}

We will assume the following hypotheses.

\begin{itemize}
\item[(H1)] $\beta  (0)=0$, $\beta\in C^\infty(\R,\R)$.
\item[(H2)] There exists a $p\in(1,5)$ such that for every
$k\ge 0$ there is a fixed $C_k$ with
$$\left| \frac{d^k}{dv^k}\beta(v^2)\right|\le C_k
|v|^{p-k-1} \quad\text{if $|v|\ge 1$}.$$

\item[(H3)]
There exists an open interval $\mathcal{O}$ such that
\begin{equation}
  \label{eq:B}
  \Delta u -\omega u+\beta(|u|^2)u=0\quad\text{for $x\in \R^3$},
\end{equation}
admits a $C^1$-family of ground states $\phi _ {\omega }(x)$ for
$\omega\in\mathcal{O}$.
\item [(H4)]
\begin{equation}
  \label{eq:1.2}
\frac d {d\omega } \| \phi _ {\omega }\|^2_{L^2(\R^3)}>0
\quad\text{for $\omega\in\mathcal{O}$.}
\end{equation}
\item [(H5)]
Let $L_+=-\Delta   +\omega -\beta (\phi _\omega ^2 )-2\beta '(\phi
_\omega ^2) \phi_\omega^2$ be the operator whose domain is $H^2
 (\R^3)$. Then we assume that $L_+$ has exactly one negative eigenvalue and
does not have kernel when restricted to $H^1_r(\R ^3)$, the subspace of $H^1 (\R ^3)$ formed by functions with radial symmetry.

\item [(H6)] Let $\mathcal{H}_\omega$ be the linearized operator around $ \phi_\omega$
(see Section \ref{section:linearization} for the precise
definition). $\mathcal{H}_\omega$ has  $m$
    positive eigenvalues $\lambda _1(\omega )\le \lambda _2(\omega )
\le ...\le \lambda _m(\omega )$
    with $0<N_j\lambda _j(\omega )<
\omega < (N_j+1)\lambda _j(\omega )$ with $N_j\ge 1$. We set
$N=N_1$. Here we are repeating each  eigenvalue a number of times
equal to its  multiplicity. We assume the multiplicity constant in $\omega$.
\item [(H7)]  There is no multi index $\mu \in \mathbb{Z}^{m}$
with $|\mu|:=|\mu_1|+...+|\mu_k|\leq 2N_1+3$ such that $\mu \cdot
\lambda =m$.

\item[(H8)] If $\lambda _{j_1}<...<\lambda _{j_k}$ are $k$ distinct
  $\lambda$'s, and $\mu\in \Z^k$ satisfies
  $|\mu| \leq 2N_1+3$, then we have
$$
\mu _1\lambda _{j_1}+\dots +\mu _k\lambda _{j_k}=0 \iff \mu=0\ .
$$
\item[(H9)] $\mathcal{H}_\omega$ has no other eigenvalues except for $0$ and
the $ \pm \lambda _j (\omega )$. The points $\pm \omega$ are not
resonances.

\item [(H10)]
The Fermi golden rule  Hypothesis (H10)   in subsection
\ref{subsec:FGR}, see \eqref{eq:FGR}, holds.

\item [(H11)] We assume that $u_0$ is a Schwartz function.

\end{itemize}
Recall that from the $\phi _\omega$ one can derive solitons
$e^{\frac i2 v\cdot x -\frac i4 |v|^2 t +it\omega +i\gamma } \phi
_\omega (x-vt- D)$. Solutions of \eqref{NLS} starting close to a
ground state, for some time can be written as
 \begin{equation} \label{eq:Ansatz1}\begin{aligned}&  u(t,x) = \tau _{D(t)}e^{i \Theta (t,x)}
(\phi _{\omega (t)} (x)+ r(t,x))  \\& \Theta (t,x)=\frac
12 v(t)\cdot  x  + \vartheta (t)
\end{aligned}
\end{equation}
with $\tau _{D }f(x):= f(x-D).$

\begin{theorem}\label{theorem-1.1}
  Let $\omega_1\in\mathcal{O}$, $v_1\in \R$   and $\phi_{\omega_1}(x)$
  a ground state of \eqref{NLS}. Let $u(t,x)$ be a solution to
\eqref{NLS}. Assume (H1)--(H10). Then, there exist an $\epsilon_0>0$
and a $C>0$ such that if $\varepsilon:=\inf_{\gamma \in \R , y\in
\R^3}\|u_0-e^{\im \gamma} e^{\   \frac  \im 2{v_1 \cdot x}  } \phi_ {\omega _1 }(\cdot -y)
\|_{H^1}<\epsilon_0,$ there exist $\omega _\pm\in\mathcal{O}$,
$v_\pm \in \R^3$, $\theta\in C^1(\R;\R)$, $ y\in C^1(\R;\R ^3)$ and
$h _\pm \in H^1$ with $\| h_\pm \| _{H^1}+|\omega _\pm
-\omega_1  |+ |v_\pm -v_1|\le C \varepsilon $ such that

\begin{equation}\label{eq:scattering}
\lim_{t\to  \pm\infty}\|u(t,\cdot)-e^{\im \theta(t)+  \frac  \im 2{v_\pm \cdot x}  }\tau _{y(t)}\phi_{\omega
_\pm} - e^{\im t\Delta }h _\pm    \|_{H^1}=0.
\end{equation}
In the notation of \eqref{eq:Ansatz1},
 we have  $\tau _{D(t)}e^{i \Theta (t,x)}
 r(t,x)=A(t,x)+\widetilde{r}(t,x)$  such that
$A(t, \cdot ) \in \mathcal{S}(\R^3, \C)$,
$|A(t,x)|\le C (t)$ with  $\lim
_{|t|\to \infty }C (t)=0$ and such that for any pair $(p,q)$ which
is admissible, by which we mean that
\begin{equation}\label{admissiblepair}  2/p+3/q= 3/2\,
 , \quad 6\ge q\ge 2\, , \quad
p\ge 2,
\end{equation}
we have
\begin{equation}\label{Strichartz} \|  \widetilde{r} \|
_{L^p_t( \mathbb{R},W^{1,q}_x)}\le
 C\|  u_0  \| _{H^1 }.
\end{equation}

\end{theorem}

\begin{remark}
\label{rem:h12} In the proof we use only bounds on the $H^1$ norm of $u_0$.
Nonetheless we use in a qualitative fashion the fact that $u_0$ is
quite regular and rapidly decaying. For simplicity we restrict attention to the case when  $u_0$  is a Schwartz function.
\end{remark}

\begin{remark}
\label{rem:parameters}   In the proof we show that we can take $\theta =\vartheta $
and $y=D$, with $( \vartheta ,D)$ the functions in Lemma \ref{lem:modulation}.
 In Lemma \ref{lem:Conv der}
 we show $\dot D =v+o(1)$ with $v$ as in Lemma \ref{lem:modulation} and $\lim _{t\to \infty}o(1)=0$. Similarly, $\dot \vartheta =\omega + \frac{v^2}{4}+o(1)$.
\end{remark}

\begin{remark}
\label{rem:embedded}  Notice that in (H6) we exclude eigenvalues of $\mathcal{H}_{\omega}$ in $(\omega , \infty )$ mainly because we think they
do not exist under (H1)--(H5). Notice that   in \cite{CPV} smoothing estimates for $\mathcal{H}_{\omega}$ are proved also in the presence of eigenvalues   in $(\omega , \infty )$, so that the theory here and in \cite{Cu1} could be developed also in that situation.
\end{remark}

\begin{remark}
\label{rem:moment} By elementary arguments, Theorem \ref{theorem-1.1} is a consequence of the
special case where the linear momenta are equal to 0, see \eqref{eq:charge}.
So we will focus only on this case.
\end{remark}

We briefly describe   the proof, which is similar in spirit to \cite{Cu1},
but departs from \cite{Cu1} in important ways.
First of all, we need to choose a system of coordinates around the ground states. There is a natural choice related to the notion of modulation
and to the spectral decomposition of the linearization. Only in a second moment we use
charge and linear moment as coordinates. Since these are invariants of motion,  we then consider a reduction of coordinates,
in an elementary fashion. We also move to canonical coordinates through
an implementation of the Darboux theorem.  We then start the
Birkhoff normal form argument, that is, we consider other canonical
coordinates where the system  looks increasingly more treatable. Finally
after a finite number of them, we settle  with coordinates where
it is possible to prove the Fermi golden rule. Then, if (H10) is true,
we conclude simultaneously that   the continuous modes
disperse and that the energy of the discrete modes leaks away
through nonlinear interaction with the continuous modes.
The most delicate and novel feature of this paper with respect to \cite{Cu1} consists in
the analysis of the flows $\phi ^t$ used for Darboux theorem and the
Birkhoff normal forms. In particular, since we are outside the realm
of ODE's, it is less obvious to conclude that for fixed $t$ the flow
 $\phi ^t$ is a differentiable map. This is where   (H11) is helpful.
As for  Birkhoff normal forms, we also  add some more material useful to understand
the homological equations, which are nonlinear, and which should help
to understand the analogous discussion in \cite{Cu1}, which is very succinct.

We end the introduction with some notation. Given two functions
$f,g:\mathbb{R}^3\to \mathbb{C}$ we set $\langle f|g\rangle = \int
_{\mathbb{R}^3}f(x)  g(x)  dx$ (with no complex conjugation). Given a matrix $A$, we denote by
  by $ A^T $,  its transpose. Given two vectors $A$ and $B$,
we denote by  $A^T B=\sum _j A_jB_j$ their inner product. Sometimes
we omit the summation symbol, and we use the convention on sum over
repeated indexes.  For any $k,s\in \mathbb{R}$ and any Banach space $K$,
 we set

 \[ H^{ k,s}(\mathbb{R}^3,K)=\{ f:\mathbb{R}^3\to K \text{ s.t.}
 \| f\| _{H^{s,k}}:=\| \langle x \rangle ^s \| (-\Delta +1)^{k} f
 \| _{K}\| _{L^2
 }<\infty \}.\]
 In particular we set
 $L^{2,s} =H^{0,s}  $, $L^2=L^{2,0}    $,  $H^k=H^{k,0}    $.
Sometimes, to emphasize that these spaces refer to spatial
variables, we will denote them by $W^{k,p}_x$, $L^{ p}_x$, $H^k_x$,
$H^{ k,s}_x$ and $L^{2,s}_x$. For $I$ an interval and $Y_x$ any of
these spaces, we will consider Banach spaces $L^p_t( I, Y_x)$ with
mixed norm $ \| f\| _{L^p_t( I, Y_x)}:= \| \| f\| _{Y_x} \| _{L^p_t(
I )}.$ In the course of the proof we will consider a fixed pair of  spaces $H^{K,S}$ and $H^{-K,-S}$,
for positive and large $K$ and $S$.

We set $ (\im  \partial _x  + \im x)  ^\alpha := \prod _{a=1}^{3}(\im  \partial _a  + \im x_a)^{\alpha _a} $   for any multiindex $\alpha$.  For any natural number $n\ge 1$ We consider the space $\Sigma _n$ defined by
 \begin{equation}\label{eq:sigma}
\begin{aligned} &
      \| U \| _{\Sigma _n} ^2:= \sum _{|\alpha | \le n}\|  (\im  \partial _x  + \im x)  ^\alpha U \| _{L^2} ^2 <\infty  .    \end{aligned}
\end{equation}

Given an operator $A$, we will denote by $R_A(z)=(A-z)^{-1}$
its resolvent. We set $\mathbb{N}_0=\mathbb{N}\cup \{0 \}$. We will
consider multi indexes $\mu \in \mathbb{N}_0^m$. For   $\mu \in
\mathbb{Z}^m$ with $\mu =(\mu _1,..., \mu _m)$ we set $|\mu |=\sum
_{j=1}^m |\mu _j|.$     For $X$ and $Y$ two Banach space, we will
denote by $B(X,Y)$ the Banach space of bounded linear operators from
$X$ to $Y$ and by $B^{\ell}(X,Y)= B (  \prod _{j=1}^\ell X ,Y)$.   Given a differential form
$\alpha$, we denote by
    $d\alpha$ its exterior differential.

\section{Linearization and set up}
\label{section:linearization}

 Let $U={^t(u,\overline{u})}$. We consider the energy

 \begin{equation} \label{eq:energyfunctional}\begin{aligned}&
 E(U)=E_K(U)+E_P(U)\\&
E_K(U):= \int _{\R ^3}
  \nabla u \cdot \nabla \overline{u} dx  \, , \quad
E_P(U):=
 \int _{\R ^3}B( u \overline{u}) dx \end{aligned}
\end{equation}
with $B(0)=0$ and $\partial _{\overline{u}}B(|u|^2)=\beta (|u|^2)u$.
We will consider the matrices \begin{equation}
\label{eq:Pauli}\begin{aligned} &\sigma _1=
\begin{pmatrix}0 &
1  \\
1 & 0
 \end{pmatrix} \, ,
\sigma _2=\begin{pmatrix}  0 &
\im  \\
-\im & 0
 \end{pmatrix} \, ,
\sigma _3=\begin{pmatrix} 1 & 0\\0 & -1 \end{pmatrix} .
\end{aligned}
\end{equation}
  For    $U\in H^1( \R ^3, \C )\times H^1( \R ^3, \C )$  we have the following charge and momenta, which yield invariants of motion of \eqref{NLS}: \begin{equation}\label{eq:charge}\begin{aligned} &Q(U)=  \int
_{\R ^3}u \overline{u} dx= \frac{1}{2}\langle U| \sigma _1 U\rangle
.\\&  \Pi _a(U)=   \Im \int _{\R ^3}  \overline{u }(x) u_{x_a} (x)
dx= -\frac{\im }{2} \langle U| \sigma _3\sigma _1
\frac{\partial}{\partial x_a}U\rangle   . \end{aligned}
\end{equation}
Sometimes we will denote $\Pi  (U)=(\Pi _1(U),\Pi _1(U),\Pi _3(U))$.
We will focus only on  solutions of the NLS \eqref{NLS}  s.t. $ \Pi  (U)=0$.

The charge $Q$  and the momenta $\Pi _a$ are   in
$C^\infty ( H^1( \R ^3, \C )\times H^1( \R ^3, \C ) , \C  )$ while $E\in C^1 ( H^1( \R ^3, \C )\times H^1( \R ^3, \C ) , \C  )$. If for any
such functional $F$ we set $dF(X)=\langle \nabla F, X\rangle $  for any $X\in \C^2$, with $dF$ the exterior differential and $\nabla F$ the gradient of $F$, then
 \begin{equation}\label{eq:charge1}\begin{aligned} &  \nabla Q(U)=\sigma _1 U
\, , \quad  \nabla \Pi _a(U)=- \im \sigma _3\sigma _1
\frac{\partial}{\partial x_a}U . \end{aligned}
\end{equation}
For later use we set
\begin{equation} \label{eq:function q} \begin{aligned} &\Phi _\omega =\begin{pmatrix} \phi _\omega
   \\ \phi _\omega
 \end{pmatrix}   , \, q(\omega )=Q(\Phi _\omega ),
\, e (\omega )=E(\Phi _\omega ), \,  p_a (\omega )= \Pi _a(\Phi
_\omega)\\& d(\omega )=e (\omega )+\omega q(\omega ) .\end{aligned}
\end{equation}
   Equation
\eqref{NLS} can be written as
\begin{equation}\label{eq:NLSvectorial} \im \dot U =
\begin{pmatrix} 0 &1
   \\ -1 & 0
 \end{pmatrix}   \begin{pmatrix} \partial _{u}E
   \\ \partial _{\overline{u}}E
 \end{pmatrix}  = \sigma _3 \sigma _1 \nabla E (U).
\end{equation}
We introduce now the \textit{linearization}
\begin{equation}  \label{eq:linearization} \begin{aligned} &
\mathcal{H}_\omega :=\sigma _3 \sigma _1 \left ( \nabla ^2 E (\Phi
_\omega
   ) +\omega \nabla ^2Q (\Phi _\omega
   )\right )   = \sigma_3(-\Delta +\omega)
  +V  _{\omega}  \\& \text{where }
  V  _{\omega} :=
  \sigma_3
\left[\beta (\phi ^2_{\omega }) +\beta ^\prime (\phi ^2_{\omega
})\phi ^2_{\omega } \right] +\im  \sigma _2 \beta ^\prime (\phi ^2
_{\omega })\phi ^2 _{\omega } .\end{aligned}
\end{equation}
The essential spectrum of $\mathcal{H}_\omega$ is $(-\infty
, -\omega ]\cup [ \omega,+\infty )$. It is well  known,
\cite{W2}, that  by  (H4)--(H5)   $0$ is an isolated eigenvalue of
$\mathcal{H}_\omega$ with $\dim N_g(\mathcal{H}_\omega)=8$ and
\begin{equation}\label{eq:Kernel} \begin{aligned} &
\mathcal{H}_\omega\sigma_3\Phi_\omega=0= \mathcal{H}_\omega\partial
_{x_j}\Phi_\omega,\\&  \mathcal{H}_\omega\partial_\omega\Phi_\omega
=-\Phi_\omega   \, , \quad \mathcal{H}_\omega x_a\sigma_3\Phi_\omega
=-\partial _{x_a}\Phi_\omega .\end{aligned}
\end{equation}
Since $\mathcal{H}_\omega^*=\sigma_3\mathcal{H}_\omega\sigma_3$, we
have $N_g(\mathcal{H}_\omega^*)=\operatorname{span}\{\Phi ,
\sigma_3\partial_\omega\Phi , \sigma_3\partial
_{x_a}\Phi,x_a \Phi\}$. We consider
  eigenfunctions $\xi _j(\omega)$   with eigenvalue   $\lambda
_j(\omega)$:
$$
\mathcal{H}_\omega\xi _j(\omega)=\lambda _j(\omega)\xi
_j(\omega),\quad \mathcal{H}_\omega\sigma_1\xi _j(\omega)=-\lambda
_j(\omega)\sigma_1\xi _j(\omega) .$$ They can be normalized so that
$\langle \sigma_3 \mathcal{H}_\omega \xi _j(\omega),\overline{\xi}
_\ell(\omega) \rangle =\delta _{j\ell }$. Furthermore, they can be chosen to be
real, that is with real entries, so $\xi _j=\overline{ \xi} _j$ for
all $j$, see Prop. 5.1 \cite{zhouweinstein1}.

Both $\phi_\omega$ and $\xi _j(\omega,x)$ are smooth in
$\omega\in\mathcal{O}$ and $x\in\R^3$ and satisfy
$$\sup_{\omega\in\mathcal{K},x\in\R^3} e^{a|x|}( |\partial ^\alpha _x\phi_\omega(x)|+
\sum _{j=1}^{m}|\partial ^\alpha _x \xi _j(\omega,x)|) <\infty$$ for
every
$a\in(0,\inf_{\omega\in\mathcal{K}}\sqrt{\omega-\lambda(\omega)})$
and every compact subset $\mathcal{K}$ of $\mathcal{O}$.

For $\omega\in\mathcal{O}$, we have the
$\mathcal{H}_\omega$-invariant Jordan block decomposition
\begin{align}  \label{eq:spectraldecomp} &
L^2(\R^3,\C^2)=N_g(\mathcal{H}_\omega)\oplus \big (\oplus _{\pm}
\oplus _{j=1}^m \ker (\mathcal{H}_\omega\mp \lambda _j(\omega))
\big)\oplus L_c^2(\mathcal{H}_\omega),
\end{align}
  $L_c^2(\mathcal{H}_\omega):=
\left\{N_g(\mathcal{H}_\omega^\ast)\oplus \big (\oplus _{\lambda \in
\sigma _d\backslash \{ 0\}}   \ker (\mathcal{H}_\omega ^*- \lambda
(\omega)) \big)\right\} ^\perp $ with $\sigma _d =\sigma _d
(\mathcal{H}_\omega)$. We also set $L_d^2(\mathcal{H}_\omega):=
N_g(\mathcal{H}_\omega)\oplus \big (\oplus _{\lambda \in \sigma
_d\backslash \{ 0\}}   \ker (\mathcal{H}_\omega  - \lambda (\omega))
\big ) .$ By $P_c(\mathcal{H}_{\omega})$ (resp.
$P_d(\mathcal{H}_{\omega})$), or simply by $P_c( {\omega})$ (resp.
$P_d( {\omega})$), we denote the projection on
$L_c^2(\mathcal{H}_\omega)$  (resp. $L_d^2(\mathcal{H}_\omega)$)
associated to the above direct sum.  The spectral decomposition
of  a vector $X$ with respect to \eqref{eq:spectraldecomp} is
\begin{equation}\label{eq:vectdec}\begin{aligned} &
 X=\big [ P_{N_g(\mathcal{H} _\omega)}
 +   \sum  _{j=1}^m  (P _{
 \ker (\mathcal{H} _\omega-\lambda _j)}
   +
 P _{
 \ker (\mathcal{H} _\omega +\lambda _j)}  )
  +
  P _c (\mathcal{H}_\omega  )
   \big ]X =\\&
   \frac{\langle   X |
   \sigma _3\partial _\omega \Phi \rangle }
  {q'(\omega )}
     \sigma _3\Phi +
  \frac{\langle   X|
   \Phi \rangle }
  {q'(\omega )}
    \partial _\omega \Phi +\sum _{a=1}^3 \frac{\langle X |   x_a   \Phi \rangle }{q (\omega )} \partial _{x_a}\Phi -\sum _{a=1}^3 \frac{\langle X |\sigma _3\partial _{x_a}   \Phi \rangle }{q (\omega )} \sigma _3 x_a  \Phi
  \\&  +
   \sum _{j=1}^m \langle X |
    \sigma _3 \xi _j \rangle
   \xi _j  +\sum _{j=1}^m     \langle X |
   \sigma _1\sigma _3 \xi _j \rangle
   \sigma _1 \xi _j   +
   P_c(\mathcal{H}_\omega   )X.
 \end{aligned}
\end{equation}

The following lemma is well known.
\begin{lemma}
  \label{lem:modulation} Fix $U_o=e^{\im \sigma _3 (\frac{v_{o}\cdot (x-D_{o})}{2} +\vartheta _{o}) }   \Phi   _{\omega _o} (x-D_{o})$. Then $\exists$ a neighborhood $\U ^{1,0}$  of $U_o$
  in $H^1$    and  functions $\omega \in C^\infty (\U^{1,0} , \mathcal{O})$,
  $\vartheta \in C^\infty (\U ^{1,0} , \R )$ and $D,v \in C^\infty (\U ^{1,0} , \R  ^3)$,
  s.t. in $U_o$ their value is $(\omega  _{o}, \vartheta _{o},D_{o}, v_{o})$
  and s.t.  $\forall U\in \U ^{1,0}$\begin{equation}\label{eq:anzatz}\begin{aligned} &
 U(x)=  e^{\im \sigma _3(\frac{v\cdot (x-D)}{2} +\vartheta ) }
  (\Phi _\omega (x-D) + R (x-D))
 \text{  and $R\in N^{\perp}_g (\mathcal{H}_\omega ^*)$.}
\end{aligned}\end{equation}
\end{lemma}

Notice that, once the functions are give, we have
\begin{equation}\label{eq:anzatz2}\begin{aligned} &R(x)= e^{-\im \sigma _3 (\frac{v \cdot x}{2} +\vartheta  ) }  U (x+D)-\Phi _ \omega (x )
\end{aligned}  \end{equation}
    with the rhs just continuous in $U$. We can further decompose $R$ using
 \eqref{eq:spectraldecomp} as
\begin{equation}
  \label{eq:decomp2}
  R (x) =\sum _{j=1}^{m}z_j \xi _j(\omega ,x )+
\sum _{j=1}^{m}\overline{z }_j\sigma_1\xi _j(\omega ,x )
+P_c(\mathcal{H}_{\omega  } )f (x), \quad   f \in L_c^2(\mathcal{H}_{\omega _0})
\end{equation}
   where we fixed  $\omega _0\in
\mathcal{O}$  such that $q(\omega _0)=\| u_0\| _2^2$.
So we have
\begin{equation}\label{eq:coordinate}\begin{aligned} & U(x)= e^{\im \sigma _3(\frac{v\cdot (x-D)}{2} +\vartheta ) }
 \big ( \Phi
_\omega (x-D)+z \cdot \xi  (\omega , x-D ) \\& +  \overline{z }\cdot
\sigma_1\xi
 (\omega , x-D)+(P_c(\mathcal{H}_{\omega}) f)(x-D)\big ).
  \end{aligned}  \end{equation}
\eqref{eq:coordinate} is a system  of coordinates because   for
$\mathcal{O}$ sufficiently small the map $P_c(\mathcal{H}_{\omega})$
is an isomorphism from $L^2_c(\mathcal{H}_{\omega _0})$ to
$L^2_c(\mathcal{H}_{\omega  })$. Notice that the maps $U\to z_j$  are smooth.

We   set $z\cdot \xi =\sum _jz_j \xi _j$ and $\overline{z}\cdot
\sigma_1\xi =\sum _j\overline{z}_j \sigma_1\xi _j$. In the sequel we
set
\begin{equation} \label{eq:partialR} \partial _\omega R=
\sum _{j=1}^{m}z_j \partial _\omega \xi _j(\omega )+ \sum
_{j=1}^{m}\overline{z }_j\sigma_1\partial _\omega \xi _j(\omega
)+\partial _\omega P_c(\mathcal{H}_{\omega } )f.
\end{equation}
Sometimes we will denote $P_c(\omega )=P_c(\mathcal{H}_{\omega } )$.

We have the following formulas:
\begin{equation}\label{eq:vectorfields} \begin{aligned} &
\frac \partial {\partial  {\omega}}   =
   e^{ \im \sigma _3\Theta } \partial _\omega ( \Phi (x-D) +R(x-D))
\, ,\quad  \frac \partial {\partial  {\vartheta}} =\im
    \sigma _3 U (x) ,\\& \frac \partial {\partial   D_a}=- \frac \partial {\partial  x_a}U (x) \, ,\quad   \frac \partial {\partial  v_a}=\frac{\im }{2}\sigma _3  (x_a-D_a)  U (x) \text{ for $a=1,2,3,$}
,\\& \frac \partial {\partial  {z_j}}  =
   e^{ \im \sigma _3\Theta }  \xi _j (x-D)  \, ,\quad
   \frac \partial {\partial  {\overline{z}_j}}   =
   e^{ \im \sigma _3\Theta}\sigma _1 \xi _j(x-D) \text{ for $j=1,...,m$} .\end{aligned}
\end{equation}

 Lemmas \ref{lem:gradient omega}--\ref{lem:gradient f} are similar to analogous ones in \cite{Cu1}.

 \begin{lemma} \label{lem:gradient omega} There is a matrix $   \mathbb{A}  $ such that

\begin{equation} \label{eq:gradient omega}\begin{aligned} &  \begin{pmatrix} \langle \tau _{D}  e^{- \im \sigma _3\Theta}  \Phi  |
\\ \langle \tau _{D} e^{- \im \sigma _3\Theta }\sigma _3\partial _\omega \Phi  | \\  2 \langle \tau _{D} e^{- \im \sigma _3\Theta }  \sigma _3\nabla  _{x } \Phi  |   \\  \langle \tau _{D} e^{- \im \sigma _3\Theta }x  \Phi  |  \end{pmatrix}  = \mathbb{A}
\begin{pmatrix} -q'(\omega )d \omega
\\ -\im  q'(\omega )d \vartheta \\ \im q(\omega )d v  \\  -q(\omega )d D \end{pmatrix}  .\end{aligned}
\end{equation}
We have $\mathbb{A} =1+ \mathbb{A} _1 (z,f,\omega ,v) $ with
$\mathbb{A} _1= (1+|v|)O( |z|+\| f\| _{H^{-K,-S}} )$ smooth in the arguments $z\in \C ^m$, $f\in H^{-K,-S}$ and
$(\omega ,v)$,  for any pair $(K,S)$,.
\end{lemma}
\proof  Consider the functions of variable $(U,\omega , \vartheta ,v,D)$

\begin{equation} \label{eq:nals}\begin{aligned} & \mathcal{F} :=\langle
e^{-\im \sigma _3 (\frac{v \cdot ( x-D)}{2} +\vartheta  ) }  U (x )-\Phi _ \omega (x-D ) | \Phi _\omega (x-D)\rangle
   \\& \mathcal{G} :=\langle e^{-\im \sigma _3 (\frac{v \cdot ( x-D)}{2} +\vartheta  ) }  U (x )-\Phi _ \omega (x-D )  |\sigma _3\partial _\omega    \Phi _\omega (x-D)
 \rangle   \\&    \mathcal{B}_a :=\langle e^{-\im \sigma _3 (\frac{v \cdot ( x-D)}{2} +\vartheta  ) }  U (x )-\Phi _ \omega (x-D )  |(x_a -D ) \Phi _\omega (x-D) \rangle \\& \mathcal{D}_a := \langle e^{-\im \sigma _3 (\frac{v \cdot ( x-D)}{2} +\vartheta  ) }  U (x )-\Phi _ \omega (x-D )  |\sigma _3 \partial _a  \Phi _\omega (x-D)
 \rangle      . \end{aligned}
\end{equation}
Notice that  differentiating in $U$  we obtain the vectors in the lhs of
\eqref{eq:gradient omega}, which therefore span a vector bundle which has as sections the gradients $\nabla \omega $, $\nabla \vartheta $, $\nabla v_a $,
$\nabla D_a $.  For $R(x)$ defined by \eqref{eq:anzatz2}, we have following partial derivatives:
\begin{equation} \label{eq:nals1}\begin{aligned} & \mathcal{G} _\omega  =\langle R |
\sigma _3 \partial _\omega ^2 \Phi \rangle   \, , \quad \mathcal{G} _\vartheta  = \im ( q'(\omega ) +\langle R |
\partial _\omega \Phi \rangle   )
  \, ,\\&
\mathcal{G} _{D_a} =\frac{\im v_aq'(\omega )}{2} +
  \langle   (\partial _ a +\frac{\im}{2} \sigma _3v_a   ) R  |
 \sigma _3\partial _\omega \Phi \rangle
\, ,\quad   \mathcal{G} _{v_a}  =-\frac{\im}{2} \langle x_aR  |
\partial _\omega \Phi \rangle \, .
 \end{aligned}
\end{equation}
Similar formulas are satisfied by the other functionals in \eqref{eq:nals}.
Substituting   decomposition \eqref{eq:decomp2},
 we see that   the functions in  \eqref{eq:nals1} satisfy the regularity required for $\mathbb{A} _1$. For the other functions in \eqref{eq:nals} it is straightforward to check that the same is true. This yields by an elementary argument Lemma \ref{lem:gradient omega}.
\qed

  \begin{lemma} \label{lem:gradient z} For $\U ^{1,0}$
   in Lemma \ref{lem:modulation}  sufficiently small,
    $z_j\in C^\infty (\U ^{1,0},\C )$. The following
formulas hold, summing over repeated index $a$:
\begin{equation}\label{ZOmegaTheta}\begin{aligned}
\nabla z_j =&  - \langle \sigma _3 \xi _j| \partial _\omega R \rangle
  \nabla \omega - \im \langle \sigma _3 \xi _j|  \sigma _3 R \rangle
  \nabla \vartheta   -\frac{\im}{2} \langle   \sigma _3 \xi _j| \sigma _3 x_a R \rangle \nabla v_a \\&+
  \langle  \sigma _3  \xi _j|  ( \partial  _{x _a}+\im \sigma _3 \frac{v_a}{2} ) R \rangle
  \nabla D_a
  + \tau _{D} (e^{-\im \sigma _3\Theta }
  \sigma
_3\xi _j) (x )\\    \nabla \overline{z}_j =&  -
\langle \sigma _1\sigma _3 \xi _j|
\partial _\omega R \rangle \nabla \omega-
   \im \langle \sigma _1\sigma _3 \xi _j|  \sigma _3 R \rangle
   \nabla \vartheta  -\frac{\im}{2} \langle \sigma _1\sigma _3  \xi _j| \sigma _3 x_a R \rangle \nabla v_a \\&+
  \langle  \sigma _1\sigma _3  \xi _j|  ( \partial  _{x _a}+\im \sigma _3 \frac{v_a}{2} ) R \rangle
  \nabla D_a + \tau _{D} (e^{-\im \sigma _3\Theta }
  \sigma _1 \sigma
_3\xi _j) (x ) .\end{aligned}\nonumber
\end{equation}
\end{lemma}
\proof   The fact that $z_j\in C^\infty (\U ^{1,0},\C )$ follows from formula \begin{equation}z_j (U)=\langle U(x), e^{-\im \sigma _3(\frac{v\cdot (x-D)}{2} +\vartheta ) } \sigma _3 \xi _j (\omega , x-D  )\rangle ,\nonumber
\end{equation}
the fact that $\omega ,\vartheta , v,D\in C^\infty (\U ^{1,0}  )$ and the properties
of  $\xi _j (\omega , x   )$.
We have
\begin{equation}\label{eq:identitiesGradZ} \begin{aligned} & \langle \nabla z_j|
\tau _{D}e^{ \im \sigma _3\Theta } \xi _\ell \rangle =\delta _{j\ell}, \quad
\langle \nabla z_j|  \tau _{D}e^{ \im \sigma _3\Theta } \sigma _1\xi _\ell
\rangle  = 0 =\langle \nabla z_j|  \tau _{D}e^{ \im \sigma _3\Theta }
\sigma _3(\Phi +R) \rangle \\& \langle \nabla z_j|  \tau _{D}e^{ \im \sigma
_3\Theta } \partial _\omega (\Phi +R) \rangle =0=  \langle \nabla z_j,\tau _{D} e^{ \sigma
_3\Theta }  x_a  (\Phi +R) \rangle \\&  \langle \nabla z_j| \tau _{D} e^{ \sigma
_3\Theta } \left (\partial_{x_a}+ \im \sigma _3\frac{v_a}{2}  \right )  (\Phi +R) \rangle
=0\\&   \langle
\nabla z_j|   \tau _{D} e^{ \im \sigma _3\Theta } P _c( \omega ) P _c( \omega
_0)g \rangle =0 \quad  \forall g\in L^2_c(\mathcal{H}_{\omega _0}).
\end{aligned}
\end{equation}
The latter implies $ P _c( \mathcal{H}_{\omega
_0}^{*} ) P _c( \mathcal{H}_{\omega  }^{*} )e^{ \im \sigma
_3\Theta  }\tau _{-D}\nabla z_j=0$ and  $   P _c(
\mathcal{H}_{\omega  }^{*} )e^{ \im \sigma _3\Theta }\tau _{-D}\nabla
z_j=0$.
By \eqref{eq:vectdec} and for $\mathbf{A}$ and unknown vector,

\begin{equation} \begin{aligned} & \nabla z_j=  {\mathbf{A}}^T \begin{pmatrix} \tau _{D}((e^{- \im \sigma _3\Theta}  \Phi )(x )
\\ \tau _{D}(e^{- \im \sigma _3\Theta }\sigma _3\partial _\omega \Phi ) (x ) \\ \tau _{D}(e^{- \im \sigma _3\Theta }x  \Phi )(x   ) \\   \tau _{D}(e^{- \im \sigma _3\Theta }  \sigma _3\nabla  _{x } \Phi ) (x   ) \end{pmatrix} +  \tau _{D} (e^{-\im \sigma _3\Theta }
  \sigma
_3\xi _j) (x )  .\end{aligned}\nonumber
\end{equation}
 Using Lemma \ref{lem:gradient omega} there is a vector $\mathbf{A}_1$ s.t.
\begin{equation} \begin{aligned} & \nabla z_j=  \mathbf{A}_1 ^T \begin{pmatrix} \nabla \omega
\\ \nabla \vartheta \\ \nabla v  \\  \nabla D \end{pmatrix} +    \tau _{D} (e^{-\im \sigma _3\Theta }
  \sigma
_3\xi _j) (x )   .\end{aligned}\nonumber
\end{equation}
By
  \eqref{eq:vectorfields} and \eqref{eq:identitiesGradZ} we obtain
\begin{equation}   \mathbf{A}_1=-
\begin{pmatrix} \langle \sigma _3 \xi _j|  \partial _\omega R \rangle
  \\ \im \langle \sigma _3 \xi _j|  \sigma _3 R \rangle \\  \frac{\im}{2} \langle   \sigma _3 \xi _j|  \sigma _3x R \rangle \\  -\langle  \sigma _3  \xi _j|  (\nabla _x+\im \sigma _3\frac{v}{2} ) R \rangle  \end{pmatrix}.\nonumber
\end{equation}
Similar formulas hold for $\nabla \overline{z}_j$ yielding Lemma \ref{lem:gradient z}.
 \qed

\begin{lemma} \label{lem:gradient f} The map
$U\to f(U)=f$, for $U$ and $f$ as in \eqref{eq:coordinate}, is continuous from  $X$
into itself, for $X= L^2, H^1 $ and $\Sigma _n$   for any $n$. Furthermore we have
$f\in C^1(\U ^{1,0}, L^{2,-1})$,  $f\in C^1(\U ^{1,0}\cap \Sigma , L^{2 })$
and $f\in C^1(\U ^{1,0}\cap \Sigma _n, \Sigma _{n-1})$    with Frech\'et derivative
$f'(U) $   defined by the following formula, summing on the repeated index $a$,
\begin{equation}\label{eq:gradient f} \begin{aligned}
f'(&U)=   (P_c(\omega )P_c(\omega _0))^{-1}  P_c(\omega )\big [
 -   \partial _\omega R \, d\omega   \\& -\im
  \sigma _3 R  \, d\vartheta - \frac{\im}{2}  \sigma _3 x_a R \, d v_a +( \partial  _{x _a}+\im \sigma _3\frac{v_a}{2} ) R \, d  D_a  +  e^{-\im
\sigma _3\Theta }\tau _{-D}\uno \big  ] .
\end{aligned}
\end{equation}
    We further have       $f\in C^n(\U ^{1,0}\cap \Sigma , H^{-n+1,-n})$

               \end{lemma}
\proof    Continuity follows from Lemmas \ref{lem:modulation}-\ref{lem:gradient z}
and formula \eqref{eq:coordinate} solved w.r.t. $f$. The latter proves also the
$C^1$   as well the $C^n$  properties.
The proof of formula \eqref{eq:gradient f}   is similar to that of Lemma \ref{lem:gradient z}, see also Lemma 4.2 \cite{Cu1}.
\qed

In the sequel given a scalar function $\psi (U)$ which is differentiable,
we will denote by $\nabla _f\psi (U)$ the only element in $L^2_c(\mathcal{H}^*_{\omega _0})$ s.t. for any $g\in L^2_c(\mathcal{H} _{\omega _0})$ we have $\langle \nabla _f\psi (U)|g\rangle =\langle \nabla \psi (U)|P_c^* (\omega )g\rangle   $.

\section{Symplectic structure}
\label{sect:symplectic}

 Our ambient space is    $ \mathbf{X} \times \mathbf{X}$ where we can have $\mathbf{X}=L^2,H^1,\Sigma _n$.
We focus only on points with $\sigma _1U=\overline{U}$, so that the space
is identified with $\mathbf{X}$. The natural
symplectic structure for our problem is
\begin{equation}\label{eq:SymplecticForm}
  \Omega (X,Y)=\langle X|  \sigma _3\sigma _1 Y \rangle .
\end{equation}

\begin{definition}\label{def:HamField} Let $F\in C^{1}( \mathcal{U},\C )$ for $ \mathcal{U}$ and open subset of  $\mathbf{X} $.     Then
the Hamiltonian vectorfield of  $F$ with respect to a symplectic form $\Omega $ is  the field $X_F $ such that $ \Omega (X_F ,Y)=-\im dF (Y)$ for any given
tangent vector  $Y\in T\mathcal{U}.$ More explicitly, $X_F=-\im \sigma _3\sigma _1 \nabla F $ for the form in \eqref{eq:SymplecticForm}.
\end{definition}

\begin{definition}\label{def:PoissonBrac}
For    $F,G \in C^{1}( \mathcal{U},\C )$  as above,  we call Poisson bracket of a pair  of   $F$
and $G$  the   function
\begin{equation}\label{eq:PoissonBracket}
  \{ F,G \}  =  dF (X_G )
= \im \Omega ( X_F ,   X_G  ) =-\im \langle \nabla F | \sigma _3\sigma _1\nabla  G   \rangle
   .
\end{equation}
Let $\mathcal{G} \in C^{1}( \mathcal{U},\mathbb{E} )$ with $\mathbb{E}$ a given
Banach space on $\C$. Then, for $F\in C^{1}( \mathcal{U},\C )$ we set, for $\mathcal{G}'$ the Frechet derivative of $\mathcal{G} $,
 \begin{equation}\label{fF}   \{ \mathcal{G}, F\} :=
\mathcal{G}'(U)X_F (U)=-\im
\mathcal{G}'(U) \sigma _3 \sigma _1 \nabla F (U).
\end{equation}
  We set
  $\{ F,\mathcal{G} \}  :=-\{ \mathcal{G}, F\} $.
\end{definition}
 Obviously our system is hamiltonian. It is important to cast it in terms
 of the Poisson brackets.
\begin{lemma}\label{lem:HamFor} In the coordinate system \eqref{eq:coordinate}, system \eqref{eq:NLSvectorial} can be written, for   $F=E$ as  $ \dot \vartheta = \{ \vartheta , F  \} $ and
\begin{equation} \label{eq:SystPoiss} \begin{aligned} &
  \dot \omega  = \{ \omega , F \} \, ,  \quad    \dot f= \{f, F
\}   \, , \\& \dot D_a  = \{ D_a , F \} \, , \quad    {\dot
 v_a }=  \{ v_a , F \}  \text{ for $a=1,2,3$} \, ,
\\&  \dot z_j  = \{ z_j , F \} \, , \quad    {\dot
{\overline{z}}_j }=  \{ \overline{z}_j , F \}
 \text{ for $j=1,...,m$}    \,   . \end{aligned}
\end{equation}
\end{lemma}
\proof The statement is not standard only for   $\dot f= \{f, E
\} $.  Notice that \eqref{eq:NLSvectorial} can be written as
\begin{equation}\label{eq:time derivative}\begin{aligned} &
\im \dot  U= -\sigma _3
\left (  \dot \vartheta
  + \frac{\dot v \cdot (x-D)}{2}\right ) U + \im \dot \omega e^{\im \sigma _3\Theta }
     \partial _\omega (\Phi   +R)  \\& - \im  \dot D\cdot  \nabla (^{\im
\sigma _3\Theta } (\Phi +R))    +  \im    e^{\im
\sigma _3\Theta }
  (
  \dot z\cdot \xi    +
   \dot {\overline{z}}\cdot \sigma _1\xi +
    P_c(\mathcal{ H}_\omega ) \dot f) .\end{aligned}
\end{equation}
So,  by
\eqref{eq:vectorfields},  system $ \im \dot  U= \sigma _3
    \sigma _1
  \nabla E  (  U)  $ is the same as
  \begin{equation}
\label{eq:sys1}
\begin{aligned} &
   \im
 \dot \vartheta
  \frac{\partial}{\partial \vartheta } +\sum _{a=1} ^{3} \dot v_a \frac{\partial}{\partial v_a }  +  \im \dot \omega
 \frac{\partial}{\partial \omega }  +\sum _{a=1} ^{3} \dot D_a \frac{\partial}{\partial D_a }  +  \im
  \sum _{j=1}^{m} \dot z_j\frac{\partial}{\partial z _j } +  \im
  \sum _{j=1}^{m}  \dot {\overline{z}}_j\frac{\partial}{\partial
  \overline{z} _j } \\& +\im
   e^{\im \Theta\sigma
_3} P_c(\mathcal{ H}_\omega ) \dot f  =  \sigma _3
    \sigma _1
  \nabla E  (  U)    .\end{aligned}
\end{equation}
When we apply the  derivative $f'(U)$ to \eqref{eq:sys1}  the first
line cancels, so we get
\begin{equation}
\begin{aligned} &  f'(U)
   e^{\im \Theta\sigma
_3} P_c(\mathcal{ H}_\omega ) \dot f  =- f'(U)\im \sigma _3
    \sigma _1
  \nabla E  (  U) =f'(U)X_E(U)=\{ f,E\} ,
\end{aligned}\nonumber
\end{equation} from
the definition of hamiltonian field  and of Poisson bracket. Notice now that
$\dot f= f'(U)
   e^{\im \Theta\sigma
_3} P_c(\mathcal{ H}_\omega ) \dot f.$ This follows from the fact that
\begin{equation}
\begin{aligned} &  f'(U)
   e^{\im \Theta\sigma
_3} P_c(\mathcal{ H}_\omega ) \dot f = \frac{d}{ds}_{|_{s=0}}f( U (\omega , \vartheta ,D,v, z , \overline{z} , f+s\dot f) )\\&  =\frac{d}{ds}_{|_{s=0}} (f+s\dot f) = \dot f   .
\end{aligned}\nonumber
\end{equation}
Hence $\dot f=\{ f,E\} $.
\qed

\section{ Reduction of variables}
\label{sec:Reduction}

The following formulas are important.
\begin{lemma}
  \label{lem:Involutions} Consider charge  $Q$ and momenta $\Pi _a$, see  defined in
  \eqref{eq:charge}. Then
\begin{equation} \label{eq:Ham.VecField} X_Q=-
\frac{\partial}{\partial \vartheta}\, , \quad X  _{\Pi _a}=-
\frac{\partial}{\partial D_a} .
\end{equation}
In particular
  \begin{equation} \label{eq:Ham.VecField1} \{ Q,\vartheta  \} =\{ \Pi _a,D_a  \}  =1
\end{equation}
  and any  Poisson bracket not in \eqref{eq:Ham.VecField1}   of any of the invariants  $Q,\Pi _a$ with any of the coordinates from \eqref{eq:coordinate} is equal to 0.
  The functions  $ Q $, $  \Pi _a     $ and $E$ Poisson commute.
\end{lemma}
\proof  \eqref{eq:Ham.VecField} follows from    \eqref{eq:vectorfields}:

  \begin{equation}
\begin{aligned} &  X_Q(U) = -\im \sigma _3 \sigma _1 \nabla Q(U)  =  -\im \sigma _3 U =-\frac{\partial}{\partial \vartheta} \, , \\&  X  _{\Pi _a} (U) = -\im \sigma _3 \sigma _1 \nabla \Pi _a (U)  =   \im \sigma _3 \sigma _1 \sigma _2\partial _a U = \partial _a U=- \frac{\partial}{\partial D_a}  . \end{aligned}\nonumber
\end{equation}
The second part of the statement follows immediately from  \eqref{eq:Ham.VecField}.
\qed

The following lemma  is elementary.
\begin{lemma}
  \label{lem:gauge} For $ U=  e^{\im \sigma _3 (\frac{v\cdot (x-D)}{2} +\vartheta ) }  ( \Phi _ \omega (x-D) + R(x-D))$ and summing on repeated indexes, we have the following formulas:
  \begin{equation} \begin{aligned} &    Q(U)= Q(\Phi _ \omega   + R )= q(\omega ) + Q(R)\, ; \\&
   \Pi _a(U)= \Pi _a(\Phi _ \omega   + R) -  {2}^{-1}v_aQ(\Phi _ \omega   + R ) \, ;
\\& E(   U)= E(   \Phi _ \omega   + R) - v_a \Pi _a(\Phi _ \omega   + R)
 + \frac{v^2 }{4}Q(\Phi _ \omega   + R ).\end{aligned}\nonumber
\end{equation}
We have $\Pi _a(\Phi _ \omega   + R) = \Pi _a(  R)$.

 \noindent  $Q(U)$ and $\Pi _a (U)$ are $C^\infty $ and
$E(U)$ is $C^1 $
in  $\omega$, $v$, $z$ and $f\in H^1$ . Furthermore both $E$ and $\nabla E$ depend
smoothly on $(\omega , z)$.
  \end{lemma}

We introduce now a new hamiltonian, for $U_0^T=(u_0, \overline{u}_0 )$,
 \begin{equation}     \label{eq:K}  \begin{aligned} &
 {K}(U):=E(U)- E\left (   \Phi _{\omega_0}\right ) +  \omega (U)  \left ( Q(U)- Q(U_0)\right ) .
\end{aligned}
\end{equation}
Lemma \ref{lem:Involutions}  implies that   the solutions of \eqref{eq:NLSvectorial} with   charge  $ q(\omega _0)$
satisfy  \eqref{eq:SystPoiss} with $F=K$ and   $\dot \vartheta -\omega = \{ \vartheta ,  {K}  \} $. We would like for $K$  to be our hamiltonian, but obviously $K$ is not a hamiltonian for \eqref{eq:NLSvectorial}.   To obviate this  we notice that  $\partial _{D_a} {K}= \partial _{\vartheta} {K}\equiv 0$ imply that the evolution    of the variables $\omega, v_a, z_j, \overline{z}_j, f$ is unchanged if we consider the following new
hamiltonian system:
\begin{equation} \label{eq:SystK} \begin{aligned} &
\dot \omega  = \{ \omega ,  {K} \} \, , \quad   \dot \vartheta   = \{ \vartheta , K  \}  \, , \quad    \dot f= \{f,  {K}
\}   \, , \\& \dot D_a   = \{ D_a , {K} \} \, , \quad    {\dot
 v_a }=  \{ v_a ,  {K} \}  \text{ for $a=1,2,3$} \, ,
\\&  \dot z_j  = \{ z_j ,  {K} \} \, , \quad    {\dot
{\overline{z}}_j }=  \{ \overline{z}_j ,  {K} \}
 \text{ for $j=1,...,m$}    \,  .  \end{aligned}
\end{equation}
It is elementary that for  \eqref{eq:SystK}   the charge $Q(U)$ and the momenta $\Pi _a(U)$ are invariants of motion.

\bigskip

We   proceed  to  a  reduction of order in \eqref{eq:SystK} such as
  as described for instance in   Theorem 6.35 p.402 \cite{Olver}.
   Here the discussion is elementary because  we have no need to prove the existence of particular variables, see p.395 \cite{Olver}.
We set  \begin{equation} \label{eq:variables} \begin{aligned} &
   Q:= Q(U)=q(\omega )+ Q(R)\quad , \\&   \Pi_a :=\Pi _a(U)=\Pi _a(R) -\frac{v_a  }{2} \left ( q(\omega )+ Q(R)\right ). \end{aligned}
\end{equation}
\begin{lemma}
\label{lem:var} After we express  $R$  in \eqref{eq:variables} in terms of $\omega$, $z$ and $f$, see \eqref{eq:decomp2},
there is an implicit function $\omega =\omega ( Q, z,f)$
  defined by the first of \eqref{eq:variables},
  with $\omega ( Q, z,f)$
smooth in $Q$, $z $ and  in $f\in L^2_c (\omega _0)$.  Similarly,
$v_a =v _a(\Pi _a ,Q, z,f)$, with the latter smooth in  $\Pi _a$, $Q$, $z $ and $f\in H^\frac{1}{2}_c (\omega _0)$.
\end{lemma}\proof For $\omega =\omega ( Q, z,f)$ the statement follows from the implicit function theorem . Write $v_a = 2 Q^{-1}(\Pi _a(R) -\Pi _a)$    and substitute  $\omega =\omega ( Q, z,f)$ in the definition of $R$.   \qed

 \bigskip Lemma \ref{lem:var}  allows to move from the variables
in the rhs of \eqref{eq:coordinate} to a new system of variables obtained replacing the functions $(\omega , v)$ with the $(Q , \Pi )$.
    \begin{lemma}
  \label{lem:reduction1} The vectorfields $\frac{\partial}{\partial \vartheta}$ and $ \frac{\partial}{\partial D_a}$ are the same for the two systems of coordinates.
  In particular, in the new system of coordinates we continue to have
   $
    \frac{\partial}{\partial \vartheta} K= \frac{\partial}{\partial D_a} K=0  .
$
 \end{lemma}
 \proof It is an immediate consequence of
$
    \frac{\partial}{\partial \vartheta} Q(U)= \frac{\partial}{\partial D_a} Q(U)=0 $ and of $
    \frac{\partial}{\partial \vartheta} \Pi _b(U)= \frac{\partial}{\partial D_a} \Pi _b(U)=0 $  in the old coordinate system, and of the chain rule.

 \qed

 In the new variables, system \eqref{eq:SystK} reduces to the pair
 of systems

        \begin{equation} \label{eq:SystK1} \begin{aligned} &
\dot Q = 0 \, , \quad   \dot \vartheta   = \{ \vartheta , K  \}  \, , \quad    \dot D_a   = \{ D_a ,K \} \, , \quad    {\dot
 \Pi _a }=  0 \quad  \text{ for $a=1,2,3$} \, ,
   \end{aligned}
\end{equation}
and
\begin{equation} \label{eq:SystK21} \begin{aligned} &
      \dot f= \{f, K
\}   \, ,\quad     \dot z_j  = \{ z_j , K \} \, , \quad    {\dot
{\overline{z}}_j }=  \{ \overline{z}_j , K \}\quad
 \text{ for $j=1,...,m$}    \,  .  \end{aligned}
\end{equation}
 Now we restrict to the set with $Q=  Q(U_0) $ and $\Pi _a= 0$. Thanks to Lemma \ref{lem:reduction1}    system  \eqref{eq:SystK21}   is closed.

\bigskip

\begin{lemma}
  \label{lem:ExpK}  Consider the restriction of the variables $(Q , \Pi )$ at  the fixed values $Q=  Q(U_0) $ and $\Pi _a= 0 $
  and set $\varrho (f):=\left ( Q(f), \Pi (f) \right )$.
   Then we have  the
expansion
\begin{equation}  \label{eq:ExpK1} \begin{aligned} & K =  \underline{{\psi}} (\varrho (f)) +K_2  +\underline{\resto} ^{(1)}
\end{aligned}
\end{equation}    where  $\underline{{\psi}} (\varrho )$  is smooth in $\varrho $
and where the following holds.
\begin{itemize}
\item[(1)]
 We have

\begin{equation}  \label{eq:ExpK2} K_2  =
\sum _{\substack{ |\mu +\nu |=2\\
\lambda ^0\cdot (\mu -\nu )=0}}
 \underline{a}_{\mu \nu} ( \varrho (f) )  z^\mu
\overline{z}^\nu + \frac{1}{2} \langle \sigma _3 \mathcal{H}_{\omega
_0} f|  \sigma _1 f\rangle .
\end{equation}

\item[(2)] We have $\underline{\resto} ^{(1)}=\widetilde{\underline{\resto} ^{(1)}} +
 \widetilde{\underline{\resto} ^{(2)}}  $, with $\widetilde{\resto ^{(1)}}=$
\begin{equation}  \label{eq:ExpK2resto} \begin{aligned} &
 =\sum _{\substack{ |\mu +\nu |=2\\
\lambda ^0\cdot (\mu -\nu )\neq 0  }} \underline{a}_{\mu \nu } (\varrho (f)
)z^\mu \overline{z}^\nu +\sum _{|\mu +\nu |  = 1} z^\mu
\overline{z}^\nu \langle \sigma _1 \sigma _3\underline{G}_{\mu \nu }(\varrho (f)
)| f\rangle \end{aligned}
\end{equation}

    \begin{equation}    \begin{aligned} &  \widetilde{\underline{\resto} ^{(2)}}=   \sum _{|\mu +\nu |= 3} z^\mu \overline{z}^\nu  \underline{a}_{\mu \nu
}( z,\varrho (f) )     +\sum _{|\mu +\nu
|= 2} z^\mu \overline{z}^\nu \langle  \underline{G}_{\mu \nu }( z,\varrho (f)
)| \sigma _3\sigma _1 f\rangle    \\& +   \sum _{d=2}^4
\langle \underline{B}_{d } (   z ,\varrho (f) )|     f  ^{   d} \rangle
      +\int _{\mathbb{R}^3}
\underline{B}_5 (x,  z , f(x),\varrho (f) )  f^{   5}(x) dx+   E_P (  f).
\end{aligned}\nonumber
\end{equation}
with $\underline{B}_2(0,0)=0$ and where, both here and   in   Theorem \ref{th:main} later, by $f^d(x)$ we schematically represent $d-$products of components of $f$.

\item[(3)] At $\varrho (f)=0$
\begin{equation}  \label{eq:ExpKcoeff1} \begin{aligned} &
\underline{a}_{\mu \nu } ( 0 ) =0 \text{ for $|\mu +\nu | = 2$  with $(\mu
, \nu )\neq (\delta _j, \delta _j)$ for all $j$,} \\& \underline{a}_{\delta _j
\delta _j }^{(r)}( 0 ) =\lambda _j (\omega _0)  , \text{ where
$\delta _j=( \delta _{1j}, ..., \delta _{mj}),$}
\\& \underline{G}_{\mu \nu }(  0 ) =0 \text{ for $|\mu +\nu | = 1$ }
\end{aligned}
\end{equation}
These $\underline{a}_{\mu \nu } ( \varrho )$ and $\underline{G}_{\mu \nu }( x,\varrho
)$ are smooth in all variables with $G_{\mu \nu }( \cdot ,\varrho )
\in C^\infty ( \mathbb{R}^4, H^{K,S} _x(\mathbb{R}^3,\mathbb{C}^2))$
for all $(K,S)$.

\item[(4)] $\underline{a}_{\mu \nu
}(   z,  \varrho   ) \in C^\infty ( \mathrm{U},
 \mathbb{C} ) $ for   a
small neighborhood $\mathrm{U}$ of $( 0,0)$ in
$ \mathbb{C}^m\times \R ^4$.

\item[(5)] $\underline{G}_{\mu \nu
}( \cdot ,  z,  \varrho ) \in C^\infty ( \mathrm{U},
H^{K,S}_x(\mathbb{R}^3,\mathbb{C}^2)) $, for $\mathrm{U}$ like in
(4) and for all $(K,S)$.

\item[(6)] $\underline{B}_{d
}( \cdot ,  z,\varrho  ) \in C^\infty ( \mathrm{U},
H^{K,S}_x(\mathbb{R}^3, B   (
 (\mathbb{C}^2)^{\otimes d},\mathbb{C} ))) $, for $2\le d \le 4$
for $\mathrm{U}$ and   $(K,S)$  like above.

\item[(7)] Let $^t\eta = (\zeta , \overline{\zeta}) $ for
$ \zeta \in \C$. Then for
  $\underline{B}_5(\cdot ,\omega  , z , \eta ,\varrho )$   we have for all $(K,S)$
\begin{equation} \label{5power2}\begin{aligned} &\text{for any $l$ ,
}  \| \nabla _{  z,\overline{z} ,\zeta,\overline{\zeta} ,\varrho  }
^l\underline{B}_5(  z,\eta ,\varrho  ) \| _{H^{K,S}_x (\mathbb{R}^3,   B   (
 (\mathbb{C}^2)^{\otimes 5},\mathbb{C} )} \le C_l (K,S).
 \end{aligned}\nonumber \end{equation}

\item[(8)] For all indexes and  for $b_{\mu \nu }=\underline{a}_{  \mu \nu  }$ and $B_{\mu \nu }=\underline{G}_{\mu \nu } $, we have:
\begin{equation}  \label{eq:ExpKcoeff2} \begin{aligned} & b_{\mu \nu }  =
\overline{b_{\nu \mu   }} \, , \quad   B_{\mu \nu } =-\sigma
_1\overline{B_{\nu \mu }} .
\end{aligned}
\end{equation}
\item[(9)] $\underline{B}_5$ depends on $f(x)$ with terms $|f(x)|^2 $ and
$\varphi (x) f(x)$, with $\varphi (x) $ Schwartz functions.

 \end{itemize}
\end{lemma}
\proof
  By
 Lemma \ref{lem:gauge}  we get from \eqref{eq:K}, for $d(\omega)$  and  $q(\omega)$ see \eqref{eq:function q},
\begin{equation}    \label{eq:ExpK3}  \begin{aligned}
K(U)& =E(\Phi _\omega +R)+  \omega Q(\Phi _\omega +R)  -d(\omega _0)- (\omega -\omega _0)q( \omega _0) \\& + \frac{v^2}{4} Q(\Phi _\omega +R) -  v_a  \Pi _a(R) .
\end{aligned}
\end{equation}
 We have \begin{equation}    \label{eq:ExpK31}  \begin{aligned} &E(\Phi _\omega +R)+  \omega Q(\Phi _\omega +R)= d(\omega )+\frac{1}{2} \langle \sigma _3 \mathcal{H}_{\omega
 } R|  \sigma _1 R\rangle + \widehat{\resto} ^{(1)},
\end{aligned}
\end{equation}
with an expansion of the following form:

 \begin{equation}    \begin{aligned} &  \widehat{\resto} ^{(1)}=   \sum _{|\mu +\nu |= 3}  z^\mu \overline{z}^\nu  \underline{a}_{\mu \nu
}(\omega , z  )    +\sum _{|\mu +\nu
|= 2} z^\mu \overline{z}^\nu \langle  \underline{G}_{\mu \nu }(\omega ,z
)| \sigma _3\sigma _1P_c(\omega )f\rangle    \\& +   \sum _{d=2}^4
\langle \underline{B}_{d } ( \omega , z  )|  f ^{  d} \rangle
+   \int _{\mathbb{R}^3}
\underline{B}_5 (x,\omega, z , f(x) )  f^{  5}(x) dx +E_P(P_c(\omega )f).
\end{aligned}\nonumber
\end{equation}
The coefficients in the   expansion of $\widehat{\resto} ^{(1)} $   satisfy appropriate smoothness
and symmetry properties. When we restrict to
$Q=  Q(U_0) $, by $q(\omega )+Q(R)= q(\omega _0 )$
we get that $\omega -\omega _0=O(|z|^2+\| f \| _{2}^2+ \| f \| _{H^{-K,-S}}^2)$. Notice that we can express the function $\omega (Q,z,f)$ as a function $\omega (Q, z,f,\varrho _0(f))$, with $\varrho _0(f)=Q(f)$, such that
$\omega (Q, z,f,\varrho _0 )$ is smooth in the variables $z\in \C^m$, $f\in H^{-K,-S}$ and $Q,\varrho _0 \in \R$.   Hence $\widehat{\resto} ^{(1)}  $ is a sum of terms which can be absorbed in
$\widetilde{\resto ^{(2)}}$, with   dependence on $\varrho (f)=(Q(f),\Pi (f))$
reduced at that on $Q(f)$.
We have by $d' =q $ \begin{equation}    \label{eq:ExpK4}  \begin{aligned} d(\omega )  -d(\omega _0)- (\omega -\omega _0)q( \omega _0) = \frac{1}{2} d''(\omega _0)  (\omega - \omega_0)^2+  o(\omega - \omega_0)^2.\end{aligned}\end{equation}
     So \eqref{eq:ExpK4} can be absorbed in part in $\underline{\psi} (\varrho (f))$ and part in the rest  of expansion \eqref{eq:ExpK1}. We have

\begin{equation}    \label{eq:ExpK5}  \begin{aligned}
\frac{1}{2} \langle \sigma _3 \mathcal{H}_{\omega
 } R|  \sigma _1 R\rangle &=
\frac{1}{2} \langle \sigma _3 \mathcal{H}_{\omega
 _0} R|  \sigma _1 R\rangle \\& +(\omega -\omega _0) Q(R) + \frac{1}{2} \langle \sigma _3 (V_{\omega
 } -V_{\omega
 _0}) R|  \sigma _1 R\rangle  ,\end{aligned}
\end{equation}
 where  $\frac{1}{2} \langle \sigma _3 \mathcal{H}_{\omega
 _0} R| \sigma _1 R\rangle =\sum _j\lambda _j(\omega _0)|z_j|^2+
\frac{1}{2} \langle \sigma _3 \mathcal{H}_{\omega
 _0} f| \sigma _1 f\rangle $.  In the last line of \eqref{eq:ExpK5}, the first term is similar to  \eqref{eq:ExpK4} and the second can be absorbed in $H_2^{(1)}
+\resto ^{(1)}.$  Since   we have $  v_a  q(\omega _0) =2\Pi _a(R)  $,  we get that also   the second line of \eqref{eq:ExpK3} can be decomposed into terms with the properties of the various terms in the rhs of \eqref{eq:ExpK1}.

\qed

\section{A vectorfield needed for     Darboux Theorem}
\label{section:Darboux}

 We introduce the  2-form,   summing on repeated indexes,
\begin{equation} \label{eq:Omega0} \Omega _0=\im d\vartheta \wedge
dQ +  \im dD_a \wedge
d\Pi _a  +   dz_j\wedge d\overline{z}_j+\langle f'  | \sigma _3
\sigma _1 f'   \rangle .
\end{equation}
Through a change of variable we will transform the  $\Omega $ defined \eqref{eq:SymplecticForm} into the  $\Omega _0$. A brief reminder of J.Moser's
scheme of proof is in Sect.7 \cite{Cu1} and will not be repeated here.
The first step in the implementation  of this scheme is the search of an  appropriate vectorfield $ \mathcal{X} ^t$, which satisfies an equation stated in Lemma \ref{lem:vectorfield0}.  Then we will have $\phi ^*\Omega =\Omega _0$, with
$\phi =\phi ^1$ the Lie transform associated to $\phi ^t$, the flow of $ \mathcal{X} ^t$. Existence and differentiability of $\phi ^t$ are not obvious. Lemma \ref{lem:vectorfield} will tell us that  to get $\phi ^t$ we need
information on a quasilinear hyperbolic symmetric system mixed to ODE's for the
discrete modes. In Lemma \ref{lem:vectorfield} we will establish some properties of the coefficients of
the system.
Sections \ref{sec:ODE}--\ref{sec:Quasilinear} contain material needed to establish existence and differentiability of $\phi ^t$. Lemma \ref{lem:1forms}
is similar to Lemma 7.2 \cite{Cu1}.

\begin{lemma}
  \label{lem:1forms}
  At the points  $\tau _D e^{\im \sigma _3 \left (\frac{v \cdot x}{2}
  +\vartheta \right )} \Phi _{\omega_0}$ for all $(\vartheta ,D,v)\in \R ^7 $  we have $\Omega _0 =\Omega   .$

  \noindent Consider the following  forms:
\begin{equation} \label{eq:1forms1}\begin{aligned} &
 \mathrm{ B} _0 := \im  \vartheta dQ  +\im  D_a d\Pi _a -  \frac{\overline{z}_j dz_j -
{z}_j d\overline{z}_j}{2} +\frac{1}{2}\langle f  | \sigma _3\sigma
_1f '  \rangle
\end{aligned}
\end{equation}
  and $\mathrm{B}:=\mathrm{B_0}+\mathrm{\Gamma} $ for $\mathrm{\Gamma}={\varphi} dQ+(\mathrm{\Gamma})_jdz_j+ (\mathrm{\Gamma})_{\overline{ j}} d\overline{z}_j+\langle  (\mathrm{\Gamma})_f| f'\rangle $, where
  \begin{equation} \label{eq:1forms2}\begin{aligned} &
\varphi :=\frac{\langle \sigma _1\sigma _3R|
  \partial _\omega R \rangle + \im v_a\langle  x_a\sigma _1 R|
  \partial _\omega R \rangle  }{2(q'+\partial _\omega Q(R))} ,\\&
(\mathrm{ \Gamma }) _j:=    \frac{\im}{2} v_a\langle x_a \sigma _1 R| \xi _j \rangle -
   {\varphi} \partial _{ j} Q(R)  \, , \,  (\mathrm{ \Gamma }) _{\overline{ j}} :=\frac{\im}{2} v _a
   \langle  x_a   R| \xi _j \rangle -
   \varphi  \partial _{\overline{ j}} Q(R),  \\& (\mathrm{ \Gamma }) _f:= P^* _c(\omega _0)\left ( \frac{\im}{2} v_a  P^* _c(\omega )\sigma _1x_a   R-  \frac{1}{2}
   \sigma _1\sigma _3  P  _d(\omega )  f
   -
   {\varphi}  P^* _c(\omega )\sigma _1 R\right ).
\end{aligned}
\end{equation}
   Then
  $ d\mathrm{B_0}=\Omega _0$ and $
  d\mathrm{B } =\Omega   .$

\end{lemma}
\proof     $ d\mathrm{B_0}=\Omega _0$ follows from the definition of exterior differential. Since  $| {\varphi} |
\le C (|z|+\| f \| _{H^{-K,-S}})^2 $ and since $P  _d(\omega )  f=
(P  _d(\omega )- P  _d(\omega _0)) f$ is  a 0 of order 2 at $\omega =\omega _0$ and $f=0$,
to get $d\mathrm{\Gamma}=0$ at $R=0$ and $\omega =\omega _0$ it is enough to show  \begin{equation}d\left ( \langle x_a \sigma _1 R| \xi _j \rangle   dz_j+
   \langle  x_a   R| \xi _j \rangle  d\overline{z}_j + \langle  P^* _c(\omega )\sigma _1x_a   R  | f'    \rangle   \right ) =0 \text{ at $R=0$}.\nonumber\end{equation}
But this is true. Formally   $d  \left ( d \langle x_a \sigma _1 R|R
\rangle - \partial _\omega \langle x_a \sigma _1 R|R
\rangle  d\omega \right )=0$ at $R=0$.

  \noindent We   prove the formula for $\mathrm{B}$. Set $\widetilde{B} :=\frac{1}{2}\langle \sigma _1\sigma _3U |$. Notice that $d \widetilde{B}=\Omega $. We get

\begin{equation} \label{eq:beta1} \begin{aligned} &
\widetilde{B} =\frac{1}{2}\langle \tau _De^{-\im \sigma _3\Theta}\sigma _1\sigma
_3 \Phi  |  + \frac{1}{2}\langle \tau _De^{-\im \sigma _3\Theta}\sigma _1\sigma _3P_c(\omega )f|  \\& +
\frac{1}{2}   z_j \langle \tau _De^{-\im \sigma _3\Theta}\sigma
_1\sigma _3\xi _j| -\frac{1}{2}\overline{z}_j \langle \tau _De^{-\im \sigma _3\Theta}\sigma _3\xi _j | ,
\end{aligned}
\end{equation}
using \eqref{eq:coordinate}.
By Lemmas \ref{lem:gradient z} and \ref{lem:gradient f},  the  sum of the last three terms equals

\begin{equation} \label{eq:beta11} \begin{aligned} &
 \frac{1}{2} z_j d\overline{z}_j -\frac{1}{2} \overline{z}_j d {z}_j +\frac{1}{2}
 \langle \sigma _1\sigma _3 f|P_c(\omega )f'\rangle +  \frac{1}{2}
 \langle  \sigma _1\sigma _3 R| \partial _{\omega}R\rangle d\omega  \\&  -  {\im }
 Q(   R ) d\vartheta -\frac{\im }{4}
 \langle  \sigma _1  R|  x_a R \rangle dv_a -\im  \left ( \Pi _a(R)-\frac{v_a}{2}
 Q(   R ) \right ) dD_a  .
\end{aligned}
\end{equation}
Applying $\frac{1}{2}  \tau _De^{-\im \sigma _3\Theta}\sigma _1\sigma
_3 $   to decomposition \eqref{eq:vectdec}  for $X=\Phi$, we get
\begin{equation} \label{eq:beta2} \begin{aligned} &
 \frac{1}{2}\langle \tau _De^{-\im \sigma _3\Theta}\sigma _1\sigma
_3 \Phi | =  -\frac{q}{q'}  \langle \tau _De^{-\im \sigma _3\Theta}\sigma _3 \partial _\omega \Phi  |  - \frac{1}{2}\langle
\tau _De^{-\im \sigma _3\Theta}P_c (\mathcal{H}_\omega ^{*})\sigma _3
\Phi |
\\& -\frac{1}{2}   \langle \sigma _3\Phi |   \xi _j
\rangle    \langle \tau _De^{-\im \sigma _3\Theta} \sigma _3\xi
_j|  +\frac{1}{2}   \langle \sigma _3\Phi |   \xi _j
\rangle \langle \tau _De^{-\im \sigma _3\Theta}\sigma
_1\sigma _3\xi _j|   .
\end{aligned}
\end{equation}
By \eqref{eq:nals1} we have
\begin{equation} \label{eq:beta3}\begin{aligned} &-\frac{q}{q'}\langle \tau _De^{-\im \sigma _3\Theta}\sigma _3 \partial _\omega \Phi |
=\frac{q}{q'}\langle R| \sigma _3 \partial _\omega ^2 \Phi \rangle
\, d \omega -\im \, \frac{q}{q'}   \, (q'   +\langle R|   \partial
_\omega \Phi \rangle ) \, d \vartheta  \\& +  \frac{q}{q'}
\left (  \frac{\im v_a}{2}   q' +
\langle (\partial _a   +  \frac{\im }{2}\sigma _3  v_a)  R|  \sigma _3   \partial
_\omega \Phi \rangle      \right ) dD_a
-\im  \frac{q}{q'} \langle x_aR,   \partial
_\omega \Phi \rangle  dv_a.\end{aligned}
\end{equation}
Set  $\psi (U):=\frac{1}{2}\langle \sigma _3\Phi
  | R\rangle .$ We have
\begin{equation} \label{eq:dPsi0} \begin{aligned} & d\psi =
\frac{1}{2}\langle \sigma _3\Phi | \partial _\omega R \rangle
d\omega +\frac{1}{2}  \langle \sigma _3\Phi | \xi _j \rangle \left (
dz_j -d\overline{z}_j\right ) +\frac{1}{2}\langle P_c^*(\omega )\sigma _3\Phi |
  f'   \rangle  .
\end{aligned}
\end{equation}
We have $\nabla _f\psi =\frac{1}{2}P_c^*(\omega _0)P_c^*(\omega )\sigma _3\Phi $. We will use the notation $\partial _{j}:=\partial _{z_j}$ and
$\partial _{\overline{j}}:=\partial _{\overline{z}_j}$.
The sum of the last three terms in the rhs of \eqref{eq:beta2} equals

\begin{equation} \label{eq:beta31}\begin{aligned} &
- \partial _j \psi dz_j -  \partial  _{\overline{j}} \psi d\overline{z}_j
-\langle \nabla _f\psi | f'\rangle  -
 \frac{1}{2}\langle P _{N^\perp _g(\mathcal{H} ^*_\omega)}\Phi | \sigma _3 \partial _\omega R \rangle     d\omega \\& -\frac{\im}{2}
    \langle P _{N  _g(\mathcal{H} _\omega)}\Phi |  \sigma _3R \rangle        \, d \vartheta   -\frac{\im }{4}   \langle P _{N  _g(\mathcal{H} _\omega)}\Phi |x_aR \rangle
    dv_a\\&  +
 \frac{1}{2} \langle P _{N  _g(\mathcal{H} _\omega)}\Phi |
 \sigma _3(\partial _a+\im \sigma _3\frac{v_a}{2})R \rangle
  dD_a  .\end{aligned}
\end{equation}
Then  we get

\begin{equation} \label{eq:beta4}\begin{aligned} &
\widetilde{B} =   \frac{{z}_j d\overline{z}_j-\overline{z}_j dz_j
}{2}+\frac{1}{2}\langle f |\sigma _3\sigma
_1P_c(\omega  )f '   \rangle
- \partial _j \psi dz_j   \\& -  \partial  _{\overline{j}} \psi d\overline{z}_j
-\langle \nabla _f\psi |  f'\rangle -\im
 \left (   Q +\left \langle  \frac{q}{q'}    \partial
_\omega \Phi   - \frac{1}{2}   \sigma _3 P _{N  _g(\mathcal{H} _\omega)}\Phi | R \right \rangle     \right )  \, d \vartheta  \\& +\left (
 \frac{1}{2}\langle \sigma _1 R| \sigma _3 \partial _\omega R \rangle  -\left \langle \frac{q}{q'}  \partial _\omega \Phi  +
 \frac{1}{2}  P _{N^\perp _g(\mathcal{H} ^*_\omega)}\Phi | \sigma _3 \partial _\omega R \right \rangle   \right ) d\omega \\&  + \left (  \im    \Pi _a+
\left \langle  \frac{q}{q'}  \partial _\omega \Phi
-\frac{1}{2} P _{N  _g(\mathcal{H} _\omega)}\Phi |
 \sigma _3(\partial _a+\im \sigma _3\frac{v_a}{2})R \right \rangle
\right ) dD_a \\& -\frac{\im }{2}
\left ( \frac{1}{2}
 \langle \sigma _1 R|x_aR\rangle  +
    \left \langle   \frac{q}{q'} \partial _\omega \Phi  -\frac{1}{2}   P _{N  _g(\mathcal{H} _\omega)}\Phi |x_aR \right \rangle
\right )   dv_a.\end{aligned}
\end{equation}
By \eqref{eq:vectdec} have the following two equalities:
\begin{equation} \label{eq:beta41}\begin{aligned} & \frac{1}{2}   P _{N  _g(\mathcal{H} _\omega)}\Phi =  \frac{\langle   \Phi |
   \Phi \rangle }
  {2q' }
    \partial _\omega \Phi =\frac{q}{q'} \partial _\omega \Phi  , \\& \Phi =\frac{q}{q'}  \partial _\omega \Phi  +
 \frac{1}{2}  P _{N^\perp _g(\mathcal{H} ^*_\omega)}\Phi .
\end{aligned}
\end{equation}
By \eqref{eq:beta41} we have   various cancelations  in \eqref{eq:beta4} yielding
\begin{equation} \label{eq:beta5}\begin{aligned} &
\widetilde{B}  =\frac{{z}_j d\overline{z}_j-\overline{z}_j dz_j
}{2}+\frac{1}{2}\langle f |\sigma _3\sigma
_1P_c(\omega  )f '  \rangle -  d\psi-\im Q  d\vartheta   -\im \Pi _a dD_a\\& +
\frac{1}{2} \langle  \sigma _1\sigma _3R|\partial _\omega R \rangle d\omega -\frac{\im }{4}
 \langle \sigma _1 R|x_aR\rangle dv_a. \end{aligned}
\end{equation}
We have $\Omega =d (\widetilde{B}+d\psi )=d \mathrm{B}$ if we define

\begin{equation} \label{eq:beta6}\begin{aligned} &
\mathrm{B}:=\im  \vartheta dQ  +\im  D_a d\Pi _a -  \frac{\overline{z}_j dz_j -
{z}_j d\overline{z}_j}{2} +\frac{1}{2}\langle f  | \sigma _3\sigma
_1P_c(\omega  )f '  \rangle \\& + \frac{1}{2}\left ( \langle  \sigma _1\sigma _3R|\partial _\omega R \rangle +\im  v_a \langle  \sigma _1 R|x_a\partial _\omega R \rangle  \right )d\omega \\& + \frac{\im}{2}  v_a \langle  x_a\sigma _1 R|\xi _j \rangle   dz_j
+ \frac{\im}{2}  v_a \langle  x_a R|\xi _j \rangle   dz_j+ \frac{\im}{2}  v_a \langle  x_a \sigma _1 R|P_c(\omega  )f'\rangle . \end{aligned}
\end{equation}
 Finally, this $\mathrm{B}$ satisfies \eqref{eq:1forms1}--\eqref{eq:1forms2} if
 we consider the formula
  \begin{equation} \label{eq:beta7}\begin{aligned} &
  (q'+\partial _\omega Q(R))d\omega =dQ-  \partial _j Q(R) dz_j -  \partial  _{\overline{j}} Q(R) d\overline{z}_j
-\langle \nabla _fQ(R) |f'\rangle
    \end{aligned}
\end{equation}
where $  \nabla _fQ(R)=P_c^*(\omega _0)P_c^*(\omega  )\sigma _1R.$ \qed

\bigskip
In general, given a function $F$, the following formula defines $\nabla _fF$:
 \begin{equation}  \begin{aligned} &
 dF=  \partial _{Q}F dQ +\partial _{\Pi _a}Fd\Pi _a +\partial _{\vartheta}F dQ +\partial _{D_a}FdD_a \\&
    +\partial _{ {j}} F dz_j+\partial _{\overline{j}} Fdz_{\overline{j}}+\langle \nabla _fF | f'    \rangle ,
    \end{aligned}\nonumber
\end{equation}
where $P_c^*(\omega _0)\nabla _fF=\nabla _fF$.
 We have the following result.

\begin{lemma}
  \label{lem:vectorfield0}   Let us denote by $\mathcal{U}_{\Sigma}$ the subset of $\Sigma$
  defined by the inequalities
   $|z|  \le   \varepsilon _0$, $  \| f\| _{H ^{-K,-S}}\le   \varepsilon _0$,
   $|\varrho (f) | \le   \varepsilon _0$ (here recall $\varrho (f):=\left ( Q(f), \Pi (f) \right )$),
   $|\omega -\omega _0|\le   \varepsilon _0$ and $|v|\le   \varepsilon _0$.
  Then there exists a  number $ \varepsilon _0>0$ such that
     there exists a unique vectorfield  $ \mathcal{X}^t :U  _{ \Sigma} \to  L^2 $       which solves
$
 i_{\mathcal{X}^t} \Omega _t=- \mathrm{\Gamma }
   $, where $\Omega _t:=\Omega _0+t(\Omega -\Omega _0)$.

\end{lemma}

\begin{remark}
\label{rem:vfield}   In \cite{Cu1} the existence of $\mathcal{X}^t$ is elementary and standard, due to the
   fact that $\Omega _0$ and $ \Omega$ are differential forms
   in $L^2$ very close to each other. This is not true any more in our setting,
   where $\Omega _0 :\Sigma \to B^2 (L^2,\C )$ and were  we cannot assume that
$\Omega _0$ and $ \Omega$ are close in $\Sigma$. This for the simple reason that while the spaces $\Sigma _n$ are invariant for  our NLS \eqref{NLS}, see Lemma \ref{lem:sigma} below, for $n>0$ the $\Sigma _n$   norms of the solutions grow in time. One of the main differences between
 \cite{Cu1}  and this paper lies here.
 \end{remark}

\proof The proof ends after Lemma \ref{lem:fred1}.
 We  are considering $i_{\mathcal{X}^t} \Omega _0+t\, i_{\mathcal{X}^t}d\mathrm{\Gamma}=- \mathrm{\Gamma }  $ where
\begin{equation} \label{eq:contr0}  \begin{aligned}
 & i_{X}\Omega _0= \im ({X})_{\vartheta}dQ + \im ({X})_{D_a}d\Pi _a-\im ({X})_{Q}d\vartheta - \im ({X})_{\Pi _a}dD_a
 \\& - ({X})_{\overline{j}}dz_j+({X})_{j}dz_{\overline{j}}+\langle\sigma _1\sigma _3({X})_{f}| f'    \rangle ,
 \end{aligned}
\end{equation}
  $({X})_{\vartheta}$   the $\vartheta$--th component of $X$, etc.
Set
\begin{equation} \label{eq:gamma-1}   \begin{aligned}
  \gamma _{ 1a}:&=\langle x_a \sigma _1 R| \xi _j \rangle   dz_j+
   \langle  x_a   R| \xi _j \rangle  d\overline{z}_j + \langle    P^* _c(\omega )\sigma _1x_a   R  | f'    \rangle \, \, , \\   W
     :&= \varphi dQ-
   \varphi  \partial _{ { j}} Q(R) d {z}_j  -
   \varphi  \partial _{\overline{ j}} Q(R) d\overline{z}_j \\& -
   {\varphi}  \langle   P^* _c(\omega )\sigma _1 R| f'    \rangle   - {2}^{-1}
      \langle   \sigma _1\sigma _3P  _d(\omega ) f| f'    \rangle .
   \end{aligned}
\end{equation}
We have $\mathrm{\Gamma} = \frac{\im }{2} v_{a}\gamma _{ 1a}+ W $   and

\begin{equation} \label{dgamma}   \begin{aligned}
 &  d\Gamma =\frac{\im }{2}dv_{a}\wedge \gamma _{ 1a} +\frac{\im }{2} v_{a}d \gamma _{ 1a}+ d W .
 \end{aligned}
\end{equation}

\begin{lemma}
\label{lem:1diff for}
Using definitions \eqref{eq:gammatilde}--\eqref{eq:dgamma13} below, we have:
\begin{equation} \label{gamma1}  \begin{aligned}
 & \gamma _{ 1a} =  \langle    P^* _c(\omega )\sigma _1x_a   f  | f'    \rangle +\widetilde{ \gamma} _{ 1a};
 \end{aligned}
\end{equation}
  \begin{equation} \label{eq:dv}  \begin{aligned}
 &   d v _{  a}  =-\frac{2}{Q}d\Pi _a +2\left [ \left (\Pi _a-\Pi _a(R)\right )Q^{-2} +    \frac{Q^{-1} \partial _\omega \Pi _a(R)}{q'+\partial _\omega Q(R)}
 \right ] dQ \\& +2Q^{-1}\left \langle \im    P^* _c(\omega )\sigma _1 \sigma _3 \partial _af - \frac{  \partial _\omega \Pi _a(R) P_c^*(\omega  )\sigma _1f}{q'+\partial _\omega Q(R)}  |f'\right \rangle
  + \widetilde{d v _{  a}};
 \end{aligned}
\end{equation}
\begin{equation}  \label{eq:dgamma1} \begin{aligned}
 &   d \gamma _{ 1a} = \widehat{\gamma }_{1a}   \wedge  (dQ- \langle P_c^*(\omega  )\sigma _1f |f'\rangle ) + \widetilde{d \gamma} _{ 1a} .
   \end{aligned}
\end{equation}
In the above formulas we have
\begin{equation} \label{eq:gammatilde}  \begin{aligned}
 &
 \widetilde{\gamma} _{ 1a}:&=\langle x_a \sigma _1 R| \xi _j \rangle   dz_j+
   \langle  x_a   R| \xi _j \rangle  d\overline{z}_j + \langle    P^* _c(\omega )x_a    (z\cdot \sigma _1\xi +\overline{z}\cdot  \xi  | f'    \rangle \, ;
\end{aligned}
\end{equation}

\begin{equation} \label{eq:dvtilde}  \begin{aligned}
 & \frac{Q}{2}  \widetilde{d v _{  a}}  =
  \left ( Q^{-1}
 \partial _j \Pi _a(R)-    \frac{ \partial _\omega \Pi _a(R)\partial _j  Q(R)}{q'+\partial _\omega Q(R)}
 \right ) dz^j\\&
 + \left ( Q^{-1}
 \partial  _{\overline{j}} \Pi _a(R)-    \frac{  \partial _\omega \Pi _a(R)\partial _{\overline{j}} Q(R)}{q'+\partial _\omega Q(R)}
 \right ) d\overline{z}^j-\\&
    \langle \im    P^* _c(\omega ) \sigma _3\partial _a(z\cdot  \sigma _1\xi + \overline{z}\cdot   \xi ) + \frac{  \partial _\omega \Pi _a(R) P_c^*(\omega  ) (z\cdot \sigma _1 \xi + \overline{z}\cdot   \xi ) }{q'+\partial _\omega Q(R)}  |f'  \rangle  ;\end{aligned}
\end{equation}
\begin{equation}  \label{eq:dgamma12} \begin{aligned}
 &  \widehat{\gamma }_{1a} := \frac{\partial _\omega \langle x_a \sigma _1 R| \xi _j \rangle dz_j + \partial _\omega \langle x_a   R| \xi _j \rangle d\overline{z}_j - \langle  \partial _\omega  P^* _d(\omega )\sigma _1x_a   R  | f'    \rangle }{q'+\partial _\omega Q(R)} ;\end{aligned}
\end{equation}
\begin{equation}  \label{eq:dgamma13} \begin{aligned}&\widetilde{d \gamma} _{ 1a}  = (\partial _k Q(R) dz_k  +  \partial  _{\overline{k}} Q(R) d\overline{z}_k +\langle P_c^*(\omega  )(z\cdot \sigma _1 \xi +\overline{z}\cdot \xi ) |f'\rangle )\wedge  \widehat{\gamma }_{1a}  .
\end{aligned}
\end{equation}
  \end{lemma}
\proof    By \eqref{eq:variables} we   get
\begin{equation}   \begin{aligned}
 &   d v _{  a} = -\frac{2}{Q}d\Pi _a +2( \Pi _a-\Pi _a(R))Q^{-2}dQ +2Q^{-1}
 \partial _j \Pi _a(R)dz^j\\& +2Q^{-1}
 \partial  _{\overline{j} }\Pi _a(R)d\overline{z}^j+2Q^{-1}\langle \nabla _f\Pi _a(R)|f'\rangle +2Q^{-1} \partial _\omega \Pi _a(R)d \omega ,
 \end{aligned}
 \nonumber
\end{equation}
where $\nabla _f\Pi _a(R)= \im  P^* _c(\omega _0) P^* _c(\omega )\sigma _1 \sigma _3 \partial _aR .$
Then by \eqref{eq:beta7} we get

\begin{equation}   \begin{aligned}
 &   d v _{  a} =- \frac{2}{Q}d\Pi _a +2\left [ \left (\Pi _a-\Pi _a(R)\right )Q^{-2} +    \frac{Q^{-1} \partial _\omega \Pi _a(R)}{q'+\partial _\omega Q(R)}
 \right ]dQ \\&
 +2Q^{-1}\left [ Q^{-1}
 \partial _j \Pi _a(R)-    \frac{ \partial _\omega \Pi _a(R)\partial _j  Q(R)}{q'+\partial _\omega Q(R)}
 \right ] dz^j\\&
 +2Q^{-1}\left [ Q^{-1}
 \partial  _{\overline{j}} \Pi _a(R)-    \frac{  \partial _\omega \Pi _a(R)\partial _{\overline{j}} Q(R)}{q'+\partial _\omega Q(R)}
 \right ] d\overline{z}^j\\& +
 2Q^{-1}\left \langle \im    P^* _c(\omega )\sigma _1 \sigma _3 \partial _aR - \frac{  \partial _\omega \Pi _a(R) P_c^*(\omega  )\sigma _1R}{q'+\partial _\omega Q(R)}  |f'\right \rangle .\end{aligned}
 \nonumber
\end{equation}
 In particular we get \eqref{eq:dv} and \eqref{eq:dvtilde}. \eqref{eq:gammatilde}   is immediate from \eqref{eq:gamma-1}.   \eqref{eq:dgamma12}--\eqref{eq:dgamma13} follow by  \eqref{eq:beta7}    from
\begin{equation}   \begin{aligned}
 &   d \gamma _{ 1a} =\left (\partial _\omega \langle x_a \sigma _1 R| \xi _j \rangle dz_j + \partial _\omega \langle x_a   R| \xi _j \rangle d\overline{z}_j - \langle  \partial _\omega  P^* _d(\omega )\sigma _1x_a   R  | f'    \rangle  \right )\\& \wedge \frac{dQ-  \partial _k Q(R) dz_k -  \partial  _{\overline{k}} Q(R) d\overline{z}_k
-\langle P_c^*(\omega  )\sigma _1R |f'\rangle}{ q'+\partial _\omega Q(R) }  \\& =\widehat{\gamma}_{1a} \wedge  ((dQ
-\langle P_c^*(\omega  )\sigma _1R |f'\rangle ) +( \partial _k Q(R) dz_k -  \partial  _{\overline{k}} Q(R) d\overline{z}_k ) ).
 \end{aligned}
 \nonumber
\end{equation}

     \qed

 \bigskip
 Set  $H_c^{K,S} =P_c(\omega_0 )H ^{K,S}$ and denote
\begin{equation}\label{eq:PhaseSpace}  {\Ph}^{K,S}=\R ^6\times \mathbb{C }^m\times H_c^{K,S}\, , \quad {\Ph} ={\Ph}^{0,0} .
\end{equation}

 The following lemma is straightforward.
 \begin{lemma}\label{lem:regForms}For $\varpi _1$=$\widetilde{d v _{  a}}$,
  $\widetilde{ \gamma} _{ 1a}$,  $\widehat{ \gamma} _{ 1a}$, we have
$ \varpi _1\in  C^\infty (\U  _{ \Sigma} , B ( {\Ph} ^{-K,-S} ,\C ) )$    with for fixed $C$ \begin{equation} \label{gamma11}  \begin{aligned}
 &\| \varpi _1
    \| _{B ( {\Ph} ^{-K,-S} ,\C ) } \le C (|z|+ \norma{f}_{H^{-K,-S}} ) . \end{aligned}
\end{equation}
For $\varpi _2$=$ \widetilde{d \gamma} _{ 1a}$  we have
$ \varpi _2\in    C^\infty (\U  _{ \Sigma} , B^2( {\Ph} ^{-K,-S} ,\C ) )$  with for fixed $C$
\begin{equation} \label{eq:star}  \begin{aligned}
 &  \| \varpi _2
    \| _{B^2( {\Ph} ^{-K,-S} ,\C ) }\le C (|z|+ \norma{f}_{H^{-K,-S}} ) ^2 . \end{aligned}
\end{equation}
Furthermore, consider the contraction operator $X\to i_X \varpi _2$ and define $Y(X)$
by $i_{Y(X)} \Omega _0 =i_X \varpi _2$. Then we have an equality of the form,
summing on repeated index, with finite sums,
\begin{equation} \label{eq:star1}  \begin{aligned}
   & (Y(X)) _{\underline{j}}= R _k^{(\underline{j})}(X) \widetilde{\lambda} _k^{(\underline{j})} \text{ for $\underline{j}=j,\overline{j}$}  \\&
   (Y(X)) _{f}= R _k (X)  {\Lambda} _k \end{aligned}
\end{equation}
with $R _k^{(\underline{j})},R _k\in    C^\infty (\U  _{ \Sigma} , B ( {\Ph} ^{-K,-S} ,\C ) )$ satisfying  \begin{equation} \label{eq:star2}  \begin{aligned} \| R _k^{(\underline{j})}
    \| _{B ( {\Ph} ^{-K,-S}  ,\C ) }+ \| R _k
    \| _{B ( {\Ph} ^{-K,-S}  ,\C ) }\le C (|z|+ \norma{f}_{H^{-K,-S}} ) ,
\end{aligned}
\end{equation}
  with   $\widetilde{\lambda} _k^{(\underline{j})}\in    C^\infty (\U  _{ \Sigma} ,  \C )$ and
$\Lambda _k \in    C^\infty (\U  _{ \Sigma} ,  H ^{ K, S}  )$ satisfying estimates
\begin{equation} \label{eq:star3}  \begin{aligned}  | \widetilde{\lambda} _k^{(\underline{j})}
     |  + \| \Lambda _k
    \| _{  H ^{ K, S}   }\le C (|z|+ \norma{f}_{H^{-K,-S}} ).
\end{aligned}
\end{equation}
 For  $E=\varpi _1$,  $\varpi _2$,  $R _k^{(\underline{j})}$,  $R _k  $, $\widetilde{\lambda} _k^{(\underline{j})} $ and ${\Lambda} _k$,
 we have $E=E (Q,\Pi , z,f, \varrho (f))$, where  $E (Q,\Pi , z,f, \varrho  )$
    is smooth w.r.t.  $(Q,\Pi)$ and $\varrho \in \R ^4,$ $z\in \C^m$ and $f\in H^{-K,-S}$.\qed
\end{lemma}

\begin{lemma}
\label{lem:doublev} Set

\begin{equation} \label{eq:fuctsG}  \begin{aligned}
 &   G_1:=\frac{1}{2} \langle \sigma _1\sigma _3R|
  \partial _\omega R \rangle \, , \quad G_ {2a}:= \frac{\im }{2}  \langle  x_a\sigma _1 R|
  \partial _\omega R \rangle \ ,
\end{aligned}\nonumber
\end{equation}
\begin{equation}  \label{eq:What}  \begin{aligned}
 &  \widehat{ W  } :=  ( \partial _j G_1+v_a\partial _j G_ {2a})    d {z}_j +
     ( \partial _{\overline{ j}} G_1+ v_a\partial _{\overline{ j}} G_ {2a})   d\overline{z}_j  + \\&
      \langle    \nabla _f G_1 +v_a\nabla _fG_ {2a}  +  {2}^{-1}\sigma _1\sigma _3\partial _\omega P  _d(\omega ) f | f'    \rangle  +G_ {2a} \widetilde{d v _{  a}} \, ,
\end{aligned}\nonumber
\end{equation}
\begin{equation}  \label{eq:dWtilde}  \begin{aligned}
 & \widetilde{d   W  } : = -  \widehat{ W  }
\wedge  \frac{\partial _j Q(R) dz_j +  \partial  _{\overline{j}} Q(R) d\overline{z}_j
+\langle P_c^*(\omega  ) (z\cdot \sigma _1\xi +\overline{z}\cdot  \xi ) |f'\rangle }{q'+ \partial _\omega Q(R)}  \, .
\end{aligned}\nonumber
\end{equation}
Then we have
\begin{equation} \label{eq:dGam2}  \begin{aligned}
 &  d   W    = \widehat{ W  }   \wedge   \frac{dQ
  -
 \langle P_c^*(\omega ) \sigma _1 f |f'\rangle }{q'+ \partial _\omega Q(R)} + \widetilde{d   W  }+\\&  G_ {2a} (dv_a -\widetilde{d v _{  a}})
\wedge \frac{dQ-  \partial _j Q(R) dz_j -  \partial  _{\overline{j}} Q(R) d\overline{z}_j
-\langle P_c^*(\omega  )\sigma _1R |f'\rangle}{ q'+\partial _\omega Q(R) },
   \end{aligned}
\end{equation}
where, for $\varpi _1=\widehat{ W  }$ and $\varpi _2= \widetilde{d   W  }$,
the conclusions of Lemma \ref{lem:regForms} are satisfied.

  \end{lemma}
\proof    Substituting  \eqref{eq:beta7} and using \eqref{eq:1forms2}  we get
\begin{equation} \label{eq:newW}  \begin{aligned}
 &  W  =   (G_1+v_aG_{2a}) d\omega    - {2}^{-1}
      \langle   \sigma _1\sigma _3P  _d(\omega ) f| f'    \rangle .
\end{aligned}
\end{equation}
By \eqref{eq:beta7} we obtain the following formula, from which we get \eqref{eq:dGam2}:
\begin{equation}   \begin{aligned}
 &  d   W     =
  \left ( dG_1+v_adG_{2a} +2^{-1} \langle   \sigma _1\sigma _3
\partial _\omega P  _d(\omega ) f| f'    \rangle   + G_ {2a} dv_a  \right) \wedge d\omega
 \\& =    \left (  \widehat{W} +G_ {2a} (dv_a -\widetilde{d v _{  a}})    \right )   \wedge d\omega .
   \end{aligned}
 \nonumber
\end{equation}
Then substitute \eqref{eq:beta7}.
     \qed

 \bigskip
We reframe the equation for $\mathcal{X}^t$.

\begin{lemma}
  \label{lem:fred1}
For $Y$ defined by  $i_Y\Omega _0 =-\Gamma$. Then equation $
 i_{\mathcal{X}^t} \Omega _t=- \mathrm{\Gamma }
   $,  is equivalent  to
\begin{equation} \label{eq:fred1}  \begin{aligned}
 & (1+ t\mathcal{K})\mathcal{X}^t=Y
\end{aligned}
\end{equation}
where the operator $\mathcal{K}$ satisfies the following properties.

\begin{itemize}
\item[(1)]
For the  component  $(\mathcal{K} X)_f$ the following facts hold.
\begin{equation}  \label{eq:fred11}   \begin{aligned}
 &   (\mathcal{K} X)_f= A _a(X) \partial _a f+ \sigma _3(B _a(X) x _a  + C(X))f+ D_i(X) \Psi _i
\end{aligned}
\end{equation}
where the last is a finite sum with $\Psi _i \in C^\infty (\mathcal{U}_{\Sigma} , H^{K,S})$ with
\begin{equation}  \label{eq:estPsii}   \begin{aligned}
 &   \| \Psi _i    \|  _{ H^{ K, S} }\le C (|z|+ \norma{f}_{H^{-K,-S}} ).
\end{aligned}
\end{equation}
For $L= A _a, B _a,C, D_i$, we have
  $ L \in C^0 (\mathcal{U}_{\Sigma} , B({\Ph} ,\C ))$, see \eqref{eq:PhaseSpace}, with
\begin{equation}   \label{eq:fred12} \begin{aligned}
 &    L(X)  =\left \langle L  _{1b} \sigma _1\sigma _3\partial _b f+ L  _{2b}  \sigma _1x _b f  +L  _{3  }  \sigma _1  f|  (X)_f  \right \rangle + \widetilde{L}  (X)
\end{aligned}
\end{equation}
with $ \widetilde{L}   \in C^\infty (\mathcal{U}_{\Sigma} , B({\Ph} ^{-K,-S},\C ))$,  $\| \widetilde{L}  \|  _{B({\Ph} ^{-K,-S},\C ) } \le C   (|z|+ \norma{f}_{H^{-K,-S}} ) $ and where $   L  _{1b} $, $  L  _{2b} $, $  L  _{3} $ are in $C^\infty (\mathcal{U}_{\Sigma} , \C )$.

\item[(2)]
The $z_j$--th component  $(\mathcal{K} X)_j  $   is of the form \eqref{eq:fred12}
with the estimates

\begin{equation}  \label{eq:discrCom}   \begin{aligned}
 &   |L  _{1b} |+|  L  _{2b} |+|  L  _{3}|+\| \widetilde{L}  \|  _{B({\Ph} ^{-K,-S},\C ) }\le C (|z|+ \norma{f}_{H^{-K,-S}} )  .
\end{aligned}
\end{equation}
\item[(3)]
For  $G$ any of the above $L  _{1b}$,  $L  _{2b}$,  $L  _{3}$,  $\widetilde{L}  $,
 we have $G=G (Q,\Pi , z,f, \varrho (f))$, where  $G (Q,\Pi , z,f, \varrho  )$
    is smooth w.r.t.  $(Q,\Pi)$ and $\varrho \in \R ^4,$ $z\in \C^m$ and $f\in H^{-K,-S}$.
\end{itemize}

\end{lemma}\proof
 We have $i_{\frac{\partial}{\partial Q}}\Omega _0= -\im d\vartheta   $
by \eqref{eq:Omega0} and
$i_{\frac{\partial}{\partial Q}}\Omega  = -\im d\vartheta  $ by Lemmas \ref{lem:Involutions} and \ref{lem:reduction1}. Similarly $i_{\frac{\partial}{\partial \Pi _a}}\Omega _0=i_{\frac{\partial}{\partial \Pi _a}}\Omega= -\im dD_a  $. So in particular $i_{\frac{\partial}{\partial Q}}d\mathrm{\Gamma}  =
i_{\frac{\partial}{\partial D_a}}d\mathrm{\Gamma}  =0$.  Then   $(\mathrm{\Gamma})_{\vartheta}= (\mathrm{\Gamma})_{D_a}=0$  implies $(\mathcal{X}^t)_{Q}=(\mathcal{X}^t)_{\Pi _a}=0$.
Then

\begin{equation}\label{eq:iX0}    \begin{aligned}
 &   i_{\mathcal{X}^t}d\Gamma =\frac{\im }{2} \left (dv_{a} (\mathcal{X}^t)\gamma _{ 1a} -  \gamma _{ 1a} (\mathcal{X}^t)dv_{a} +  v_{a} i_{\mathcal{X}^t}d \gamma _{ 1a}\right ) + i_{\mathcal{X}^t}d W  \\&  = \frac{\im }{2}  dv_{a} (\mathcal{X}^t) \left ( \langle    P^* _c(\omega )\sigma _1x_a   f  | f'    \rangle +\widetilde{ \gamma} _{ 1a} \right ) -\frac{\im }{2}\gamma _{ 1a} (\mathcal{X}^t)  \left (  \widetilde{d v _{  a}}+ {d v _{  a}}- \widetilde{d v _{  a}}   \right ) \\&+ \frac{\im }{2}
  v_{a}  \widehat{\gamma }_{1a} (\mathcal{X}^t)   (dQ- \langle P_c^*(\omega  )\sigma _1f |f'\rangle ) + \frac{\im }{2}
  v_{a} \langle P_c^*(\omega  )\sigma _1f |(\mathcal{X}^t)_f\rangle \widehat{\gamma }_{1a} \\&+
  \widehat{ W  } (\mathcal{X}^t) \frac{ dQ
  -
 \langle P_c^*(\omega ) \sigma _1 f |f'\rangle }{q'+ \partial _\omega Q(R)}   +
 \frac{ \langle P_c^*(\omega ) \sigma _1 f |(\mathcal{X}^t)_f\rangle}{q'+ \partial _\omega Q(R)} \widehat{ W  }
  +i_{\mathcal{X}^t}\widetilde{d   W  }
  \\& +
 G_ {2a} (dv_a(\mathcal{X}^t) -\widetilde{d v _{  a}}(\mathcal{X}^t))
  \frac{dQ-  \partial _j Q(R) dz_j -  \partial  _{\overline{j}} Q(R) d\overline{z}_j
-\langle P_c^*(\omega  )\sigma _1R |f'\rangle}{ q'+\partial _\omega Q(R) }
   \\& + G_ {2a}\frac{  \partial _j Q(R) (\mathcal{X}^t)_j +  \partial  _{\overline{j}} Q(R) (\mathcal{X}^t)_{\overline{j}}
+\langle P_c^*(\omega  )\sigma _1R |(\mathcal{X}^t)_f\rangle}{ q'+\partial _\omega Q(R) }
  (dv_a -\widetilde{d v _{  a}}).
  \end{aligned}\nonumber
\end{equation}
So for    $\widetilde{ Q}$ and $\widetilde{\Pi }_a$ two  1--forms irrelevant in the sequel (since we are not interested about
$(\mathcal{X}^t)_{\vartheta}$ and $(\mathcal{X}^t)_{D _a} $), after a tedious but elementary computation we have

\begin{equation}\label{eq:iX}    \begin{aligned}
 &   i_{\mathcal{X}^t}d\Gamma =  \widetilde{ Q} (\mathcal{X}^t) dQ  +\widetilde{\Pi }_a (\mathcal{X}^t) d\Pi _a  +\frac{\im }{2}  dv_{a} (\mathcal{X}^t)\widetilde{\gamma} _{ 1a} - \frac{\im }{2} \gamma _{ 1a} (\mathcal{X}^t)\widetilde{d v _{  a}} \\&+     \frac{\im}{2} v_a  i_{\mathcal{X}^t}\widetilde{d \gamma}_{1a}  + i_{\mathcal{X}^t}\widetilde{d W }  +\langle  P_c^*(\omega )\sigma _1f |(\mathcal{X}^t)_f\rangle \left( \frac{\im}{2}v_a\widehat{\gamma }_{1a}
  + \frac{ \widehat{W } }{q'+\partial _\omega Q(R) }\right )\\& -(q'+\partial _\omega Q(R))^{-1}  G_{2a} (dv_a(\mathcal{X}^t)-\widetilde{dv_a}(\mathcal{X}^t)) (\partial _j Q(R)
  dz_j+\partial  _{\overline{j}} Q(R)
  d\overline{z}_j)
  \\& + \widehat{\Gamma}_3(\mathcal{X}^t)
  \left \langle P_c^*(\omega )\sigma _1f|f' \right \rangle   +\widehat{\Gamma}_{2a}(\mathcal{X}^t) \left \langle P_c^*(\omega )\sigma _1x_af|f' \right \rangle \\& +
  \widehat{\Gamma}_{1a}(\mathcal{X}^t)  \left \langle P_c^*(\omega )\sigma _1\sigma _3 \partial _af|f' \right \rangle ,
 \end{aligned}
\end{equation}
with by \eqref{eq:dv}
\begin{equation}\label{eq:iX1}    \begin{aligned}
 &     \widehat{\Gamma}_3:= \frac{-    \widehat{W } }{q'+ \partial _\omega Q(R)}    -\frac{\im}{2}v_a \widehat{\gamma} _{1a}   +\im Q^{-1}\frac{  \partial _\omega \Pi _a(R) }{q'+\partial _\omega Q(R)} \gamma _{1a}  \\& -    \frac{2G_{2a}}{Q(q'+\partial _\omega (R))} \left \langle \im    P^* _c(\omega )\sigma _1 \sigma _3 \partial _af - \frac{  \partial _\omega \Pi _a(R) P_c^*(\omega  )\sigma _1f}{q'+\partial _\omega Q(R)}  |f'\right \rangle \\& -2 G_{2a}\partial _\omega \Pi _a(R) \frac{  \partial _k Q(R) dz_k +  \partial  _{\overline{k}} Q(R) d\overline{z}_k
+\langle P_c^*(\omega  )\sigma _1R |f'\rangle}{ Q(q'+\partial _\omega Q(R))^2 }  \\&
 \widehat{\Gamma}_{2a}:=\frac{\im }{2} dv_{a}\\&
 \widehat{\Gamma}_{1a}:=
  \frac{\gamma _{1a}   }{Q}  +2  G_{2a}\frac{  \partial _k Q(R) dz_k +  \partial  _{\overline{k}} Q(R) d\overline{z}_k
+\langle P_c^*(\omega  )\sigma _1R |f'\rangle}{Q( q'+\partial _\omega Q(R)) }  .
 \end{aligned}
\end{equation}
The operator $\mathcal{K}$ is defined by
  $i_{\mathcal{K} X}\Omega _0= i_{X}d\Gamma$ for any $X$.
By \eqref{eq:contr0} this implies that
\begin{equation} \label{eq:KXf} \begin{aligned}
&  (\mathcal{K} X)_f   =   \widehat{\Gamma}_3(X)
 P_c  (\omega _0 )  P_c (\omega )\sigma _3f +\widehat{\Gamma}_{2a}(X)
 P_c  (\omega _0 )  P_c (\omega )\sigma _3x_af +(\widehat{\resto} ( X ) )_f  \\&  +\widehat{\Gamma}_{1a}(X)
 P_c  (\omega _0 )  P_c (\omega )\partial _af+ \langle P_c^*(\omega ) \sigma _1f| (X)_f\rangle \sigma _3\sigma _1 \left ( \frac{\im}{2} v_a(\widehat{\gamma}_{1a})_f+(\widehat{W} )_f\right ) ,\end{aligned}
\end{equation}
where $ i_{\widehat{\resto} ( X ) }\Omega _0:=    \frac{\im}{2} v_a  i_{X}\widetilde{d \gamma}_{1a}  + i_{X}\widetilde{d W} $
and where $(\widehat{\gamma}_{1a})_f$ resp. $(\widehat{W} )_f$ are the analogues of
$(\Gamma) _f$ of the expansion of $\Gamma$ under \eqref{eq:1forms1}.
They  are like the $\Psi _i$ of the statement by Lemma \ref{lem:regForms} resp. \ref{lem:doublev}.  By Lemma
\ref{lem:regForms} we have    $(\widehat{\resto} ( X ) )_f = R_k(X) \Lambda _k $,   a sum  of the form
\eqref{eq:star1} which satisfies \eqref{eq:star2}--\eqref{eq:star3}.
Claim (1) follows   from \eqref{eq:KXf} after expansions like
$P_c  (\omega _0 )  P_c (\omega )\sigma _3f =\sigma _3f +(1-P_c  (\omega _0 )  P_c (\omega ))\sigma _3f $, where the second term on the right is like a $\Psi _i$.

The   terms in \eqref{eq:iX} contributing to $(\mathcal{K} X)_j  $
   are the last two in the first line and those in  the second  and third lines. Then Claim (2) follows from
   Lemma \ref{lem:regForms}. Claim (3) follows from Lemmas \ref{lem:1forms} and \ref{lem:regForms}.
\qed

\bigskip

  We can apply  Fredholm alternative to prove existence and uniqueness of
a solution to \eqref{eq:fred1}.We check that $(1+t\mathcal{T} )X=0$ admits no solution, with $\mathcal{T} $   the adjoint of $\mathcal{K}$ with respect to $\Omega _0$.   $\mathcal{T}  $ is like  $\mathcal{K}$ and satisfies \eqref{eq:fred11}--\eqref{eq:discrCom}. We show that   $X=-t\mathcal{T} X $  does not have nontrivial
solutions for $|t| \le 3$ for the $\varepsilon _0 $,  see the statement of Lemma
\ref{lem:vectorfield0},
sufficiently small. We prove this by showing that there is a fixed constant $\kappa$ such that the following holds:
\begin{equation}  \label{eq:fred2}   \begin{aligned}
 &  ( \mathcal{T}  ^nX)_f= A _a^{(n)}(X) \partial _a f+ \sigma _3[B _a^{(n)}(X) x _a  + C^{(n)}( X)]f+ D_i^{(n)}(X) \Psi _i  ,\\&  ( \mathcal{T}  ^nX) _{\underline{j}}=Z_{\underline{j}}^{(n)}(X) , \\& | L^{(n)}(X)|\le \kappa ^{n-1}   (|z|+\| f\| _{H^{-K,-S}} + | \varrho (f) |  )^{n-1}| L (X)|,
\end{aligned}
\end{equation}
   for $L^{(n)}= A _a^{(n)}$, $B _a^{(n)}$, $C^{(n)},$ $D _1^{(n)},$ $Z_{\underline{j}}^{(n)}$. This will imply that the solutions of $X=-t\mathcal{T} X $ are trivial for $ \kappa \varepsilon _0$ sufficiently small.   We have
\begin{equation}\label{eq:fred3}   \begin{aligned}
 &  L^{(n+1)}(X)=  \widetilde{L}  (\mathcal{K} ^nX) \\& +\left \langle L  _{1b} \sigma _1\sigma _3\partial _b f+ L  _{2b}  \sigma _1x _b f  +L  _{3  }  \sigma _1  f|   (\mathcal{T}  ^nX)_f  \right \rangle.
 \end{aligned}
\end{equation}
We have $| \widetilde{L}  (\mathcal{K} ^nX) | \le c_0\kappa ^{n-1}   (|z|+\| f\| _{H^{-K,-S}} +| \varrho (f) | )^{n }| L (X)|$ for a fixed constant $c_0=c_0(\widetilde{L})$, by induction and Lemma \ref{lem:fred1}. Substituting \eqref{eq:fred2}, the second line in   \eqref{eq:fred3}
  becomes
\begin{equation}\label{eq:fred4}   \begin{aligned}
 &    \langle L  _{1b} \sigma _1\sigma _3\partial _b f+ L  _{2b}  \sigma _1x _b f  +L  _{3  }  \sigma _1  f|    A _a^{(n)}(X) \partial _a f+\cdots + D_i^{(n)}(X) \Psi _i   \rangle \\ =   &\delta _{a b} ( L  _{1b}B _a^{(n)}(X) -L  _{2b}  A _a^{(n)}(X)) \, Q(f)\\  +& D_i^{(n)}(X) \langle L  _{1b} \sigma _1\sigma _3\partial _b f+ L  _{2b}  \sigma _1x _b f  +L  _{3  }  \sigma _1  f|  \Psi _i   \rangle ,
 \end{aligned}
\end{equation}
where we used
\begin{equation}\label{eq:fred41}   \begin{aligned}
 &      \langle \sigma _1\sigma _3\partial _a f | \partial _b f \rangle = \langle \sigma _1 \partial _a f |   f \rangle =\langle \sigma _1\sigma _3x _a f | x _b f \rangle = \langle \sigma _1\sigma _3x _a f |   f \rangle \\& =\langle \sigma _1\sigma _3  f |   f \rangle =0 \text{ and  }  2\langle \sigma _1 \partial _a f | x _b f \rangle =  - \delta _{ab}\| f\| _{L^2}^2 .
 \end{aligned}
\end{equation}
The absolute value of the rhs of \eqref{eq:fred4} is by induction
  \begin{equation}  \begin{aligned}
 &    \le c_1\kappa ^{n-1}   (|z|+\| f\| _{H^{-K,-S}} +| \varrho (f) | )^{n }| L (X)|.
 \end{aligned}\nonumber
\end{equation} for a fixed constant $c_1=c_1(L)$.
So
\begin{equation} \label{eq:fred5}    \begin{aligned}
 &   | L^{(n+1)}(X)|\le \kappa ^{n+1}  (|z|+\| f\| _{H^{-K,-S}} +| \varrho (f) | )^{n }| L (X)|
 \end{aligned}
\end{equation}
  if   the constant $\kappa$ is chosen sufficiently large. The proof of Lemma \ref{lem:vectorfield0} is concluded.
The fact that $| L (X)|\le c_3  (|z|+\| f\| _{\Sigma}) \| X\| _{L^2}  $ does not need to be used.

\qed

The following one is another   most important step  in the proof. We
need to show that the flow of $\mathcal{X}^t $ corresponds to the
flow of a quasilinear hyperbolic symmetric system. To study this
system, well posedness and regularity with respect to the initial
data, we need more information on $\mathcal{X}^t $. We remark that $
H^{ K , S }$   has been fixed with any preassigned pair  $( K , S
)$. We will take both $K$ and $S$ very large.

\begin{lemma}
  \label{lem:vectorfield}  For   $\mathcal{X}^t $ the   vectorfield
        of Lemma \ref{lem:vectorfield0}, we have
\begin{equation}\label{eq:quasilin1}
\begin{aligned} &
    (\mathcal{X}^t)_f = \mathcal{L} f  + \mathcal{D} \quad , \quad     (\mathcal{X}^t)_{\underline{j}}  =Z_{\underline{j}}   \text{ for $\underline{j}\in\{ j, \overline{j} \}$},  \\& \mathcal{L}:= \mathcal{A }_a \partial _a  +( \mathcal{B}_ax_a + \mathcal{C})  \sigma _3 \end{aligned}
\end{equation}
where the coefficients satisfy the  following properties.

 \begin{itemize}
\item[(1)]$ \mathcal{A }_a $ are real valued functions.  $\mathcal{B }_a$ and $\mathcal{C} $ are imaginary valued. $Z_{\underline{j}}$ has values in $\C  $  with $Z_{\overline{j}} =\overline{ Z}_{j}$. $\mathcal{D}$ has values in $  H^{ K , S }$.

\item[(2)]      We have $G=G(t,z,f,\varrho (f))$
for $G=\mathcal{A }_a,\mathcal{B }_a,\mathcal{C },\mathcal{D},Z_{\overline{j}}$,
     for functions  $G(t,z,f,\varrho  )$    smooth in $t$, $z$, in $f\in H^{-K,-S}$ and in
$\rho  $.

 \item[(3)] We have   \begin{equation} \label{eq:vectorfield1}
\begin{aligned} &|\mathcal{A }_a|  \le C(   K, S) (|z| ^2+\| f \| _{H^{-K,-S}}^2 +|\varrho (f)| ) \, ,  \\& |Z|  + |\mathcal{C }|  +\norma{\mathcal{D}}_{H^{ K, S}}   \\& \le C(   K, S) (|z|+\| f \| _{H^{-K,-S}} +|\varrho (f)| ) (|z|+\| f \| _{H^{-K,-S}})\ .
\end{aligned}
\end{equation}

     \item[(4)] We have  $\mathcal{B }_a=-\frac{ \im }{2} v_a-\frac{ \im }{2}tdv_a(\mathcal{X}^t) $.

\end{itemize}

\end{lemma}
\proof  Let us start assuming  that  $\mathcal{X}^t $ is of the form \eqref{eq:quasilin1}.
Let $X$ be a vector such that $\sigma _1X=\overline{X}$. Then $\mathrm{\Gamma} (X)$ is imaginary.
We have $\overline{\Omega _t (\mathcal{X}^t,X)} =-\Omega _t (\sigma _1\overline{\mathcal{X}^t},X)  $. For $t=1$ is straightforward
and for $t=0$ can be checked using Lemmas \ref{lem:gradient z}
 and \ref{lem:gradient f}.
Since also $\overline{\Omega _t (\mathcal{X}^t,X)} =-\mathrm{\Gamma} (X)=\Omega _t (\mathcal{X}^t,X)$ we get $\sigma _1\overline{\mathcal{X}^t}=\mathcal{X}^t$.
 From this discussion we can conclude  that Claim (1) holds if \eqref{eq:quasilin1}
 is true.

Let $Y$ be defined by $i_Y\Omega _0 =-\Gamma$. Then

\begin{equation} \label{eq:ExprY}\begin{aligned} & (Y)_j =- \mathrm{\Gamma}_{\overline{ j}} \, , \quad (Y)_{\overline{ j}} =  \mathrm{\Gamma}_{j}   \\&
 (Y)_{f}=\sigma _3\sigma _1
\mathrm{\Gamma}_f =  -\frac{\im}{2}
\sigma _3v_ax_af+\varphi \sigma _3f+\widetilde{Y}_f,
\end{aligned}
\end{equation}
  with $  ({Y})_{\underline{j}}   $ and  $ \widetilde{Y}_f$, like $G=G(t,z,f,\varrho (f))$ in the statement above,
   smooth in  $z$, in $f\in H^{-K,-S}$ and in
$\rho  $ and s.t., by  \eqref{eq:1forms2},
\begin{equation}\label{eq:ExprY1}  |   ({Y}) _{\underline{j}}  |  + \|  \widetilde{Y}_f \|  _{H^{K,S}}
\le C (|z|+\|   f \|  _{H^{-K,-S}})(|z|+\|   f \|  _{H^{-K,-S}}+|\varrho (f)|) . \end{equation}
  Our \textit{first  claim}  is that  the following series converge:

\begin{equation}\label{eq:neumann}\begin{aligned} &
(\mathcal{X}^t)_f=\sum _{n=0}^{\infty} (-1)^n t^n (\mathcal{K}^nY )_f \\& =: \widehat{A} _a (Y) \partial _a f+ \sigma _3(\widehat{B} _a (Y) x _a  + \widehat{C} ( Y))f+ \widehat{D}_i (Y) \Psi _i,\\& (\mathcal{X}^t) _{\underline{j}}=\sum _{n=0}^{\infty} (-1)^n t^n (\mathcal{K}^nY )_{\underline{j}}=:  \widehat{Z}_{\underline{j}} (Y)  ,\end{aligned}
\end{equation}
with,  for $L= A _a, B _a,C, D_i, Z_{\underline{j}}$,
\begin{equation}\label{eq:neumann1}\begin{aligned} &
\widehat{L}   (Y)=\sum _{n=0}^{\infty} (-1)^n t^n L  ^{(n)}(Y),
\end{aligned}
\end{equation}
   for $L^{(n)}= A _a^{(n)}$, $B _a^{(n)}$, $C^{(n)},$ $D _i^{(n)},$ $Z_{\underline{j}}^{(n)}$ defined as in \eqref{eq:fred2}
   but with $\mathcal{K}$ instead of $\mathcal{T}$. To prove the \textit{first  claim}, notice that   by the proof
   of Lemma \ref{lem:fred1} we can conclude that there exists a fixed
   constant  $\kappa$ such that the following  analogue of  \eqref{eq:fred5}  holds:
   \begin{equation} \label{eq:fred6}    \begin{aligned}
 &   | L^{(n+1)}(Y)|\le \kappa ^{n+1}  (|z|+\| f\| _{H^{-K,-S}} +| \varrho (f) | )^{n }| L (Y)|\le  \\&  \kappa ^{n+1}  (|z|+\| f\| _{H^{-K,-S}} +| \varrho (f) | )^{n  }(|z|^2+\| f\| _{H^{-K,-S}}^2 +| \varrho (f) | ) ,
 \end{aligned}
\end{equation}
where in the second inequality we exploited \eqref{eq:ExprY}--\eqref{eq:ExprY1}  and Lemma \ref{lem:fred1} which yield
\begin{equation} \label{eq:fred7}    \begin{aligned}
 &    | L (Y)|  \le C_L(|z|^2+\| f\| _{H^{-K,-S}}^2 +| \varrho (f) | ).
 \end{aligned}
\end{equation}
 In the case of $L=   \mathcal{C},Z , D_i $  we have   a better
 estimate:
\begin{equation} \label{eq:fred71}    \begin{aligned}
 &    | L (Y)|    \le C_0(|z| +\| f\| _{H^{-K,-S}}    )(|z| +\| f\| _{H^{-K,-S}}  +| \varrho (f) | ).
 \end{aligned}
\end{equation}
This yields the \textit{first claim}, proves that $\mathcal{X}^t$ is of the
form of \eqref{eq:quasilin1} and that \eqref{eq:vectorfield1} holds.
A  \textit{second  claim}  is  that
\begin{equation}\label{eq:fred31}   \begin{aligned}
 &  L^{(n+1)}(Y)=  \sum _{i_1\in I_1 } \varsigma _{i_1}  (Y) \sum _{ \sigma , \tau }   \varepsilon _{i_1 }(\sigma , \tau )\prod _{j=0}^{n } \ell  _{\sigma (j)} \varpi _{\tau (j)}          ,
 \end{aligned}
\end{equation}
where we have what follows.
\begin{itemize}
\item[(i)] The sums on $ \sigma$ and $\tau$ are over all
maps $\sigma :\Z _n\to I_2$ and $\tau :\Z _n\to I_3$
with
$I_j$, for $j=1,2,3$, three finite sets  described below.   For $ i_1\in I_1$ are functions $ (\sigma , \tau )\to \varepsilon _{i_1 }(\sigma , \tau )$  with values in $\{ 0,1 \}$.

\item[(ii)]
$\varsigma  _{i_1} (Y)$ varies
in all possible ways among $ \langle \sigma _1\sigma _3\partial _a f|Y\rangle $, $  \langle  \sigma _1x _a f|Y\rangle $, $  \langle  \sigma _1  f|Y\rangle $ and  $ \widetilde{L} (Y)$,
for $L= A _a, B _a,C, D_i,(\mathcal{K}\quad ) _{\underline{j}}$. We denote by $I_1$ the set of these functions.

\item[(iii)] $ \ell _{i_2}  $ varies   among  $ L_{1b}$, $ L_{2b}$, $ L_{ 3}$,
for $L= A _a, B _a,C, D_i,(\mathcal{K}\quad ) _{\underline{j}}$. We denote by $I_2$ the set of these functions.

 \item[(iv)]
 $ \varpi _{i_3} $ varies among  1, $\| f\| _2^2$, $ \langle \sigma _1\sigma _3\partial _a f|\Psi _i\rangle $, $  \langle  \sigma _1x _a f|\Psi _i\rangle $, $  \langle  \sigma _1  f|\Psi _i\rangle $,  $ \widetilde{L} (\Psi _i)$,
 $ \widetilde{L} (\Psi _i)$,  $ \widetilde{L} (\partial _a f)$, $ \widetilde{L} (\sigma _3x _af)$, $ \widetilde{L} (\sigma _3 f)$ for $L= A _a, B _a,C, D_i,(\mathcal{K}\quad ) _{\underline{j}}$. We denote by $I_3$ the set of these functions.

   \end{itemize}
  The \textit{second  claim} is true for $n=0$ since by \eqref{eq:fred12}
   \begin{equation}    \begin{aligned}
 &    L(Y)  = L  _{1b}  \langle \sigma _1\sigma _3\partial _b f|  Y   \rangle + L  _{2b} \langle \sigma _1x _b f |  Y   \rangle  +L  _{3  } \langle \sigma _1  f|  Y   \rangle  + \widetilde{L}  (Y).
\end{aligned}\nonumber
\end{equation}
  In this case $\varpi _{i_3} =1$. Suppose that the \textit{second  claim} holds for $n-1$.
  For definiteness we will look at the case of $\mathcal{A}^{(n+1)}_a(Y)$. We have
  \begin{equation}    \begin{aligned} &
     \mathcal{A}^{(n+1)}_a(Y) = \widetilde{\mathcal{A} _a}  (\mathcal{K} ^nY) +  (\mathcal{A} _a)  _{1b}  \langle  \sigma _1\sigma _3\partial _b f|   (\mathcal{K}  ^nY)_f    \rangle\\&  + (\mathcal{A} _a)  _{2b} \langle \sigma _1x _b f |   (\mathcal{K}  ^nY)_f    \rangle +(\mathcal{A} _a)  _{3  } \langle \sigma _1  f|   (\mathcal{K}  ^nY)_f    \rangle
   .
 \end{aligned}\nonumber
\end{equation}
For definiteness let us look   at the last term on the first line.  Then, substituting
 \begin{equation}     \begin{aligned}
 &  ( \mathcal{K}  ^nY)_f= A _a^{(n)}(Y) \partial _a f+ \sigma _3[B _a^{(n)}(Y) x _a  + C^{(n)}( Y)]f+ D_i^{(n)}(Y) \Psi _i
\end{aligned}\nonumber
\end{equation}
 and using the induction hypothesis, we obtain the desired result.
 By proceeding in the same  way in all the other cases we get the \textit{second  claim}.

 What is left in the proof of Lemma \ref{lem:vectorfield}
 is  the regularity in Claim (2). The functions in (i)--(iv) are like
 the $G(t,z,f,\varrho )$ in the statement of Claim (2). But then
 Claim (2) follows by the   elementary  fact that, if $f_0$,...,$f_N$ are smooth scalar functions, if $\sigma   $ varies  in the set $\mathfrak{F}(\Z _n,\Z _N)$ of  all the
 maps   $\Z _n\to \Z _N$ and we consider arbitrary maps $\varepsilon _n : \mathfrak{F}(\Z _n,\Z _N)
 \to \{ 0 ,1 \}$, then there exists a fixed $\delta >0$ such that if
 $ |f_j|_{\infty}< \delta $ for all $j$, then the following series
 converges to a smooth function:
\begin{equation}\label{eq:series}
 \sum _{n=0}^{\infty} \sum _{\sigma \in \mathfrak{F}(\Z _n,\Z _N)}\varepsilon _n (\sigma ) \prod _{j=0}^{n} f_{\sigma (j)}.
\end{equation}
We sketch a proof assuming that the $f_j$ are functions of one scalar variable.
It is enough to show that the series obtained differentiating term by term in  \eqref{eq:series} are totally convergent.
This is immediate for the 0--th derivative. The $k$--th  derivative yields, for $ | \mu |= \sum _j |\mu (j)|,$ a series of the form \begin{equation}\label{eq:series1}
 \sum _{n=0}^{\infty} A_n \text{ with } A_n=\sum _{\substack{\sigma \in \mathfrak{F}(\Z _n,\Z _N) \\ \mu \in \mathfrak{F}(\Z _n,\Z _{k }) \text{ s.t. $| \mu |=k$} }}\varepsilon _n (\sigma ) \prod _{j=0}^{n} \,   f_{\sigma (j)}  ^{( \mu (j ))}.
\end{equation}
Then we have the bound
   \begin{equation} \begin{aligned} &
| A_n|\le  (N+1) ^{n+1}  (n+1) ^{k }    \sup   \left \{|\prod _{j=0}^{n} \,   f_{\sigma (j)}  ^{( \mu (j ))}|  \text{ s.t. $( {\sigma , \mu })$ as in \eqref{eq:series1}} \right \}\\& \le  (N+1) ^{n+1}  (n+1) ^{k } \delta ^{n-k}  \sup _{j}\| f_j\| _{W^{k,\infty}} ^k  .\end{aligned}\nonumber
\end{equation}
So  $
| A_n  |\le  \delta ^{n +1-k} (N+1) ^{n+1} (n+1) ^{k } C    _k^k  $     for   $C_k =\sup _{j}\| f_j\| _{W^{k,\infty}} $. Then for  $ (N+1)\delta _0<1$ and for $\delta \in (0,\delta _0)$ the series \eqref{eq:series1} is convergent for any $k$.

\qed

 Having established  the existence and a number of properties
 of $\mathcal{X}^t$, in Sect. \ref{sec:ODE}  we prove in an abstract set up a number
 of results on vectorfields. After a preliminary section on the spaces
 $\Sigma _\ell$, in  Sect. \ref{sec:Quasilinear} we check that it is possible
 to apply the theory  in Sect. \ref{sec:ODE} to  appropriate generalizations of
   $\mathcal{X}^t$.

\section{Some results on abstract ODE's}
\label{sec:ODE}

We collect a number of results needed for Darboux Theorem and the method
of normal forms.
We will consider a system

\begin{equation} \label{eq:ODE}\begin{aligned} &
  \dot x= f(t,x) \, , \quad x(0)=\underline{x}.
\end{aligned}   \end{equation}
We assume the following set up.
\begin{itemize}
\item[(1)]  We consider five separable Hilbert spaces $\mathbb{E}_i$ with $i=0,4$ s.t. $\mathbb{E}_{i }\subset \mathbb{E}_{i+1}$ for $i<4$,
$\mathbb{E}_{i }$ is dense in $\mathbb{E}_{i+1}$ and the   immersion $j ^{(i)}$ of $\mathbb{E}_{i}$ in $ \mathbb{E}_{i+1}$ is compact. We denote by $( \, ,\, )_i$ resp. $\| \, \| _{i}$ the inner product resp. the norm in $\mathbb{E}_i$.

\item[(2)] We assume the existence of $\{ j_{\epsilon}  : \epsilon >0  \} \subset
B(\mathbb{E} _{i+1}, \mathbb{E} _{i})$  for $i=0,...,3$
s.t.:  $\| j_{\epsilon} \circ  j^{(i)} \| _{B(\mathbb{E}_{i} , \mathbb{E}_{i} )}\le C$ for fixed $C$ for all $ \epsilon >0$;
$\lim _{\epsilon \searrow 0}j^{(i)}\circ j_{\epsilon}  =\uno _{\mathbb{E}_{i+1}}$ in $B(\mathbb{E} _{i+1}, \mathbb{E}_{i+1} )$;  $\lim _{\epsilon \searrow 0}j _{\epsilon}\circ j^{(i)} =\uno _{\mathbb{E}_{i }}$ in $B(\mathbb{E} _{i }, \mathbb{E}_{i } )$.

\item[(3)] Let  $\mathfrak{B}_i$ be the
  neighborhood of $0\in    \mathbb{E}_i$ defined by $\| x\| _4 < c_0$, for some fixed $c_0>0$. Then   $f\in C^n ((-3,3)\times \mathfrak{B}_i , \mathbb{E}_{i+1}) $ for a $n\ge 1$.

\item[(4)] The following inequalities hold for a positive constant  $C(\lambda )$
which  is increasing  functions of  $\lambda $:

\begin{eqnarray} \label{eq:ODE1}&
 \| f(t,x) \|  _{ {i+1}} \le C (\| x \|  _{ {4}} ) \| x \|  _{ {i }};
 \\  \label{eq:ODE2}&  | ( j _\epsilon f(t,x) , x ) _{ {i }} |\le C (\| x \|  _{ 4} ) \| x \|  _{ {i }} ^2  \text{ $\forall \, \epsilon \in (0,1)$ and $  i$}; \\ \nonumber &  | (   \partial _x^kj  _{\epsilon}f(t,x)(u,v)  ,  v ) _{i+1} |\le C (\| x \|  _{i}  )  \| u\| _{\mathbb{E}_i^{k-1}}  \|    v\|  _{i+1} ^2     \\ \label{eq:ODE6}&  \text{$\forall$ $1\le k\le n$,  $\epsilon >0$ ,  $i=0,...,3$ and $v\in \mathbb{E}_i$}.\end{eqnarray}

\end{itemize}
The main results of this section are the three Proposition
\ref{prop:flows1}  and \ref{prop:flows3}.

\begin{proposition}\label{prop:flows1}   $\exists$ a neighborhood $\U$ of 0 in $\mathfrak{B}_0\subseteq \mathbb{E}_0$ s.t. $\forall \, \underline{x}\in \U$
  system \eqref{eq:ODE} has exactly one solution
  $x(t)\in \cap _{i=1}^{2}  C^{i-1} ([-2,2], \mathbb{E}_i)$.
  $\U$  can be chosen to be defined by $ \| \underline{x} \| _{4}<\varepsilon _0$. For $\varepsilon _0>0$ small enough we have for a fixed $C$
  \begin{equation}\label{eq:fl1} \begin{aligned} &
  \| x \| _{L^{\infty}([-2,2], \mathbb{E}_0)  } \le C \| \underline{x} \| _{0}  \, , \quad  \| x \| _{  W^{1,\infty}([-2,2], \mathbb{E}_1)} \le C \| \underline{x} \| _{0}\, ,\\&  \| x \| _{L^{\infty}([-2,2], \mathbb{E}_4)  } \le C \| \underline{x} \| _{4}.\end{aligned}
\end{equation}
Furthermore, denoting by $ \phi ^t $ the flow associated to \eqref{eq:ODE}, we
have $\phi ^t(\underline{x})\in  C([-2,2],C^{n } (  \U, \mathbb{E}_2))  $.
\end{proposition}
We will use also a second version of the above result.

\begin{proposition}\label{prop:flows3}  Assume  that hypotheses
(1)--(4) hold,  but only with four    spaces $\mathbb{E}_i$ with $i=0,3$  and  with $\| \, \|  _3$ replacing the $\| \, \|  _4$ norm.
Then there exists an $\varepsilon _0>0$ such that if    $\U$ is the subset of $\mathbb{E}_0$ defined by  $ \| \underline{x} \| _{3}<\varepsilon _0$, system \eqref{eq:ODE} has exactly one solution
  $x(t)\in \cap _{i=1}^{2}  C^{i-1} ([-2,2], \mathbb{E}_i)$. The following inequalities hold for a fixed $C$
  \begin{equation}\label{eq:fl1bis} \begin{aligned} &
  \| x \| _{L^{\infty}([-2,2], \mathbb{E}_0)  } \le C \| \underline{x} \| _{0}  \, , \quad  \| x \| _{  W^{1,\infty}([-2,2], \mathbb{E}_1)} \le C \| \underline{x} \| _{0}\, ,\\&  \| x \| _{L^{\infty}([-2,2], \mathbb{E}_3)  } \le C \| \underline{x} \| _{3}.\end{aligned}\end{equation}
Furthermore, if we denote by $ \phi ^t $ the flow associated to \eqref{eq:ODE},
then $\partial _t^i\phi ^t(\underline{x})\in  C([-2,2],C^{n } (  \U, \mathbb{E} _{2+i} ))  $ for $i=0,1$.

\end{proposition}
\proof The  proof is the same  of that of Proposition \ref{prop:flows1}
with one minor modification.
That is, in Lemma \ref{lem:ODEmoll} below, set $X=\mathbb{E}_3$ instead of
$X=\mathbb{E}_4$. The corresponding inequalities and their
proofs are exactly the same.
\qed

\textit{Proof of Proposition \ref{prop:flows1}}. The proof tailored on standard arguments, see  \cite{Taylor} after p. 360,
but in the absence of an obvious reference we review it.

We  consider   systems

\begin{equation} \label{eq:ODEmoll}\begin{aligned} &
  \dot x _{\epsilon}= j _\epsilon   f(t,x_{\epsilon}) \, , \quad x_{\epsilon}(0)=\underline{x}.
\end{aligned}   \end{equation}
We have the following lemma.

\begin{lemma} \label{lem:ODEmoll} There is  $\varepsilon _0>0$ s.t., $\forall \, \underline{x}\in\mathbb{E} _0$  with $ \| \underline{x} \| _{4}<\varepsilon _0$, system   \eqref{eq:ODEmoll} has exactly one solution $x_{\epsilon}(t)\in C^1 ([-2,2], \mathbb{E}_0)  .$ In particular
  there is a fixed $C$ s.t. for all $ \epsilon \in (0,1)$ \begin{equation}\label{eq:fl2} \begin{aligned} & \| x _{\epsilon}\| _{L ^{ \infty } ([-2,2], \mathbb{E}_0)\cap W ^{1,\infty } ([-2,2], \mathbb{E}_1)}
  \le C \| \underline{x} \| _{0},  \\& \| x _{\epsilon}\| _{  L ^{ \infty } ([-2,2], \mathbb{E}_4)}
  \le C \| \underline{x} \| _{4} . \end{aligned}
\end{equation}
\end{lemma}
\proof  By hypothesis have $\| j  _{\epsilon }  \partial ^k _xf(t,x) \|  _{0} \le C (   \epsilon ,x) $ for all $k=0,...,n$ and $x\in \mathfrak{B}_{0} $.
This implies that for  $\underline{x}\in \mathfrak{B}_{0} $   the conclusions of Lemma \ref{lem:ODEmoll}
hold in some interval $(-a_ \epsilon  (\underline{x}),a_ \epsilon (\underline{x}))$. We want to show
that $a_ \epsilon (\underline{x})>2$ if $ \| \underline{x} \| _{X}<\varepsilon _0$ for $\varepsilon _0>0$ small enough, where $X=\mathbb{E}_4$. We consider by \eqref{eq:ODE2}

\begin{equation}  \begin{aligned} & \left | \frac{d}{dt}
 \|     x_{\epsilon }   \|  _{0}^2\right | = 2  |( j  _\epsilon  f(t,x_{\epsilon }) , x_{\epsilon } ) _{0}|\le 2 C (\| x _{\epsilon }\|  _{X} ) \| x _{\epsilon }\|  _{0} ^2   .
\end{aligned}  \nonumber \end{equation}
 By Gronwall
we get $\|     x_{\epsilon }   \|  _{0} \le  e^{|t|C (1)}\|     \underline{x}    \|  _{0} $ as long as  $\| x _{\epsilon }\|  _{X}\le 1$ in $[-|t|,|t|]$.
  We claim that the latter holds for $|t|\le \frac{1-\log \varepsilon  _{0} }{C (1)}$.  By  \eqref{eq:ODE2} we have

\begin{equation}  \begin{aligned} &  \frac{d}{dt}
 \|     x_{\epsilon }   \|  _{X}^2 = 2  ( j  _\epsilon f(t,x_{\epsilon }) , x_{\epsilon } ) _{X}\le 2 C (\| x _{\epsilon }\|  _{X} ) \| x _{\epsilon }\|  _{X} ^2   .
\end{aligned}  \nonumber \end{equation}
For  $\varepsilon _0$ small enough, using Gronwall we obtain \eqref{eq:fl2} and   $|t|\ge 2$. By   equation \eqref{eq:ODEmoll} and by \eqref{eq:ODE1} we have  $\| \dot x _{\epsilon }\|  _{1}\le C' \| \underline{x}  \|  _{0}  $ for some fixed $C'$.
So Lemma \ref{lem:ODEmoll} is proved.
\qed

We now exploit Lemma \ref{lem:ODEmoll} to prove the first part of Proposition
\ref{prop:flows1}
 estimates \eqref{eq:fl1}. The argument is routine.
Given a sequence $\epsilon _{\nu}\searrow 0 $ then  $\{  x_{\epsilon _\nu}  \}$
admits a subsequence convergent in $ C  ([-2,2], \mathbb{E}_1) $
 by Ascoli Arz\'{e}la. We can assume it is the whole sequence.  We denote by $x(t)$ the limit.   Since  we have necessarily $x_{\epsilon _\nu} \to x$ in $ C  ([-2,2], \mathbb{E}_j) $ for $j>1$,   \eqref{eq:fl2} yields the third inequality in \eqref{eq:fl1}.
 Notice that $x(t)$ is the weak limit of $ x_{\epsilon _\nu} (t)$ in $\mathbb{E}_0$ for any $|t|\le 2$. By Fathou this yields the first inequality in \eqref{eq:fl1}.
 We have $f(t, x_{\epsilon _\nu} )\to f(t, x  )$ in $ C  ([-2,2], \mathbb{E}_2) $. We claim we
 have $ j  _{\epsilon _\nu}f(t,x_{\epsilon _\nu}  )  \to f(t, x  )$
 in $ C  ([-2,2], \mathbb{E}_2) $. Indeed, by  $\lim _{\epsilon \searrow 0}j _{\epsilon  } =j $ in $B(\mathbb{E}_1,\mathbb{E}_2)$ and by
 $\|  f(t,x_{\epsilon _\nu}  )\| _{\mathbb{E}_1}\le C \| \underline{x} \| _{\mathbb{E}_0} $, it follows that  in $\mathbb{E}_2$  and uniformly in $t\in [-2,2]$, we have
\begin{equation*}\begin{aligned} & j  _{\epsilon _\nu}f(t,x_{\epsilon _\nu}  )=
   (j  _{\epsilon _\nu} -j     f(t,x_{\epsilon _\nu}  ) + j  ^{(1)}     f(t,x_{\epsilon _\nu}  )  \to  f(t, x  ).\end{aligned}
  \end{equation*}
This implies $\dot x \in C  ([-2,2], \mathbb{E}_2) $ where   is the limit of $\dot x_{\epsilon _\nu} $ and satisfies the equation \eqref{eq:ODE}. By
 Fathou as before this yields the second inequality in \eqref{eq:fl1}.
 Suppose
$y \in \cap _{i=1}^{2}  C^{i-1} ([-2,2], \mathbb{E}_i)$   is  a solution of \eqref{eq:ODE}. Then by \eqref{eq:ODE6}  for $k=1$ we have for $\delta x =y-x$
\begin{equation}  \label{eq:fl31}\begin{aligned} &  \frac{d\|     \delta x  \|  _{2}^2 }{dt}
 = 2 \int _0^1 (   \partial _xf(t,(1-\tau) x +\tau  y ) \delta x, \delta x ) _{2} \le C (\| x  \|  _{1} , \| y  \|  _{1} ) \| \delta x\|  _{2} ^2  .
\end{aligned}  \end{equation}
Notice indeed that, by $\lim _{\epsilon \searrow 0}j  _{\epsilon  } =j  ^{(i)}  $ in $B(\mathbb{E}_i,\mathbb{E} _{i+1})$, \eqref{eq:ODE6} implies
\begin{equation} \label{eq:fl32}
 | (   \partial _x^k f(t,x)(u,v)  , v ) _{i+1} |\le C (\| x \|  _{i}  )  \| u\| _{\mathbb{E}_i^{k-1}}  \|   v\|  _{i+1} ^2.
\end{equation}
By Gronwall, \eqref{eq:fl31} implies  $x\equiv y$. This concludes the proof of
the first part of
Prop. \ref{prop:flows1} and of  \eqref{eq:fl1}.

We now turn to the proof of the last sentence of Prop. \ref{prop:flows1}.
$\U$ will be the neighborhood of 0 in  $ \mathbb{E}_0$  defined by  $ \| \underline{x} \| _{4}<\varepsilon _0$.
We have proved that $\phi ^t(\underline{x})= \lim _{\epsilon \searrow 0}\phi ^t_{\epsilon}(\underline{x})$ in $C([-2,2],\mathbb{E}_1)$,
with  $\phi ^t_{\epsilon}$ the flow associated to \eqref{eq:ODEmoll}.
We have $\phi ^t_{\epsilon}(\underline{x})\in C^{n } ([-2,2] \times \U, \mathbb{E}_0) $.

 \begin{lemma} \label{lem:contf}  $\phi ^t\in C( \U , \mathbb{E}_1)$  for all $t\in [-2,2]$.
\end{lemma}
\proof Given $\underline{x}, \underline{y}\in \U $ set $ \delta \phi  _{\epsilon} = \phi ^t_{\epsilon}(\underline{x})-\phi ^t_{\epsilon}(\underline{y})$. By \eqref{eq:ODE6}
and the first part of Proposition \ref{prop:flows1}
\begin{equation}  \begin{aligned} &  \frac{d}{dt}
 \|    \delta \phi  _{\epsilon}  \|  _{1}^2 = 2  ( j  _{\epsilon} ( f(t,\phi ^t_{\epsilon}(\underline{x}) ) -f(t,\phi ^t_{\epsilon}(\underline{y}) ) ) , \delta \phi  _{\epsilon} ) _{1}\le C   \| \delta \phi  _{\epsilon}\|  _{1} ^2   .
\end{aligned}  \nonumber \end{equation}
This implies $\| \phi ^t_{\epsilon}(\underline{x})-\phi ^t_{\epsilon}(\underline{y})\|  _{1}\le C' \| \underline{x}-\underline{y}\|  _{1}$   for a fixed $C'$. For $ \epsilon \searrow 0$ this yields for fixed $C''$
\begin{equation} \label{eq:fl3} \begin{aligned} &\| \phi ^t (\underline{x})-\phi ^t  (\underline{y})\|  _{1}\le C' \|\underline{ x}-\underline{y}\|  _{1}\le C''\| \underline{x}-\underline{y}\|  _{0}.\end{aligned}    \end{equation}

\qed

\begin{lemma} \label{lem:flows21} We have $\|\partial ^l_{ \underline{y}} \phi ^t_{\epsilon}(\underline{y})\| _{B^l (\mathbb{E}_0 , \mathbb{E}_1 )} +\| \partial _t\partial ^l_{ \underline{y}} \phi ^t_{\epsilon}(\underline{y})\| _{B^l (\mathbb{E}_0 , \mathbb{E}_2 )}\le C$ for a fixed $C$ and
for    all $1\le l\le n$,  $\underline{y}\in \U$ and $\epsilon \in (0,1) $.
\end{lemma}
\proof We have
\begin{equation}\label{eq:flows22} \begin{aligned} &  \partial _t \partial _{ \underline{y}} \phi ^t_{\epsilon}(\underline{y})= j   _{\epsilon}\partial _{\phi} f(t, \phi ^t_{\epsilon}(\underline{y})) \partial _{ \underline{y}} \phi ^t_{\epsilon}(\underline{y}) \, , \quad \partial _{ \underline{y}} \phi ^0_{\epsilon}(\underline{y})=\uno .
\end{aligned}
\end{equation}
  For $l>1$ we have  $\partial _{ \underline{y}}^l   \phi ^0_{\epsilon}(\underline{y})=0$ and, succinctly,
  \begin{equation}\label{eq:flows23} \begin{aligned} &
  \partial _t \partial _{ \underline{y}}^l \phi ^t_{\epsilon}(\underline{y})= j _{\epsilon}\partial _{\phi} f(t, \phi ^t_{\epsilon}(\underline{y})) \partial _{ \underline{y}}^l \phi ^t_{\epsilon}(\underline{y}) +  \\&  \text{Sym}  \sum _{k=2} ^{l} \sum _{| \alpha  |=l }   j   _{\epsilon}\partial _{\phi}^k  f(t, \phi ^t_{\epsilon}(\underline{y})) \frac{l!}{\alpha !}   \partial _{ \underline{y}} ^{\alpha _1} \phi ^t_{\epsilon}(\underline{y})  ...  \partial _{ \underline{y}} ^{\alpha _k}\phi ^t_{\epsilon}(\underline{y})  ,
\end{aligned}
\end{equation}
with $| \alpha | =\sum _j \alpha _j$ and $  \alpha  ! =\prod _j \alpha _j!$
and Sym an appropriate symmetrization, see \cite{nelson} p.7.  Fix $v\in \mathbb{E}_0$. By \eqref{eq:ODE6} and for $\partial _{ \underline{y}} \phi ^t_{\epsilon}=\partial _{ \underline{y}} \phi ^t_{\epsilon}(\underline{y})$,

\begin{equation}\label{eq:flows24} \begin{aligned} & \frac{d}{dt} \|   \partial _{ \underline{y}} \phi ^t_{\epsilon} v\| ^{2}_{1}=2 ( j   _{\epsilon}\partial _{\phi} f(t, \phi ^t_{\epsilon} ) \partial _{ \underline{y}} \phi ^t_{\epsilon} v,\partial _{ \underline{y}} \phi ^t_{\epsilon} v )_{1}     \le C \|   \partial _{ \underline{y}} \phi ^t_{\epsilon} v\| ^{2}_{1}.
\end{aligned}
\end{equation}
Since $C$ is independent of $v$, we obtain  $\|\partial  _{ \underline{y}} \phi ^t_{\epsilon}(\underline{y})\| _{B  (\mathbb{E}_0 , \mathbb{E}_1 )}\le C_1$ for some fixed $C_1$   by Gronwall.
  We assume
$\|\partial ^k_{ \underline{y}} \phi ^t_{\epsilon}(\underline{y})\| _{B^k (\mathbb{E}_0 , \mathbb{E}_1 )}\le C$
for $k<l$  by induction. Let  $ K(t, \epsilon ,\underline{y})$ be the second line of  \eqref{eq:flows23}. Then $\|K(t, \epsilon ,\underline{y})\| _{B^l (\mathbb{E}_0 , \mathbb{E}_1 )}\le C$ by induction. Fix $v\in \mathbb{E}_0^{l}$. Then, proceeding as for $l=1$ we get

\begin{equation}\label{eq:flows25} \begin{aligned} & \frac{d}{dt} \|   \partial _{ \underline{y}} ^l \phi ^t_{\epsilon}(\underline{y})v\| ^{2}_{1}\le      C \|   \partial _{ \underline{y}} ^l\phi ^t_{\epsilon}(\underline{y})v\|  _{1} + C\|   K(t, \epsilon ,\underline{y})v\|  _{1} .
\end{aligned}
\end{equation}
By Gronwall we get $\| \partial _{ \underline{y}} ^l \phi ^t_{\epsilon}(\underline{y})v \| _{ 1  }\le C _l\| v\| _{ \mathbb{E}_0^{l}}$ and so
$\| \partial _{ \underline{y}} ^l \phi ^t_{\epsilon}(\underline{y})v \| _{B^l (\mathbb{E}_0 , \mathbb{E}_1 )}\le C_l$ since the constants $C$ in \eqref{eq:flows25} and $C_l$ do not depend on $v$. By equations \eqref{eq:flows22}--\eqref{eq:flows23}
we obtain $\| \partial _t \partial _{ \underline{y}} ^l \phi ^t_{\epsilon}(\underline{y})v \| _{B^l (\mathbb{E}_0 , \mathbb{E}_2 )}\le C_l'$.
\qed

The natural embedding $B^l(\mathbb{E}_i,\mathbb{E} _{i+1}) \hookrightarrow B^l(\mathbb{E}_i,\mathbb{E}_{ i+2})$ is compact. This implies that for any fixed $y$ and any sequence
$ \epsilon _{\nu}\searrow 0$   there is a   subsequence
$\partial _{ \underline{y}} ^l \phi ^t_{\epsilon _\nu}(\underline{y}) $ convergent in $  C  ([-2,2], B^l(\mathbb{E}_0,\mathbb{E}_2)) $
to a
 $
 g^{ (l)} (t,y)$.   $\partial _t g^{ (1)} (t,\underline{y})$ exists  by \eqref{eq:flows22}, with $\partial _t\partial _{ \underline{y}} \phi ^t_{\epsilon _\nu}(\underline{y}) $ convergent to it  in $  C  ([-2,2], B (\mathbb{E}_0,\mathbb{E}_3)) $ and with

\begin{equation}\label{eq:flows26} \begin{aligned} &  \partial _t g^{ (1)} (t,\underline{y})=  \partial _{\phi} f(t, \phi ^t (\underline{y})) g^{ (1)} (t,\underline{y}) \, , \quad g^{ (1)} (0,\underline{y})=\uno  .
\end{aligned}
\end{equation}
  Suppose that
$
 g^{ (1)}_1 (t,\underline{y})$ and $
 g^{ (1)}_2 (t,\underline{y})$ are two such solutions of \eqref{eq:flows26}. Fix $v\in \mathbb{E} _0$ and
 set $\delta  g  (t,\underline{y}) =  g^{ (1)}_1 (t,\underline{y})-
 g^{ (1)}_2 (t,\underline{y})$.
  Then

\begin{equation}  \begin{aligned} &  \frac{d}{dt}
 \|    \delta  g  (t,\underline{y})v  \|  _{3}^2 = 2  (    \partial _{\phi} f(t, \phi ^t (\underline{y})) \delta  g  (t,\underline{y})v,  \delta  g  (t,\underline{y})v ) _{3}\le C   \| \delta  g  (t,\underline{y})v\|  _{3} ^2  ,
\end{aligned}  \nonumber \end{equation}
by \eqref{eq:ODE6} and the first part of Proposition \ref{prop:flows1}.
This implies $\delta  g  (t,\underline{y})v=0$ and so $
 g^{ (1)}_1 (t,\underline{y})=
 g^{ (1)}_2 (t,\underline{y})$. For $l>1$, by induction and a similar argument, we get a function $g^{ (l)} (t,\underline{y}) \in  C  ([-2,2], B^l(\mathbb{E}_0,\mathbb{E}_2)) $  satisfying $g^{ (l)} (0,\underline{y})=0$    and such that $\partial _tg^{ (l)} (t,\underline{y}) \in  C  ([-2,2], B^l(\mathbb{E}_0,\mathbb{E}_3)) $ satisfies
  \begin{equation}\label{eq:flows27} \begin{aligned} &
  \partial _t g^{ (l)} (t,\underline{y})=  \partial _{\phi} f(t, \phi ^t (\underline{y})) g^{ (l)} (t,\underline{y}) + \\&
  \text{Sym}  \sum _{k=2} ^{l} \sum _{| \alpha  |=l }   j   _{\epsilon}\partial _{\phi}^k  f(t, \phi ^t_{\epsilon}(\underline{y})) \frac{l!}{\alpha !}
    g^{ (\alpha _1)} (t,\underline{y})  ...g^{ (\alpha _k)} (t,\underline{y})
     ,
\end{aligned}
\end{equation}
 where in the second line we have only $g^{ (\ell)} (t,\underline{y})$ with $\ell <l$, which  can be assumed uniquely defined by induction.
 By repeating the previous argument we get   uniqueness also for   $k=l$.

 \begin{lemma} \label{lem:firstd} For any $t$ the map $\phi ^t :\U \to \mathbb{E}_2$ is Frech\'{e}t differentiable with $ \partial _{ \underline{y}}\phi ^t ( \underline{y})=  g^{ (1)} (t,\underline{y})$. We have $\partial _t^i\partial _{ \underline{y}}\phi ^t\in C ([-2,2]\times \U , B (\mathbb{E}_0,\mathbb{E}_{i+2}) )$ for $i=0,1$.
\end{lemma}
\proof For fixed $\underline{x}, \underline{y}\in \U $  set $ \delta \underline{x} = \underline{y}-\underline{x}$  and $\delta \phi = \phi   ^t (\underline{y})-   \phi    ^t (\underline{x}) $. Then

\begin{equation}  \begin{aligned} &
\partial _t [ \delta \phi
- g^{ (1)} (t,\underline{x}) \delta \underline{x}]=    f(t,\phi   ^t (\underline{y}) )  -        f(t,\phi  ^t (\underline{x}) )  -\partial _{\phi} f(t,\phi  ^t (\underline{x}) ) g^{ (1)} (t,\underline{x}) \delta \underline{x}  .
\end{aligned}  \nonumber \end{equation}
Then, for $ g^{ (1)}= g^{ (1)} (t,\underline{x})$, we have
\begin{equation}  \begin{aligned} &  \frac{d}{dt}
 \|   \delta \phi
- g^{ (1)}  \delta \underline{x}  \|  _{2}^2 = 2  ( \partial _{\phi} f(t,\phi  ^t (\underline{x}) ) (\delta \phi -  g^{ (1)}   \delta \underline{x} ) , \delta \phi
- g^{ (1)}   \delta \underline{x}  ) _{2} \\& +  \int _0^1 ( \partial _{\phi}^2 f(t,(1-\tau )\phi  ^t (\underline{x})+\tau \phi  ^t (\underline{y}) ) (\delta \phi )^2   , \delta \phi
  -g^{ (1)}   \delta \underline{x}  ) _{2} \, d\tau .
\end{aligned}  \nonumber \end{equation}
The last line is less than $C\|   \delta \phi
   \|  _{1}^2\|   \delta \phi
- g^{ (1)}  \delta \underline{x}  \|  _{2}$. Hence by  \eqref{eq:fl3}
\begin{equation}  \begin{aligned} &  \frac{d}{dt}
 \|   \delta \phi
- g^{ (1)}  \delta \underline{x}  \|  _{2}  \le C  \|   \delta \phi
- g^{ (1)}  \delta \underline{x}  \|  _{2}   +  C \|   \delta \phi
   \|  _{1}^2    \\& \le  C  \|   \delta \phi
- g^{ (1)}  \delta \underline{x}  \|  _{2}   +  C_1 \| \delta \underline{x}\|  _{1} ^2  .
\end{aligned}  \nonumber \end{equation}
This yields $
 \|   \delta \phi
- g^{ (1)}  \delta \underline{x}  \|  _{2}  = o ( \| \delta \underline{x}  \| _{0})  .
$
Hence   $\phi ^t :\U \to \mathbb{E}_2$ is Frech\'{e}t differentiable with $ \partial _y\phi ^t ( \underline{y})=  g^{ (1)} (t,\underline{y})$. Set   $\delta  g^{ (1)} = g^{ (1)} (t,\underline{y})-   g^{ (1)} (t,\underline{x}) $.
  Set $ \delta \underline{x} = \underline{y}-\underline{x}$ like above. Then
\begin{equation}   \begin{aligned} &  \partial _t \delta  g^{ (1)} =  \partial _{\phi} f(t, \phi ^t (\underline{y})) g^{ (1)} (t,\underline{y}) -  \partial _{\phi} f(t, \phi ^t (\underline{x})) g^{ (1)} (t,\underline{x}).
\end{aligned} \nonumber
\end{equation}
Then for fixed $v\in \mathbb{E}_0$
\begin{equation}  \begin{aligned} &  \frac{d}{dt}
 \|   \delta  g^{ (1)} v \|  _{2}^2 = 2  ( \partial _{\phi} f(t,\phi  ^t (\underline{x}) )  \delta  g^{ (1)}v, \delta  g^{ (1)}v ) _{2} \\&+  2  ( (\partial _{\phi} f(t,\phi  ^t (\underline{y}) ) -\partial _{\phi} f(t,\phi  ^t (x ))   g^{ (1)} (t,\underline{y}) v, \delta  g^{ (1)}v ) _{2}  .
\end{aligned}  \nonumber \end{equation}
By  \eqref{eq:fl32},   for a fixed $C$
we have
\begin{equation}  \begin{aligned} &  \frac{d}{dt}
 \|   \delta  g^{ (1)} v \|  _{2}  \le C    \|   \delta  g^{ (1)} v \|  _{2} +C \| \delta \underline{x} \| _{0}  \| v \| _0    .
\end{aligned}  \nonumber \end{equation}
 Hence $  \|  g^{ (1)} (t,\underline{y})-   g^{ (1)} (t,\underline{x}) \|  _{B (\mathbb{E }_0,\mathbb{E }_2 )}\le C' \| \underline{y}-\underline{x}\| _0 $ for a fixed $C'$. This and
$g^{ (1)} (t,\underline{x}) \in  C  ([-2,2], B (\mathbb{E}_0,\mathbb{E}_2)) $ for
 all $\underline{x}$ imply   $g^{ (1)}   \in  C  ([-2,2]\times \U, B (\mathbb{E}_0,\mathbb{E}_2)) $.  We obtain $\partial _tg^{ (1)}   \in  C  ([-2,2]\times \U, B (\mathbb{E}_0,\mathbb{E}_3)) $ by \eqref{eq:flows26}.
\qed

 \begin{lemma} \label{lem:highd} $\partial _ {\underline{y}}^ {l-1}\phi ^t :\U \to \mathbb{E}_2$ is Frech\'{e}t differentiable
 for any $t$ and $1\le l\le n$   with $ \partial _{\underline{y}}^l\phi ^t ( \underline{y})=  g^{ (l)} (t,\underline{y})$. We have $\partial _t^i\partial _{\underline{y}}^l\phi ^t\in C ([-2,2]\times \U , B ^l(\mathbb{E}_0,\mathbb{E}_{i+2}) )$ for $i=0,1$.
\end{lemma}
\proof Case $l=1$ is proved in
  Lemma \ref{lem:firstd}. Consider $l>1$.
By induction we can assume: $\phi ^t (\underline{x})$ admits Frech\'{e}t derivatives of order $k<l$; $\partial  _{\underline{x}}^k\phi ^t (\underline{x})=g^{ (k)} (t,\underline{x})$;
$\partial ^i_t g^{ (k)}   \in  C  ([-2,2]\times \U, B^k (\mathbb{E}_0,\mathbb{E}_{i+2})) $ for $i=0,1$.
   Fix now $\underline{x}, \underline{y}\in \U $ and $v\in \mathbb{E}_0^{l-1}.$ Set $ \delta \underline{x} = \underline{y}-\underline{x}$, $g^{ (l )}=g^{ (l )} (t,\underline{x}) $   and $\delta g^{ (l-1)} =  g^{ (l-1)} (t,\underline{y})-   g^{ (l-1)} (t,\underline{x}) $. Then
\begin{equation} \label{highd}  \begin{aligned} &
\partial _t [\delta g^{ (l-1)}v-g^{ (l )}   (v,\delta \underline{x})]
=  \partial _{\phi} f(t, \phi ^t (\underline{x})) [\delta g^{ (l-1)}v-g^{ (l )} (v,\delta \underline{x})]  \\& + (\partial _{\phi} f(t, \phi ^t (\underline{y})) -\partial _{\phi} f(t, \phi ^t (\underline{x})))
\partial _y^ {l-1}\phi ^t ( \underline{y})\cdot  v
+ F(z)  \big ] _{z=\underline{x}}^{z=\underline{y}}\cdot v \\&  -
\text{Sym}  \sum _{k=2} ^{l} \sum _{| \alpha  |=l  }    \partial _{\phi}^k  f(t, \phi ^t (\underline{x})) \frac{l!}{\alpha !}   \partial _{ \underline{x}} ^{\alpha _1} \phi ^t (\underline{x})  ...  \partial _{ \underline{x}} ^{\alpha _k}\phi ^t (\underline{x})   \cdot  (v,\delta \underline{x}),
\end{aligned}    \end{equation}
where $|\alpha |= \sum _{j=1}^{k}\alpha _j$, $ \alpha !=\prod _{j=1}^{k}\alpha _j!$ and where

 \begin{equation} \label{highd1}   \begin{aligned} & F(z):=\text{Sym}  \sum _{k=2} ^{l-1} \sum _{| \alpha  |=l-1 }    \partial _{\phi}^k  f(t, \phi ^t (z)) \frac{(l-1)!}{\alpha !}   \partial _{ z} ^{\alpha _1} \phi ^t (z)  ...  \partial _{ z} ^{\alpha _k}\phi ^t (z).
\end{aligned}    \end{equation}
 Notice that   $|\alpha _j|\le l-1$ (resp. $|\alpha _j|\le l-2$) for all $j$  in \eqref{highd} (resp.  \eqref{highd1}).
The last two lines  in \eqref{highd} can be written as
\begin{equation} \label{eq:highd3} \begin{aligned} &
  F(z)  \big ] _{z=\underline{x}}^{z=\underline{y}}\cdot v  +o(\delta \underline{x}) +   \partial _{\phi}^2  f(t, \phi ^t (\underline{x}))       (\partial _{\underline{x}}  \phi ^t ( \underline{x}) \delta \underline{x})\,  \partial _{\underline{x}}^{ l -1}  \phi ^t ( \underline{x})     \cdot   v   \\&-
\text{Sym}  \sum _{k=2} ^{l} \sum _{| \alpha  |=l  }       \partial _{\phi}^k  f(t, \phi ^t (\underline{x})) \frac{l!}{\alpha !}   \partial _{ \underline{x}} ^{\alpha _1} \phi ^t (\underline{x})  ...  \partial _{ \underline{x}} ^{\alpha _k}\phi ^t (\underline{x})   \cdot  (v,\delta \underline{x})   ,
\end{aligned}    \end{equation}
where $ \| o(\delta \underline{x})\| _2= \| v\| _{\mathbb{E}_0^{l-1}}o(\|  \delta \underline{x} \| _0)$. The last term in the first line cancels with a corresponding one in the second. What remain  in the last line \eqref{eq:highd3}
   is $ -(F'(\underline{x})\delta  \underline{x} )\cdot v,$   where we know by induction that $F(z) $ is Frech\'{e}t differentiable. Hence all \eqref{eq:highd3}  is
$o(\delta \underline{x})$.
 Then by \eqref{eq:ODE6} and Gronwall we get
\begin{equation}  \begin{aligned} &
 \|  \delta g^{ (l-1)} -g^{ (l )} (t,\underline{x})  \delta \underline{x}  \|  _{B  ^{l-1} (\mathbb{E }_0,\mathbb{E }_2)}  =  o ( \|\delta \underline{x} \| _{0})  .
\end{aligned}  \nonumber \end{equation}
So   $g^{ (l-1)}(t,\cdot )  :\U \to \mathbb{E}_2$ is Frech\'{e}t differentiable and $ \partial _xg^{ (l-1)}(t,\underline{x} )=  g^{ (l)} (t,\underline{x} )$.
 By the usual argument one shows that
  $  \|  g^{ (l)} (t,\underline{y})-   g^{ (l)} (t,\underline{x}) \|  _{B^l (\mathbb{E }_0,\mathbb{E }_2 )}\le C' \| \underline{y}-\underline{x}\| _0 $ for a fixed $C'$. This and
$g^{ (l)} (t,\underline{x}) \in  C  ([-2,2], B^l (\mathbb{E}_0,\mathbb{E}_2)) $ for
 fixed $\underline{x}$ imply that $g^{ (l)}   \in  C  ([-2,2]\times \U, B ^l (\mathbb{E}_0,\mathbb{E}_2)) $.  We obtain $\partial _tg^{ (l)}   \in  C  ([-2,2]\times \U, B^l (\mathbb{E}_0,\mathbb{E}_3)) $ by \eqref{eq:flows27}.

\qed

  The last statement of Proposition \ref{prop:flows1} is proved.

\qed

\section{Some facts on the spaces $\Sigma _l$}
\label{sec:sigma}

We need some preliminary information on the  space $\Sigma _\ell$  defined by   the norm \eqref{eq:sigma}.
Consider the space  $ \Sigma _\ell '$   with norm
 \begin{equation}\label{eq:sigma1}
\begin{aligned} &
    \| U \| _{\Sigma _\ell'} ^2 :=  \| U \| _{H^n} ^2 + \sum _{|\alpha |\le \ell} \|  \, \, x ^{\alpha } U \| _{L^2} ^2<\infty . \end{aligned}
\end{equation}
\begin{lemma}
  \label{lem:sigma} We have    $\Sigma _\ell =\Sigma _\ell'$.  $\Sigma _{\ell  }  $  is preserved by the flow of   \eqref{eq:NLSvectorial}.
  \end{lemma}
  \proof  The second statement follows from the first   by the fact that $\Sigma _{\ell  } ' $  is preserved by the flow of   \eqref{eq:NLSvectorial}, see  \cite{MS}.

  The fact that   $\Sigma _\ell' \subseteq\Sigma _n  $
 follows from Proposition 2 \cite{MS}. We prove $\Sigma _\ell' \supseteq \Sigma _\ell  $ by induction.   Set $T_a:=\im  \partial _a  + \im x_a$. We have
 $T_a^*T_a=-\partial _a^2+x_a^2-1$. Then
    \begin{equation}\label{eq:sigma2}
\begin{aligned} &
      \| U \| _{\Sigma _1} ^2 \approx  \| U \| _{L^2} ^2+ \langle (-\Delta + |x|^2 )U|\overline{U}\rangle \approx  \| U \| _{\Sigma _1 '} ^2     \end{aligned}
\end{equation}
 by  the fact that $ \langle (-\Delta + 1 )U|\overline{U}\rangle$ defines
 $H^1$. So $\Sigma _1 =\Sigma _1 '$. Suppose  by induction
 $\Sigma _{\ell-1} =\Sigma _{\ell-1} '$. Let $|\alpha |=\ell-1$. Then

 \begin{equation}\label{eq:sigma3}
\begin{aligned} &
      \| \partial ^\alpha _x U \| _{H^1} ^2 \lesssim   \| \partial ^\alpha _x U \| _{L^2} ^2+ \langle (-\Delta + |x|^2 )\partial ^\alpha _x U| \partial ^\alpha _x\overline{U}\rangle \approx  \| \partial ^\alpha _x U \| _{L^2} ^2\\& + \sum _{a=1}^{3}\|  (\im  \partial _a  + \im x_a)    \partial ^\alpha _x U \| _{L^2} ^2  \le   \| \partial ^\alpha _x U \| _{L^2} ^2\\& + \sum _{a=1}^{3}\| \partial ^\alpha _x (\im  \partial _a  + \im x_a) U \| _{L^2} ^2  + \sum _{a=1}^{3}\| [   \im  \partial _a  + \im x_a  ,\partial ^\alpha _x] U \| _{L^2} ^2 \lesssim  \| U \| _{\Sigma _\ell} ^2 ,  \end{aligned}
\end{equation}
by induction and by $\sum _{a=1}^{3}\| [   \im  \partial _a  + \im x_a  ,\partial ^\alpha _x] U \| _{L^2} \lesssim \| u\| _{H^{\ell-1}}.$   For the same
$\alpha$,
using the Fourier transform we see that \eqref{eq:sigma3} implies
also
\begin{equation}\label{eq:sigma4}
\begin{aligned} &
    \sum _{a=1}^{3}  \| x_a x ^\alpha   U \| _{L^2} ^2     \lesssim  \| U \| _{\Sigma _\ell} ^2  . \end{aligned}
\end{equation}
Then \eqref{eq:sigma3}--\eqref{eq:sigma4} imply $ \| U\| _{ \Sigma '_\ell}
\lesssim \| U\| _{ \Sigma  _\ell}$ and so $\Sigma _{\ell }\subseteq \Sigma _{\ell  } '$. Since we know already $\Sigma _{\ell } \supseteq\Sigma _{\ell  } '$, they are equal.
\qed

\begin{definition}\label{def:sigma} For any positive integer $\ell $
we denote by $\Sigma  _{-\ell} $ the space formed by the $V$ such that
the map $U\to \langle U|\overline{V}\rangle $  is
in  $B(\Sigma  _{ \ell},\C ) $. We we also set $\Sigma  _{ 0}:=L^2 $.
\end{definition}
Definition \ref{def:sigma} yields a natural Banach structure on $\Sigma  _{-\ell} $.
Notice that we have found two distinct but equivalent  norms  in $\Sigma  _{ \ell} $. A third one comes from Claim (4) in Lemma  \ref{lem:sigma1}.  In the subsequent proofs we will pick in the proofs from time to time the norms which are most convenient, the statements being unsensible to the particular choice. There will also be a corresponding implicit choice of inner products.

\begin{lemma}
  \label{lem:sigma1}  Set $h=-\Delta +|x|^2$.
  The following facts hold.

  \begin{itemize}
\item[(1)]
  The   map $
     (h +1) ^\ell :\Sigma  _\ell \to   \Sigma _{-\ell} $
  is  an isomorphism between $\Sigma _\ell$  and $\Sigma _{-\ell}$
  for any $\ell \in  \mathbb{N}$.

 \item[(2)]  $\| U\| _{\Sigma  _\ell}\approx \| (h +1)  ^{\frac{\ell}{2}} U\| _{L^2} $  for any $\ell \in  \Z$.

     \item[(3)] The   map $
     (h +1) ^{ \frac{j}{2}} :\Sigma  _\ell \to   \Sigma _{ \ell -j} $
  is  an isomorphism
  for any $(\ell ,j)\in  \mathbb{Z}^2$.

      \item[(4)] The   map $
     (h +1) ^{ \frac{j}{2}} :\Sigma  _\ell \to   \Sigma _{ \ell -j} $
  is  an isomorphism
  for any $(\ell ,j)\in  \mathbb{R}^2$  if we define $\Sigma  _\ell$ for $\ell \in \R $ setting  $\| U\| _{\Sigma  _\ell}:= \| (h +1)  ^{\frac{\ell}{2}} U\| _{L^2} $.
\end{itemize}

  \end{lemma}
 \proof The proof of (1) is an easy consequence of $T_a^*T_a=-\partial _a^2+x_a^2-1$
 and $[T_a, T_b^*]=2 \delta _{ab}$ and is skipped.  Proofs of the other claims are elementary.

 \qed

We will consider the following  mollifier:
\begin{equation}\label{eq:mollif}
\begin{aligned} &
      J_{\varepsilon} = (1-\varepsilon \Delta +\varepsilon |x|^2  )^{-1}  \text{ for $ \varepsilon>0$}. \end{aligned}
\end{equation}

\begin{lemma}
  \label{lem:sigma12} Let $s>0$. the following facts hold.

   \begin{itemize}
\item[(1)] Denote by $j: \Sigma _\ell\to \Sigma  _{\ell - s} $ the natural embedding. Then $ \lim _{\epsilon\searrow 0}J _{\epsilon}  ^{\frac{s}{2}} =j$ in $B(\Sigma _\ell, \Sigma  _{\ell- s})$.

  \item[(2)]   We have $ J _{\epsilon}  ^{\frac{s}{2}}:   \Sigma _\ell\to \Sigma  _{\ell + s} $ with $\| J _{\epsilon}  ^{\frac{s}{2}} \| _{B(\Sigma _\ell, \Sigma  _{\ell+ s})}\le C  \epsilon   ^{-\frac{s}{2}}.$

   \end{itemize}
  \end{lemma}
\proof Set $h = -\Delta +|x|^2$. For any
$U \in \Sigma _\ell$  there is a  fixed $C$ s.t.
\begin{equation}  \begin{aligned} & \| (1-J _{\epsilon} ^{\frac{s}{2}}) U\| _{\Sigma  _{\ell- s} } =  \| (h +1)  ^{\frac{\ell -s}{2}  } (1-J _{\epsilon} ^{\frac{s}{2}}) U\| _{L^2 } \\& = \| (h +1)  ^{\frac{\ell -s}{2}  } \left (1- \left (1-  \frac{\epsilon h}{\epsilon h +1} \right )  ^{\frac{s}{2}}\right ) U\| _{L^2 }\le C \epsilon \| (h +1)  ^{\frac{\ell}{2}  }   U\| _{L^2 }.
\end{aligned}\nonumber
\end{equation}
So $J _{\epsilon} ^{\frac{s}{2}} -j= O(\epsilon)$ in $B(\Sigma _\ell, \Sigma  _{\ell- s})$. The second claim follows by the Spectral Theorem
and $(1+r) ^{\frac{\ell +s}{2}}(1+\epsilon r) ^{-\frac{ s}{2}}\le  \epsilon   ^{-\frac{s}{2}}(1+r) ^{\frac{\ell  }{2}} $ for any $r\ge 0$.
\qed

\begin{lemma} \label{lem:sigma3} For any $ \ell \in \Z$ we have $x_a,\partial _a \in B(\Sigma _{ \ell},\Sigma _{ \ell -1}) $
\end{lemma}
\proof For $\ell >0$,  $f\in \Sigma _{ \ell}$   and $|\alpha |\le \ell -1$ we have
\begin{equation}
\begin{aligned} & \|  (\im \partial   +\im x  )^ \alpha \partial _a f\|  _{L^2}\le \| \partial _a (\im \partial   +\im x  )^ \alpha  f\|  _{L^2}+\|  [(\im \partial   +\im x  )^ \alpha  ,\partial _a] f\|  _{L^2}\le C \| f \| _{\Sigma _{ \ell}},
\end{aligned} \nonumber
\end{equation}
where we are using \eqref{eq:sigma3} and the fact that
$[(\im \partial   +\im x  )^ \alpha  ,\partial _a]$ is a linear combination of
$(\im \partial   +\im x  )^ \beta$ with $|\beta |=|\alpha |-1.$ The case with $x_a$ is seen to be equivalent, through the Fourier transform. The case with
$\ell \le 0$ follows by duality.
\qed

\begin{lemma} \label{lem:sigma2}There is a fixed $C>0$ s.t.   $\forall$
    $ \varepsilon \in (0,1)$ for $T=x_a, \partial _a$ and $\forall$  $\ell \in \Z$ we have
$ \varepsilon  ( \|    T J_\varepsilon\|  _{B( \Sigma_{\ell},\Sigma_{\ell} )} +\|  J_\varepsilon  T \|  _{B( \Sigma_{\ell},\Sigma_{\ell} )})  <C   . $
\end{lemma}
  \proof
    By Lemma \ref{lem:sigma12} and Lemma \ref{lem:sigma3},
for $T=x_a,\partial _a$ we have $   \epsilon \| J _{\epsilon}  Tf \|  _{ \Sigma _{\ell}} \le C \|   Tf \|  _{ \Sigma _{\ell -1}}\le C'\|    f \|  _{ \Sigma _{\ell  }} .$
Similarly $   \epsilon \| TJ _{\epsilon}   f \|  _{ \Sigma _{\ell}} \le C \epsilon \|  J _{\epsilon} f \|  _{ \Sigma _{\ell +1}}\le C'\|    f \|  _{ \Sigma _{\ell  }} .$
\qed

\begin{lemma} \label{lem:sigma4} For any $ \ell \in \Z$  and $n\in \mathbb{N}$ there is a fixed $C>0$ s.t. $\forall$     $ \varepsilon \in (0,1)$ we have $   \|    [ T, J_\varepsilon ^n]\|  _{B(\Sigma _{2\ell},\Sigma _{2\ell} )} <C    $
with $T=x_a, \partial _a$.
\end{lemma}
  \proof We have
$
  [ \partial _a , J_\varepsilon ^n] =-(\epsilon h+1) ^{-n}[ \partial _a , (\epsilon h+1) ^{-n}]
  (\epsilon h+1) ^{-n}.
$
It is elementary that this is a sum of terms $-2\epsilon (\epsilon h+1) ^{j-n}
x_a(\epsilon h+1) ^{k-n}$ with $j+k=n-1$. Then for $\ell \ge 0$
\begin{equation} \begin{aligned} &
 \| \epsilon (\epsilon h+1) ^{j-n}
x_a(\epsilon h+1) ^{k-n}f \|_{\Sigma _{\ell}}\le C \| \epsilon
x_a(\epsilon h+1) ^{k-n}f \|_{\Sigma _{\ell}} \\& \le C_1\| \epsilon
 (\epsilon h+1) ^{k-n}f \|_{\Sigma _{\ell +1}} \le C_2\|
 (\epsilon h+1) ^{k-n+1}f \|_{\Sigma _{\ell  }}\le C_3\|
 f \|_{\Sigma _{\ell  }}.
\end{aligned}\nonumber
\end{equation}
The other estimates can be proved similarly. The case $\ell <0$ follows by
duality.

  \qed

\begin{lemma} \label{lem:sigma41} For any   $n\in \mathbb{N}$ there is a fixed $C>0$ s.t. $\forall$     $ \varepsilon \in (0,1)$ we have $   \|    [ T, J_\varepsilon ^n]\|  _{B(H^1,H^1 )} <C    $
with $T=x_a, \partial _a$.
\end{lemma}
  \proof We know that $   \|    [ T, J_\varepsilon ^n]\|  _{B(L^2,L^2 )} <C    $
  from Lemma \ref{lem:sigma4} for $\ell =0$. Proceeding as above we need to show that  terms like the following ones are in $B(H^1,L^2)$:
\begin{equation} \begin{aligned} &
  \epsilon \partial _b (\epsilon h+1) ^{j-n}
x_a(\epsilon h+1) ^{k-n} =    \epsilon [ \partial _b ,(\epsilon h+1) ^{j-n}]
x_a(\epsilon h+1) ^{k-n} \\& + \delta _{ab} \epsilon   (\epsilon h+1) ^{-n-1}
 +\epsilon  (\epsilon h+1) ^{j-n}
x_a[\partial _b,(\epsilon h+1) ^{k-n}]
 \\&+
  \epsilon  (\epsilon h+1) ^{j-n}
x_a(\epsilon h+1) ^{k-n} \partial _b.
\end{aligned}\nonumber
\end{equation}
All the terms in the rhs except for the last one are in $B(L^2,L^2)$
with norm bounded uniformly in $\epsilon$.For the last one the same holds in   $B(H^1,L^2)$.
\qed

\begin{lemma} \label{lem:sigma5} Let $\ell \in \Z$ and let
$n\in \N$  such that $n+\ell \ge 1$. Then there are fixed constants  $C_\ell $
and $C$
such that for $T=\partial _a, x_a$

\begin{equation}
 \begin{aligned}\label{eq:qlin11} &
  |( {f},J_\epsilon ^{ 2n} Tf)_{\Sigma_{ 2\ell}}| \le C    \|  f \| _{\Sigma_{2\ell}} ^{2}\, , \quad
 |( {f},J_\epsilon ^{ 2n} Tf)_{H^{1}}| \le    \|  f \| _{H^{1}}^{2} .
\end{aligned}
\end{equation}
\end{lemma}
  \proof Set $h=-\Delta +|x|^2$ and $X= \Sigma_{2\ell},H^1 $.
 \begin{equation}\label{eq:qlin111}
 \begin{aligned} &    ( (\epsilon h+1)^{- 2n} \partial _af,f) _X =( [(\epsilon h+1)^{-  n} , \partial _a]f,(\epsilon h+1)^{-  n}f) _X \\&+ (  \partial _a(\epsilon h+1)^{- n}f,(\epsilon h+1)^{- n}f) _X.
\end{aligned}
\end{equation}
  We have
\begin{equation}\label{eq:qlin112}
 \begin{aligned} &    |( [(\epsilon h+1)^{- n} , \partial _a]f,(\epsilon h+1)^{- n}f) _X  |\le \| f\| _{X} \| [(\epsilon h+1)^{- n} , \partial _a]f \| _{X} \le C  \| f\| _{X}^2
\end{aligned}
\end{equation}
 by Lemmas \ref{lem:sigma4} and \ref{lem:sigma41}.
    We have for $\ell \in \Z$
\begin{equation}
 \begin{aligned} &
 (  \partial _a(\epsilon h+1)^{- n}f,(\epsilon h+1)^{- n}f) _{\Sigma_{2\ell }}  \\& = (   (  h+1)^{  \ell}\partial _a(\epsilon h+1)^{- n}f,   (  h+1)^{  \ell}  (\epsilon h+1)^{- n}f) _{L^2} \\&  =  (   [(  h+1)^{  \ell},\partial _a](\epsilon h+1)^{- n}f,   (  h+1)^{  \ell}  (\epsilon h+1)^{- n}f) _{L^2},
 \end{aligned} \nonumber
\end{equation}
where we exploited $ (\partial _a g,g)_{L^2}=(\sigma _3 x _a g,g)_{L^2}=0 $ for $\overline{g} =\sigma _1 g $. The rhs is in absolute value less than $\|f \| _{\Sigma _{\ell}}^2.$ The proof for the case   $X=H^1$ is   simpler.
  \qed

\section{Quasilinear systems}
\label{sec:Quasilinear}
We will apply the theory developed in Sect. \ref{sec:ODE}
in two distinct  forms to quasilinear systems

\begin{equation}\label{eq:quasilin10}
\begin{aligned} & \dot f = \mathcal{L} f  + \mathcal{D} \quad , \quad     \dot z =Z \\& \mathcal{L}:=\mathcal{A }_a \partial _a  +( \mathcal{B}_ax_a + \mathcal{C}) \sigma _3.   \end{aligned}
\end{equation}
with $\mathcal{L}$,  $\mathcal{D}$ and $Z$  satisfying hypotheses which we will state below.

\subsection{First type of  system}
\label{subsec:Ftype}

We consider
$ 4 {n}_i = 4 {n}_0-i4n\ge n+1\gg 1$ with $i=0,1,2,3$.
We  denote   $ \mathbb{E}_i= \C ^m \times P_c(\omega _0)\Sigma  _{4{n}_i}$  with $i=0,1,2,3$.
  Set also   $ \mathbb{E}_4= \C ^m \times P_c(\omega _0) H^{1}$ and $j _{\epsilon} =J_{\epsilon}  ^{2n }$.
We assume:

 \begin{itemize}
\item[(A1)]$ \mathcal{A }_a $ are real valued functions.  $\mathcal{B }_a$ and $\mathcal{C} $ are imaginary valued. $Z_{\underline{j}}$ has values in $\C  $  with $Z_{\overline{j}} =\overline{ Z}_{j}$.

\item[(A2)]    $\mathcal{D}$ has values in $  \Sigma  _{4 {n}_0}$
   For $G=\mathcal{A }_a,\mathcal{B }_a,\mathcal{C },\mathcal{D},Z_{\overline{j}}$  we have $G=G(t,z,f,\varrho (f))$ where $G(t,z,f,\varrho  )$ is  $C^n$ in $t$, $z$, in  $f\in \Sigma  _{-4 {n}_0} $ and in
$\rho  $.

 \item[(A3)] We have   \begin{equation} \label{eq::Ftype1}
\begin{aligned} &|\mathcal{A }_a|+|Z|  + |\mathcal{C }|  +\norma{\mathcal{D}}_{\Sigma  _{4{n}_0}} \\&  \le C  (|z|+\| f \| _{\Sigma  _{-4 {n}_0}} +|\varrho (f)| ) (|z|+\| f \| _{\Sigma  _{-4 {n}_0}}).
\end{aligned}
\end{equation}

     \item[(A4)] We have either  $\mathcal{B }_a=-\frac{ \im }{2} v_a-\frac{ \im }{2}tdv_a(\mathcal{X}^t) $  for $a=1,2,3$ or $\mathcal{B }_a\equiv 0$.

\end{itemize}
 The coefficients of  Lemma \ref{lem:vectorfield} satisfy (A1)--(A4) for any choice of $n$ and $n_0$,
 by our freedom of choice of space  $  H^{ K, S}$ and by the fact that
$  H^{ K, S}\subset \Sigma  _{4{n}_0}$ for $K> 4 {n}_0$ and $S>
4{n}_0$.  We have:

\begin{proposition}\label{prop:quasilin1}    The following facts hold.
\begin{itemize}
\item[(1)]
$\exists$ a neighborhood $\U$ of $0\in  \mathbb{E}_0$
defined by
$|z|<\varepsilon _0$ and    $\| f \| _{H^{1}}  <\varepsilon _0$,
s.t.  $\forall \, (\underline{z},\underline{f})\in \U$
  system \eqref{eq:ODE} has exactly one solution $(z(t), f(t))\in \cap _{i=1}^{2} (C^{i-1} ([-2,2], \mathbb{E}_i)\cap W ^{i-1,\infty } ([-2,2], \mathbb{E}_{i-1})).$
\item[(2)] Call $ \phi ^t $ the  flow of \eqref{eq:quasilin10}.
Then $\partial _t^i\phi ^t \in  C^{i}([-2,2],C^{n} (  \U, \mathbb{E} _{2+i}))  $ for $i=0,1$.
\item[(3)] For $(z^t,f^t)=\phi ^t (z,f)$ we have
\begin{equation}\label{eq:fl11} \begin{aligned} &
\|  z ^t    \| _{L^{\infty}( -2,2   )  }  + \|   f^t    \| _{L^{\infty}( [-2,2], H^1   )  }    \le C  (|  z|+ \| f \| _{H^1})  \, , \\&
  \| f^t    \| _{L^{\infty}( [-2,2], \Sigma _{4n_0}   )  }   \le C  (|  z|+ \| f \| _{\Sigma _{4n_0}}) .\end{aligned}
\end{equation}
\end{itemize}
\end{proposition}
\proof We will need to check that we are in the framework and the
hypotheses of Sect. \ref{sec:ODE}.
Proposition \ref{prop:quasilin1} is a consequence of Proposition  \ref{prop:flows1}  if we can prove
the hypotheses (1)--(4) in Sect.\ref{sec:ODE}. Specifically we need to prove the inequalities
in (4) Sect. \ref{sec:ODE}.
By Hypotheses (A3)--(A4) and by Lemma \ref{lem:sigma3}
we see immediately that the analogue of
\eqref{eq:ODE1} is satisfied.  \eqref{eq:ODE2} is a consequence of the
following lemma.

\begin{lemma}\label{lem:lin11} For a fixed constant $C$ and for $ C _{|z|+|\varrho (f)|}$ an increasing positive function in $|z|+|\varrho (f)|$,
we have:

\begin{equation}
 \begin{aligned}  & | Z\cdot \overline{z} | \le C |z| (|z|+\| f \| _{\Sigma  _{-4 {n}_0}} +|\varrho (f)| ) (|z|+\| f \| _{\Sigma  _{-4 {n}_0}}) ,
  \\&  |( {f}, \mathcal{D})_{\Sigma_{4{n}_i}}| \le C   \|  f \| _{\Sigma_{4n_i}} (|z|+\| f \| _{\Sigma  _{-4 {n}_0}} +|\varrho (f)| ) (|z|+\| f \| _{\Sigma  _{-4 {n}_0}})\, , \\&  |( {f}, \mathcal{D})_{H^{1}}| \le C   \|  f \| _{L^2}(|z|+\| f \| _{\Sigma  _{-4 {n}_0}} +|\varrho (f)| ) (|z|+\| f \| _{\Sigma  _{-4 {n}_0}}),
  \\&
  |( {f},J_\epsilon ^{2n} \mathcal{L}f)_{\Sigma_{ {4n}_i}}| \le C _{|z|+|\varrho (f)|}  \|  f \| _{\Sigma_{4{n}_i}} ^{2},\\&
 |( {f},J_\epsilon ^{2n} \mathcal{L}f)_{H^{1}}| \le  C _{|z|+|\varrho (f)|} \|  f \| _{H^{1}}^{2} .
\end{aligned}\nonumber
\end{equation}
 \end{lemma}
 \proof The first three inequalities follow immediately  from   (A2)--(A3).
   The last two inequalities are an immediate consequence of    the following two inequalities
 for $\widehat{T}=\uno , \partial _a, x_a$ for any $a=1,2,3$: there is a fixed $C$ s.t.
 \begin{equation}
 \begin{aligned}\label{eq:qlin110} &
  |( {f},J_\epsilon ^{2n} \widehat{T}f)_{\Sigma_{ {4n}_i}}| \le C \|  f \| _{\Sigma_{4{n}_i}} ^{2}\, ,\quad
 |( {f},J_\epsilon ^{2n} \widehat{T}f)_{H^{1}}| \le  C   \|  f \| _{H^{1}}^{2} .
\end{aligned}
\end{equation}
  \eqref{eq:qlin110} follows from Lemma \ref{lem:sigma5}.

 \qed

Finally, to finish with the proof of Proposition \ref{prop:flows1} we need to prove
that  \eqref{eq:ODE6} is true. But like for \eqref{eq:ODE2} this too is an easy consequence of
(A2)--(A4) and of Lemma \ref{lem:sigma5}.
Hence  Proposition \ref{prop:quasilin1} is
proved.

\qed

\subsection{Second type of  system}
\label{subsec:Stype}
Before setting up the system we  notice that for solutions of \eqref{eq:quasilin10} satisfying (A1)--(A4) we have
 \begin{equation}\label{eq:quasilin401}
\begin{aligned} &  \frac{d}{dt}Q(f)=\langle \sigma _1 f | \dot f\rangle = \langle \sigma _1 f |  \mathcal{D}\rangle \, , \\& \frac{d}{dt}\Pi _a (f)=\im \langle \sigma _1 \sigma _3 \partial _a f | \dot f\rangle = \im  \mathcal{B}_b \langle \sigma _1 \sigma _3 \partial _a f | \sigma _3 x _b f\rangle +\im    \langle \sigma _1 \sigma _3 \partial _a f |  \mathcal{D}\rangle \\& =2\im Q(f) \mathcal{B}_a -\im    \langle \sigma _1 \sigma _3  f |  \partial _a\mathcal{D}\rangle .
   \end{aligned}
\end{equation}
 We set $ \varrho _0(f):= Q(f)$ and $ \varrho _a(f):= \Pi _a (f)$.
 We consider the system
\begin{equation} \label{eq:quasilin402}
\begin{aligned} &  \dot   \varrho _0 = \langle \sigma _1 f |  \mathcal{D}\rangle  \quad , \quad  \dot   \varrho _a = 2\im \varrho _0 \mathcal{B}_a -\im    \langle \sigma _1 \sigma _3  f |  \partial _a\mathcal{D}\rangle , \\&\dot f = \mathcal{L} f  + \mathcal{D} \quad , \quad     \dot z =Z.
\end{aligned}
\end{equation}
We  denote   $ \mathbb{E}_i= \R ^4\times \C ^m \times P_c(\omega _0)\Sigma  _{-4{n}_{3-i}} $  with $i=0,1,2,3$.
We assume:

 \begin{itemize}
\item[(B1)] same as (A1);

\item[(B2)]    $\mathcal{D}$ has values in $  \Sigma  _{4 {n}_3}$. For $G=\mathcal{A }_a,\mathcal{B }_a,\mathcal{C },\mathcal{D},Z_{\overline{j}}$  we have $G=G(t,z,f,\varrho  )$ where $G(t,z,f,\varrho  )$ is  $C^n$ in $t$, $z$, in  $f\in \Sigma  _{-4 {n}_0} $ and in
$\rho  $;

 \item[(B3)] we have   \begin{equation} \label{eq::Stype1}
\begin{aligned} &
|\mathcal{A }_a|   \le C  (|z|^2+\| f \| _{\Sigma  _{-4{n}_0} }^2 +|\varrho  | ) \, , \\&
 |Z|  + |\mathcal{C }|  +\norma{\mathcal{D}}_{\Sigma  _{4{n}_0}}  \le C  (|z|+\| f \| _{\Sigma  _{-4{n}_0} } +|\varrho  | ) (|z|+\| f \| _{\Sigma  _{-4 {n}_0}});
\end{aligned}
\end{equation}

\item[(B4)] we have either  $|\mathcal{B }_a|\le C  |\varrho |$    for $a=1,2,3$ or $\mathcal{B }_a\equiv 0$.

\end{itemize}
Notice that, starting from   $E=\mathcal{A }_a,\mathcal{B }_a,\mathcal{C },\mathcal{D},Z_{\overline{j}}$,
    with $E=E(t,z,f,\varrho (f))$  which satisfy hypotheses (A1)--(A4) in
    Subsect.  \ref{subsec:Ftype} and substituting   $ \varrho  (f) $ with
an external parameter
$ \varrho    $, we obtain functions satisfying (B1)--(B4).

 We have:

\begin{proposition}\label{prop:quasilin2}    The following facts hold.
\begin{itemize}
\item[(1)]
$\exists$ a neighborhood $\U$ of $0\in  \mathbb{E}_0$
defined by
$|\varrho |+|z|<\varepsilon _0$ and    $\| f \| _{\Sigma  _{-4{n}_{0 }}}  <\varepsilon _0$,
s.t.  $\forall \, (\underline{z},\underline{f})\in \U$
  system \eqref{eq:ODE} has exactly one solution $(\varrho (t),z(t), f(t))\in \cap _{i=1}^{2} (C^{i-1} ([-2,2], \mathbb{E}_i)\cap W ^{i-1,\infty } ([-2,2], \mathbb{E}_{i+1})).$
\item[(2)] Call $ \Phi ^t $ the  flow of \eqref{eq:quasilin402}.
Then $\partial _t^i\Phi ^t(y)\in  C([-2,2],C^{n} (  \U, \mathbb{E}_{2+i}))  $ for $i=0,1$.
\item[(3)] Set $  (z _t,f_t, \varrho _t):=\Phi ^t (z,f,\varrho )$. Then we have for $i=0,3$
 \begin{equation}\label{eq:fl12} \begin{aligned} &
\|  ( z ^t , \|   f^t    \| _{ \Sigma _{-4n_i}       }   , \varrho ^t)   \| _{L^{\infty}( -2,2   )  }   \le C  (|  z|+ \| f \| _{\Sigma _{-4n_i}} +|\varrho |)   .\end{aligned}
\end{equation}

\item[(4)]  Call $ \phi ^t$ the  flow of \eqref{eq:quasilin10}
satisfying (A1)--(A4) and suppose that \eqref{eq:quasilin402} is the corresponding system substituting   $ \varrho  (f) $ with
an external parameter
$ \varrho    $.
When $f\in H^{1}$  and for $ (z^t,f^t):=\phi ^t (z,f)$, we have $\Phi ^t ( \varrho (f), z,f) = ( \varrho (f^t), z^t,f^t) . $

\end{itemize}
\end{proposition}
\proof  Like in Subsect. \ref{subsec:Ftype}   we can reduce to
material in Sect. \ref{sec:ODE}. Specifically, by the arguments of  Subsect. \ref{subsec:Ftype} we can apply  Proposition \ref{prop:flows3}.
This yields the    Claims (1)--(3).
Claim (4) follows from \eqref{eq:quasilin401} and the uniqueness of solutions in \eqref{eq:quasilin402}.

\qed

\subsection{Structure of the Lie transform}
\label{subsec:Ltrans}

Consider system \eqref{eq:quasilin10} such that (A1)--(A4) hold. Consider the corresponding system \eqref{eq:quasilin402} satisfying (B1)--(B4).
We denote by $\phi =\phi ^1$ the \textit{Lie transform}.

\begin{lemma}\label{lem:quasilin5}
Set      $(z',f')=\phi  (z,f)$.   Then we have
\begin{equation} \label{eq:quasilin51}
\begin{aligned} &     z'  =   z  +
 \mathcal{Z}  \, \, \quad
  & f'(x) = e^{   \sigma _3  ( \textbf{B}\cdot x +  {\gamma} ) }\tau_{\mathbf{A}}f
    +\mathcal{G} (x)
\end{aligned}
\end{equation}
with $\tau_{\mathbf{A}}f
   (x)= f(x-\mathbf{A})$, $\mathbf{A}=- \int  _0^1\mathcal{A}(\tau )d\tau ,$
   $\mathcal{Z}=\int _0^1Z(\tau )d\tau  $, $\textbf{B}=\int _0^1\mathcal{B} (\tau )d\tau$, $\gamma =\int _0^1 \left (\mathcal{C}(s)+ \mathcal{A}_a(s)\int _0^s \mathcal{B} _a(\tau ) d\tau  \right ) ds$ and $\mathcal{G}$  functions of   $(z,f, \varrho (f))$

\begin{eqnarray}\mathcal{G}(x) =  \int _0^1 e^{-\sigma _3  \int _s^1 \left (x_a \mathcal{B}_a (\tau ) +\mathcal{C} (\tau )  -\mathcal{A}_a (\tau )   \int _0 ^{\tau}\mathcal{B}_a (\tau ')d\tau '  \right ) d\tau} \mathcal{D}  (s, x+\int _s^1\mathcal{A}(\tau )d\tau   ) ds.
 \nonumber
\end{eqnarray}
We have the following estimates for a fixed constant $C>0$:

\begin{eqnarray}
\label{eq:fl13} &\|  ( z   , \|   f     \| _{ \Sigma _{-4n_0}       }   )   \| _{L^{\infty}( -2,2   )  }   \le C  (|  z|+ \| f \| _{\Sigma _{-4n_0}}  )
 \, ,  \\
  &  \label{eq:fl131}  |\mathcal{Z}|    +\norma{\mathcal{G}}_{\Sigma  _{4{n}_0}}   \le C  (|z|+\| f \| _{\Sigma  _{-4{n}_0} } +|\varrho  | ) (|z|+\| f \| _{\Sigma  _{-4 {n}_0}}), \\
  &  \label{eq:fl132}  |\textbf{A}|    + |\textbf{B}| +|\gamma |  \le C  (|z|^2+\| f \| _{\Sigma  _{-4{n}_0} }^2 +|\varrho  | )  .
\end{eqnarray}
 We have $\mathcal{Z}= \mathcal{Z}(z,f,\varrho (f)) $,  $\textbf{A}= \textbf{A}(z,f,\varrho (f)) $,   $\textbf{B}= \textbf{B}(z,f,\varrho (f))$ and $\gamma= \gamma(z,f,\varrho (f))$,
with $C^n$ dependence in $z\in \C ^m$, $f\in \Sigma _{-4n_3}$ and $\varrho \in \R ^4$.    $ \mathcal{G}= \mathcal{G}(z,f,\varrho (f)) $  has   $C^n$ dependence in $z\in \C ^m$, $f\in \Sigma _{-4n_3}$ and $\varrho \in \R ^4$
with values in $\Sigma  _{ 4{n}_0-n}$. Here  $n_0=n_3+2n$ and $4n_3\ge n+1$.

\end{lemma}
\proof The formulas follow  by the use of the integrating factor. The last claim follows from Proposition \ref{prop:quasilin2}. We now consider
the estimates \eqref{eq:fl13}--\eqref{eq:fl132}, which improve \eqref{eq:fl12}. Notice that \begin{equation} \label{eq:fl14}
\begin{aligned} &     z^t  =   z  +
 \mathcal{Z}_t  \, \, \quad
  &  f^t(x) =  e^{   \sigma _3  ( \textbf{B}_t\cdot x +  {\gamma}_t ) }\tau_{\mathbf{A}_t}f
    +\mathcal{G} _t(x),
\end{aligned}
\end{equation}
 defined similarly to \eqref{eq:quasilin51}
  but with integrals in $[0,t]$ (resp. $[s,t]$) rather than in $[0,1]$  (resp. $[s,1]$).  We have by \eqref{eq::Stype1}--\eqref{eq:fl12}
 \begin{equation}
\begin{aligned} &     |
 \mathcal{Z}_t | \le \int _0^t|Z(\tau )|d\tau \le C (|z |+\| f \| _{\Sigma _{-4n_0}}+ |\varrho (f )|) \int _0^t(|z ^{\tau}|+\| f^{\tau} \| _{\Sigma _{-4n_0}}) d\tau .
\end{aligned}\nonumber
\end{equation}
Similarly
\begin{equation}\label{eq:fl141}
\begin{aligned} &       \int _0^t|\mathcal{C}(\tau )|d\tau \le C (|z |+\| f \| _{\Sigma _{-4n_0}}+ |\varrho (f )|) \int _0^t(|z ^{\tau}|+\| f^{\tau} \| _{\Sigma _{-4n_0}}) d\tau .
\end{aligned}
\end{equation}
For $\Upsilon =\mathcal{A},\mathcal{B}$ we have
 $ \int _{s}^{t}|
  \Upsilon (\tau ) | d\tau\le C (|z |+\| f \| _{\Sigma _{-4n_0}}+ |\varrho (f )|)  ^2 $.
  As a consequence we  get  $ \|  e^{   \sigma _3  ( \textbf{B}_t\cdot x +  {\gamma}_t ) }\tau_{\mathbf{A}_t}f \| _{\Sigma _{-4n_0}}\le C  \| f\| _{\Sigma _{-4n_0}}$ and
  \begin{equation}
\begin{aligned} &    \| \mathcal{G} _t  \| _{\Sigma _{ 4n_0}}\le C
\int _0^t  \| \mathcal{D}  (s, x+\int _s^t\mathcal{A}(\tau )d\tau   )  \| _{\Sigma _{ 4n_0}} ds \le C' \int _0^t  \| \mathcal{D}  (s, x    )  \| _{\Sigma _{-4n_0}} ds
\\& \le C ''   (|z |+\| f \| _{\Sigma _{-4n_0}}+ |\varrho (f )|) \int _0^t(|z ^{s}|+\| f^{s} \| _{\Sigma _{-4n_0}}) ds .
\end{aligned}\nonumber
\end{equation}
Then  \eqref{eq:fl13} follows by   Gronwall inequality and implies \eqref{eq:fl131}--\eqref{eq:fl132}.
\qed

\begin{remark}
\label{rem:lossreg} Notice  that the theory in Sect \ref{sec:ODE} entails
loss of regularity, in the sense that   $\mathbf{A}$ and the other functions are regular in
$f\in \Sigma _{-4 n_3}$ and not $f\in \Sigma _{-4 n_0}$. Since we consider many flows, we have big losses of regularity. Fortunately
we consider no more that $2N+2$ transformations  and we have
a lot of regularity to begin with.
\end{remark}

\begin{lemma}\label{lem:quasilin6}
Consider the system $\dot f =(\mathcal{X}^t)_f$ and $\dot z_j =(\mathcal{X}^t)_j$. Then the conclusions of Lemma \ref{lem:quasilin5}
continue to hold and we have also
\begin{equation} \label{eq:quasilin511}
\begin{aligned} &    f'(x) =  e^{  \sigma _3  (  - \frac{\im}{2} v'\cdot (x-\mathbf{A} )   +  \widetilde{\gamma} ) } f(x-\mathbf{A} )
    +\mathcal{G} (x)
\end{aligned}
\end{equation}
with $v'$ the velocity associated to the $t=1$ vector.
 \end{lemma}
\proof
The starting point is formula \eqref{eq:quasilin51}.
By Lemma \ref{lem:vectorfield} we have
\begin{equation} \mathbf{B}=\int_0^1\mathcal{B}_a(t )dt = -\frac{\im}{2}\int_0^1 ( v_a(t ) +t  dv_a(\mathcal{X}^{t})  dt=-\frac{\im}{2} v_a(1 ).\nonumber
\end{equation}
Recalling the $\gamma =\int _0^1 \left (\mathcal{C}(s)+ \mathcal{A}_a(s)\int _0^s \mathcal{B} _a(\tau ) d\tau  \right ) ds$  in Lemma \ref{lem:quasilin5},
we have

\begin{equation}
\begin{aligned} &    \int _0^1 ds  \mathcal{A}_a(s)\int _0^s \mathcal{B} _a(\tau ) d\tau    =-\frac{\im}{2}\int _0^1   \mathcal{A}_a(s)v_a(s )\, ds=\\& \frac{\im}{2}\int _0^1   v_a(s )\frac{d}{ds}\mathbf{A}_a (s)\, ds=\frac{\im}{2}   \mathbf{A}_a (1)v_a(1 )-\frac{\im}{2}\int _0^1    \mathbf{A}_a (s)\frac{d}{ds}v_a (s )\, ds.
\end{aligned}\nonumber
\end{equation}
We get \eqref{eq:quasilin511}   setting
$\widetilde{ \gamma}:= \gamma -\frac{\im}{2}   \mathbf{A}\cdot v '$.

\qed

\begin{lemma}\label{lem:quasilin7} $z'$ is $C^n$   in $z\in \C ^m$, $f\in \Sigma _{-4n_3}$ and $\varrho \in \R ^4$. $f'$ is $C^n$   in $z\in \C ^m$, $f\in \Sigma _{-4n_3}$ and $\varrho \in \R ^4$ with values in $  \Sigma _{-4n_3-n}$,   where $4n_3\ge n+1$.  $\varrho (f')$ is $C^n$   in $z\in \C ^m$, $f\in \Sigma _{-4n_3}$
and $\varrho \in \R ^4$.

\end{lemma}
\proof With the notation of Lemmas \ref{lem:quasilin5} and \ref{lem:quasilin6},
 for $\mathcal{Z}$, $\mathbf{A}$, $\widetilde{\gamma}=\gamma -\frac{\im}{2} \mathbf{A}\cdot v'$ and  $\mathcal{G}$ we have the result of   Lemma  \ref{lem:quasilin5}. So for $z'$ the claim follows immediately while for  $f'$ is a consequence of formula \eqref{eq:quasilin511}, the chain rule and Lemma \ref{lem:sigma3}.
We have
\begin{equation}  \label{eq:Hcoeff1}  \begin{aligned}
 &   Q(f')=  Q(f )+2\langle  e^{  \sigma _3  (  - \frac{\im}{2} {v'\cdot x}   +  \widetilde{\gamma} ) } f|\sigma _1\tau _{-\mathbf{A}} \mathcal{G} \rangle + Q( \mathcal{G} ) \\&  \Pi _a (f')= \Pi _a(f )- \frac{v_ a'}{2} Q(f )+\im \langle e^{  \sigma _3  (  - \frac{\im}{2} {v'\cdot x}   +  \widetilde{\gamma} ) } f|\tau _{-\mathbf{A}}  \sigma _1\sigma _3
 \partial _a\mathcal{G} \rangle +  \Pi _a ( \mathcal{G} ).
\end{aligned}
\end{equation}
By Lemma \ref{lem:quasilin5} we have that $\sigma _1\sigma _3^ie^{  \sigma _3  (    \frac{\im}{2} {v'\cdot x}  -  \widetilde{\gamma} ) } \tau _{-\mathbf{A}}\partial _a^i \mathcal{G}$ for $i=0,1$  is   $C^n$   in $z\in \C ^m$, $f\in \Sigma _{-4n_3}$ and $\varrho \in \R ^4$
with values in $\Sigma  _{ 4{n}_0-n-i}$ with $n_0=n_3+2n$. Then for  $f\in \Sigma _{-4n_3}$ and by $4{n}_0-n-i\ge 4n_3$   it follows that the mixed terms are $C^n$
in $z\in \C ^m$, $f\in \Sigma _{-4n_3}$
and $\varrho \in \R ^4$.   $Q( \mathcal{G} )$  and $\Pi _a ( \mathcal{G} )$ are   of  the desired type.
\qed

\section{Reformulation of \eqref{eq:SystK} in the new coordinates}
\label{section:reformulation}
 Denote by  $ \mathcal{F}_t$ the flow of the system \eqref{eq:quasilin10} associated to the field of
  Lemma \ref{lem:vectorfield}. We set
\begin{equation}    \label{eq:newH}  \begin{aligned} &
H=K\circ \mathcal{F}_1.
\end{aligned}
\end{equation}
In the new coordinates \eqref{eq:SystK}  becomes

\begin{equation} \label{eq:SystK2} \begin{aligned} &
 \im \dot z_j  =   \frac{\partial H}{\partial   \overline{z}_j  }
\, , \quad    \im \dot f=    \sigma _3    \sigma _1 \nabla _f  H. \end{aligned}
\end{equation}
   For system
\eqref{eq:SystK2} we prove:
\begin{theorem}\label{theorem-1.2}
  There exist $\varepsilon _0 >0$
and $C>0$ such that for   $ |z(0)|+\| f (0) \| _{H^1 }\le \epsilon
<\varepsilon _0 $  with $f(0)$ in the Schwartz class,  then the corresponding  solution of \eqref{eq:SystK2} is
globally defined and there are    $f_\pm  \in H^1$ with $\| f_\pm\|
_{H^1 }\le C \epsilon $ and functions   $\widehat{\vartheta}
\in C^1(\R , \R )$ and
  $\widehat{D}\in C^1(\R , \R ^3 )$
such that
\begin{equation}\label{scattering1}\lim _{t\to \pm \infty }
\left \|  \tau _{\widehat{D}(t)}e^{\im \widehat{\vartheta} (t) \sigma _3}f (t) -
 e^{  \im t \Delta \sigma _3  }  {f}_\pm   \right
\|_{H^{1}}=0
.
\end{equation}
  We   have

  \begin{equation}\label{decay}
  \lim _{t\to   \infty } z(t)=0.
  \end{equation}
   It is possible to write
$f(t,x)=A(t,x)+\widetilde{f}(t,x)$ with $|A(t,x)|\le C_N(t) \langle
x \rangle ^{-N}$ for any $N$, with $\lim _{t\to \infty }C_N(t)=0$
and such that for any admissible  pair $(p,q)$, i.e.
\eqref{admissiblepair}, we have
$
 \|  \widetilde{f} \| _{L^p_t(  \mathbb{R},
W^{ 1,q}_x)}\le
 C\epsilon .$

\end{theorem}

We now move to the proof of Theorem \ref{theorem-1.2}.
First of all we remark that in the sequel we need that the hamiltonians and
the Lie transforms be sufficiently regular. The amount of regularity
needed depends on $N=N_1$.  We will consider  a total of $2(N+1)$
Lie transforms, considering both the implementation of Darboux theorem
and the  Birkhoff normal
forms. We need to end up with a final hamiltonian which is
at least $C^1$. We can make sure that all the
  hamiltonians and  Lie transforms are sufficiently regular
  by picking $H^{K,S}$ with $K\gg 2N$ and $S\gg 2N$
  in Lemmas
\ref{lem:ExpK} and  \ref{lem:vectorfield}.

The first step in the proof of Theorem \ref{theorem-1.2}.
is a preliminary
discussion    of $H=K\circ  \mathcal{{F}}_1$. In   Sections \ref{sec:Normal form}--\ref{sec:Canonical} we implement the method of Birkhoff normal
forms, looking for other coordinates. Finally we will settle in the right system of coordinates and in Sect. \ref{sec:dispersion} we will finally prove
estimates.
\begin{lemma}
  \label{lem:ExpH}  Fix a large number $M\in \mathbb{N}$  with $M\gg  2N$.
  We have  the
expansion
\begin{equation}  \label{eq:ExpH1} \begin{aligned} & H =  \psi (\varrho (f)) +H_2 ^{(1)}+\resto ^{(1)}
\end{aligned}
\end{equation}  where  $\ {{\psi}} (\varrho )$  is $C^{ M}$  in $\varrho $  and where:
\begin{itemize}
\item[(1)]
 We have for $\ell =1$

\begin{equation}  \label{eq:ExpH2} H_2 ^{(\ell )}=
\sum _{\substack{ |\mu +\nu |=2\\
\lambda ^0\cdot (\mu -\nu )=0}}
 a_{\mu \nu}^{(\ell )}( \varrho (f) )  z^\mu
\overline{z}^\nu + \frac{1}{2} \langle \sigma _3 \mathcal{H}_{\omega
_0} f| \sigma _1 f\rangle .
\end{equation}

\item[(2)] We have $\resto ^{(1)}=\widetilde{\resto ^{(1)}} +
 \widetilde{\resto ^{(2)}}  $, with $\widetilde{\resto ^{(1)}}=$
\begin{equation}  \label{eq:ExpH2resto} \begin{aligned} &
 =\sum _{\substack{ |\mu +\nu |=2\\
\lambda ^0\cdot (\mu -\nu )\neq 0  }} a_{\mu \nu }^{(1)}(\varrho (f)
)z^\mu \overline{z}^\nu +\sum _{|\mu +\nu |  = 1} z^\mu
\overline{z}^\nu \langle \sigma _1 \sigma _3G_{\mu \nu }(\varrho (f)
)|f\rangle \end{aligned}
\end{equation}

    \begin{equation}    \begin{aligned} &  \widetilde{\resto ^{(2)}}=   \sum _{|\mu +\nu |= 3} z^\mu \overline{z}^\nu  a_{\mu \nu
}( z,f,\varrho (f) )     +\sum _{|\mu +\nu
|= 2} z^\mu \overline{z}^\nu \langle  G_{\mu \nu }( z,f,\varrho (f)
)|\sigma _3\sigma _1 f\rangle    \\& +   \sum _{d=2}^4
\langle B_{d } (   z ,f,\varrho (f) )|    f  ^{   d} \rangle
      +\int _{\mathbb{R}^3}
B_5 (x,  z ,f, f(x),\varrho (f) )  f^{   5}(x) dx\\& + \widehat{\resto}^{(1 )}_2(   z ,f,\varrho (f) )+   E_P (  f).
\end{aligned}\nonumber
\end{equation}
with $B_2(0,0,0)=0$ and where, both here    and in   Theorem \ref{th:main} later, by $f^d(x)$ we schematically represent $d-$products of components of $f$.

\item[(3)] At $\varrho (f)=0$ with $\ell =1$
\begin{equation}  \label{eq:ExpHcoeff1} \begin{aligned} &
a_{\mu \nu }^{(\ell )}( 0 ) =0 \text{ for $|\mu +\nu | = 2$  with $(\mu
, \nu )\neq (\delta _j, \delta _j)$ for all $j$,} \\& a_{\delta _j
\delta _j }^{(\ell )}( 0 ) =\lambda _j (\omega _0)  , \text{ where
$\delta _j=( \delta _{1j}, ..., \delta _{mj}),$}
\\& G_{\mu \nu }(  0 ) =0 \text{ for $|\mu +\nu | = 1$ }
\end{aligned}
\end{equation}
These $a_{\mu \nu }^{(\ell )}( \varrho )$ and $G_{\mu \nu }( x,\varrho
)$ are $C ^{ M }$ in all variables with $G_{\mu \nu }( \cdot ,\varrho )
\in C^M ( \mathbb{R}^4, \Sigma _{4M  } (\mathbb{R}^3,\mathbb{C}^2))$.

\item[(4)] For   a
small neighborhood $\mathrm{U}$ of $( 0,0,0)$ in
$ \mathbb{C}^m\times \Sigma _{-4M  }\times \R ^4$, we have
$a_{\mu \nu
}(   z,  \varrho   ) \in C^{ M }( \mathrm{U},
 \mathbb{C} ) $ .

\item[(5)] $G_{\mu \nu
}( \cdot ,  z,  \varrho ) \in C^{ M } ( \mathrm{U},
\Sigma _{4M }(\mathbb{R}^3,\mathbb{C}^2)) $.

\item[(6)] $B_{d
}( \cdot ,  z,f,\varrho  ) \in C^{ M } ( \mathrm{U},
\Sigma _{4M }(\mathbb{R}^3, B   (
 (\mathbb{C}^2)^{\otimes d},\mathbb{C} ))) $, for $2\le d \le 4$.

\item[(7)] Let $^t\eta = (\zeta , \overline{\zeta}) $ for
$ \zeta \in \C$. Then for
  $B_5(\cdot ,  z ,f, \eta ,\varrho )$   we have
\begin{equation} \label{H5power2}\begin{aligned} &\text{for  $|l|\le  M $ ,
}  \| \nabla _{  z,\overline{z} ,f,\zeta,\overline{\zeta} ,\varrho  }
^lB_5(  z,f,\eta ,\varrho  ) \| _{\Sigma _{4M } (\mathbb{R}^3,   B   (
 (\mathbb{C}^2)^{\otimes 5},\mathbb{C} )} \le C_l .
 \end{aligned}\nonumber \end{equation}

\item[(8)] We have for all indexes and for $\ell =1$
\begin{equation}  \label{eq:ExpHcoeff2} \begin{aligned} & a_{\mu \nu }^{(\ell )} =
\overline{a}_{\nu \mu   }^{(\ell )}\, , \quad a_{\mu \nu }
=\overline{a}_{\nu \mu   }\, , \quad  G_{\mu \nu } =-\sigma
_1\overline{G}_{\nu \mu } .
\end{aligned}
\end{equation}
\item[(9)]
\begin{equation}\label{eq:Rhat0}\begin{aligned} &
\widehat{\resto}^{(1 )}_2
 \in C^{M} ( \mathrm{U} ,\R ),  \\&    | \widehat{\resto}^{(1 )}_2 (z ,f, \varrho
)|  \le C (|z|+|\varrho |+ \| f \| _{ \Sigma _{-4M}}) \| f \| _{
\Sigma _{-4M}}^2;
\end{aligned}\end{equation}
 \end{itemize}
\end{lemma}
\proof We consider the notation of Lemma \ref{lem:quasilin6}.
Thanks to Lemma  \ref{lem:quasilin5} and by the freedom of choice of
$H^{K,S}$
  in Lemmas
\ref{lem:ExpK} and  \ref{lem:vectorfield}, we can assume that $\mathbf{A} $ and $v'
\in C ^{\widetilde{M}}(\widetilde{\mathrm{U}}, \R^3)$, $\mathcal{G}\in C ^{\widetilde{M}}( \widetilde{{\mathrm{U}}} , \Sigma _{4\widetilde{M}})$ and $\widetilde{\gamma}  \in C ^{\widetilde{M}}(\widetilde{\mathrm{U}},  \im \R )$,  with $\widetilde{\mathrm{U}}$
 a neighborhood of the origin in the space
 $\varrho (f)=\varrho \in \R ^4$, $z\in \C ^m$ and $f\in \Sigma _{-4\widetilde{M}} $ and with $\widetilde{M}\gg M  $ .
Having in mind \eqref{eq:Hcoeff1},
we notice that $e^{\sigma _3  (  - \frac{\im}{2} {v'\cdot x}   +  \widetilde{\gamma} ) }
\partial ^i_a\tau _{-\mathbf{A}} \mathcal{G}\in   C ^{ \widehat{M}}( \widetilde{{\mathrm{U}}} , \Sigma _{4\widetilde{M}-4\widehat{M}-i})$ for $i=0,1$.
Set $F:=   e^{  \sigma _3  (  - \frac{\im}{2} {v'\cdot x}   +  \widetilde{\gamma} ) } \tau _{-\mathbf{A}} \mathcal{G}  $. Then we have, for $\varrho =\varrho (f)$
and $G(t) =F (tz,tf,\varrho  )$,
\begin{equation}  \label{eq:Hcoeff11} \begin{aligned}
 &  \langle f|F\rangle =  \langle f|  (z_j\partial _j +\overline{z}_j\partial _{\overline{j}}  ) F (0,0,\varrho  ) \rangle \\&  +\langle f| \partial_fF (0,0,\varrho  ) f\rangle    + \frac{1}{2}\int _0^1\langle f|   \frac{d^2}{dt^2} G (t   )     \rangle  dt .
\end{aligned}
\end{equation}
The first term in the right  is like the $ z^\mu
\overline{z}^\nu \langle \sigma _1 \sigma _3G_{\mu \nu }(\varrho (f)
)|f\rangle $ in \eqref{eq:ExpH2resto} with $G_{\mu \nu }(\varrho
) \in C ^{ \widehat{M}-1}( \R ^4 , \Sigma _{4\widetilde{M}-4\widehat{M} })$.
By taking appropriate $\widetilde{M}$ and $\widehat{M}$, the conditions in the statement will hold.
The second term in the rhs of \eqref{eq:Hcoeff11} is like $\widehat{\resto}^{(1 )}_2 (z ,f, \varrho
)$, satisfying \eqref{eq:Rhat0}   for appropriate choices of $\widetilde{M}$ and $\widehat{M}$ by \eqref{eq:fl131}. The last term in \eqref{eq:Hcoeff11} is higher order,
and again is like  $\widehat{\resto}^{(1 )}_2$ or can be absorbed in $\widetilde{\resto}^{(2 )}$.
Similar expansions hold for  $Q( \mathcal{G} )$ and analogous
terms on the rhs of the second equality in \eqref{eq:Hcoeff1}.
  Expanding like in \eqref{eq:Hcoeff1} we obtain
\begin{equation}  \label{eq:Hcoeff2}  \begin{aligned}
 & (z')^\mu
(\overline{z} ')^\nu \langle \sigma _1 \sigma _3\underline{G}_{\mu \nu }(\varrho (f')
)|f'\rangle =\\&
(z+Z)^\mu
(\overline{z} +\overline{Z})^\nu \langle \sigma _1 \sigma _3\underline{G}_{\mu \nu }(\varrho (f')
)| \tau _{ \mathbf{A}}e^{  \sigma _3  (  - \frac{\im}{2} {v'\cdot x}   +  \widetilde{\gamma} ) }  f
   +\mathcal{G}\rangle .
\end{aligned}
\end{equation}
 Analogous formulas hold for $\langle \underline{B}_{d } (   z' ,\varrho (f') )|    (f')  ^{   d} \rangle $ for $d=2,3,4$ and $E_P (  f') $. With Taylor expansions
similar to   \eqref{eq:Hcoeff11}  we get in an elementary fashion that
\eqref{eq:Hcoeff2}  expands into  terms
falling in one of the cases in the statement. For   $d=5$,   $0\le j\le 5$  and with some exponentials in absorbed  $\underline{B}$,   schematically,  we have  terms like
\begin{equation}  \label{eq:Hcoeff3}  \begin{aligned}
 & \int _{\R ^3}\underline{\widehat{ {B}}}_5(x+\mathbf{A}  ,  z' ,e^{  \sigma _3  (  - \frac{\im}{2} {v'\cdot x}   +  \widetilde{\gamma} ) }
  f(x )+  \mathcal{G} (x+\mathbf{A}) ,\varrho (f') )  f^j(x)\mathcal{G}^{5-j}(x+\mathbf{A} )  dx
\end{aligned}\nonumber
\end{equation}
 which by Taylor expansion and by the fact that $\widetilde{M}\gg  M$ can be absorbed in $\widetilde{\resto ^{(2)}}$. We next look at
\begin{equation}  \label{eq:Hcoeff31}  \begin{aligned}
 & \langle \sigma _3 \mathcal{H}_{\omega
_0} f'| \sigma _1 f'\rangle = \langle (-\Delta +\omega _0) f ' | \sigma _1 f'\rangle + \langle \sigma _3 V_{\omega
_0} f'| \sigma _1 f'\rangle  .  \end{aligned}\nonumber
\end{equation}
We have
\begin{equation}  \label{eq:Hcoeff5}  \begin{aligned}
 &   \langle (-\Delta +\omega _0) f ' |  \sigma _1 f'\rangle =  \langle (-\Delta +\omega _0) f   | \sigma _1 f\rangle  +\frac{v^2}{2}Q(f)  + 2v_a \Pi _a(f)\\& + \langle (-\Delta +\omega _0) \mathcal{G}| \sigma _1 \mathcal{G}\rangle   + 2\langle (-\Delta +\omega _0) \mathcal{G}| \sigma _1\tau_{\mathbf{A}}e^{  \sigma _3  (   - \frac{\im}{2} {v'\cdot x}  +  \widetilde{\gamma} ) }f\rangle .\end{aligned}
\end{equation}
In \eqref{eq:Hcoeff5} the last line  can be treated as above and absorbed in the $\resto ^{(1)}$ while the last two terms of the first line go in part in $\psi (\varrho (f))$ and in part in $\resto ^{(1)}$ by $v=2\Pi (R)/Q(U)$.
The term $\langle \sigma _3 V_{\omega
_0} f', \sigma _1 f'\rangle $ can be expanded as the sum of $\langle \sigma _3 V_{\omega
_0} f , \sigma _1 f \rangle $ plus a reminder term treating it as
\eqref{eq:Hcoeff2}.

\qed

\section{Normal forms and homological equation}
\label{sec:Normal form}

  We set $\mathcal{H}=\mathcal{H }_{\omega _{0}}
P_c(\mathcal{H} _{\omega _{0}} )$. Consider   $\C $ valued functions $a_{\mu \nu}^{(\ell )}( \varrho   )$ such that  $a_{\nu \mu}^{(\ell )} \equiv \overline{a_{\mu \nu}^{(\ell )}}$. We assume that $a_{\mu \nu}^{(\ell )}\in C ^{ k_0} (U,\C )$ for $k_0\in \N$ a fixed number and $U$ a neighborhood of 0 in $\R ^4$.
Then we set

\begin{equation}  \label{eq:H2} H_2^{(\ell )}  (\varrho ):=
\sum _{\substack{ |\mu +\nu |=2\\
\lambda (0 ) \cdot (\mu -\nu )=0}}
 a_{\mu \nu}^{(\ell )}( \varrho   )  z^\mu
\overline{z}^\nu + \frac{1}{2} \langle \sigma _3 \mathcal{H}_{\omega
_0} f| \sigma _1 f\rangle .
\end{equation}

\begin{equation} \label{eq:lambda}
\lambda _j^{(\ell )} ( \varrho  ) :=   a
_{\delta _j\delta _j} ^{(\ell )}(\varrho   ), \quad \lambda ^{(\ell )}=
(\lambda _1^{(\ell )}, \cdots, \lambda _m^{(\ell )}).\end{equation}
We assume  $\lambda _j^{(\ell )}(0 ) =  \lambda _j (\omega _0 )   $ and   $a_{\mu \nu}^{(\ell )}( 0 ) =0 $ if $(\mu ,\nu)\neq (\delta _j,\delta _j)$ for
all $j$, with $\delta _j$ defined in \eqref{eq:ExpHcoeff1}.

\begin{definition}

\label{def:normal form} A function $Z(z,f, \varrho )$ is in normal form if it
is a sum
\begin{equation}
\label{e.12} Z=Z_0+Z_1
\end{equation}
where $Z_0$ and $Z_1$ are finite sums of the following type:
\begin{equation}
\label{e.12a}Z_1= \sum _{|\lambda ( \varrho ) \cdot(\nu-\mu)|>\omega _{0}}
z^\mu \overline{z}^\nu \langle  \sigma _1\sigma _3 G_{\mu \nu}( \varrho   )|f\rangle
\end{equation}
with $G_{\mu \nu}( x,\omega ,\varrho )\in  C^{k_0} (   U,\Sigma _{k_1})$ for  fixed $k_0\in \N$;
\begin{equation}
\label{e.12c}Z_0= \sum _{\lambda  (0)\cdot(\mu-\nu)=0} a_{\mu   \nu}
( \varrho  )z^\mu \overline{z}^\nu
\end{equation}
and $a_{\mu   \nu} ( \varrho  )\in  C^{k_0} (   U,
\mathbb{C})$. We will always assume the symmetries
\eqref{eq:ExpKcoeff2}. \qed\end{definition}

For $G=G(x)$, by elementary computations we have
 \begin{equation} \label{PoissBra2} \begin{aligned} &  \frac{1}{2} \{ \langle \sigma _3 \mathcal{H} f| \sigma _1 f\rangle , \langle \sigma _1\sigma
_3 G ,f\rangle \} = -\im \langle f |\sigma _1\sigma _3 \mathcal{H}G
\rangle  ,
\\&  \frac{1}{2}\{  \langle \sigma _3 \mathcal{H} f| \sigma _1 f\rangle, Q(f) \} = \im \langle \mathcal{H}f|
 \sigma _1 f   \rangle = -\im \langle
 \beta ^\prime (\phi ^2  )\phi ^2 \sigma _3 f|
  f   \rangle ,\\&  \frac{1}{2}\{ \langle \sigma _3 \mathcal{H} f| \sigma _1 f\rangle, \Pi _a(f) \} =   \langle \sigma _3\mathcal{H}f|
 \sigma _1 \partial _a f   \rangle = -\frac{1}{2}  \langle
  \sigma _3(\partial _a V_{\omega _0}) f|
  \sigma _1 f   \rangle     , \\& \{ f , Q(f) \}  =-\im P_c(\omega _0) \sigma _3 f ,  \quad  \{ f , \Pi _a (f) \}  =  P_c(\omega _0) \partial  _a f .
\end{aligned}
\end{equation}

 We now discuss the homological equations.  We start   by assuming that
 $ \varrho$ is an external parameter

\begin{lemma}
\label{lem:NLhom1} We consider $\chi =\chi (b,B)$ with
\begin{equation}
\label{eq:chi}\chi (b,B) =\sum _{|\mu +\nu |=M_0+1} b_{\mu\nu}   z^{\mu} \overline{z}^{\nu} + \sum _{|\mu +\nu |=M_0} z^{\mu}
\overline{z}^{\nu}
 \langle  \sigma _1\sigma _3  B_{\mu   \nu
}
  | f \rangle
\end{equation}
for $b_{\mu\nu} \in \C $ and $B_{\mu   \nu
} \in  \Sigma  _{2k_1}$ with $k_1\in \N$. Here we interpret
the polynomial  $\chi$ as a function with parameters $b=(b_{\mu\nu})$ and $B=(B_{\mu\nu})$. Denote by $X_{2k_1}$ the space of the pairs $(b,B)$.
Let us also consider given polynomials with  $K=K( \varrho  ) $  and $\widetilde{K}=K( \varrho ,b,B  ) $ where:

\begin{equation}\label{eq:Krho}
K(\varrho ) :=\sum _{|\mu +\nu |=M_0+1} k_{\mu\nu} ( \varrho   ) z^{\mu} \overline{z}^{\nu} + \sum _{|\mu +\nu |=M_0} z^{\mu}
\overline{z}^{\nu}
 \langle  \sigma _1\sigma _3  K_{\mu   \nu
}(\varrho  )
  | f \rangle ,
\end{equation}
with $k_{\mu\nu} ( \varrho   ) \in C^ {k_0} ( U, \C )$ and $K_{\mu\nu} ( \varrho   ) \in C^{k_0} ( U,    \Sigma  _{2k_1}\cap P_c(\omega _0)L^2)$ for $U$ a neighborhood of 0 in $\R ^4$;    let

\begin{equation}\label{eq:tildeKrho} \begin{aligned} &\widetilde{K} ( \varrho ,b,B  ) :=
 \sum _{|\mu +\nu |=M_0+1} \widetilde{k}_{\mu\nu}   ( \varrho , b,B  ) z^{\mu} \overline{z}^{\nu} \\&+
\sum _{i=0}^1\sum _{a=1}^3
 \sum _{|\mu +\nu |=M_0} z^{\mu}
\overline{z}^{\nu}
 \langle  \sigma _1\sigma _3 \partial ^i_a  K_{a\mu   \nu
}^i(\varrho , b,B   )
  | f \rangle ,\end{aligned}
\end{equation}
with $\widetilde{k}_{\mu\nu}  \in C^{k_0} ( U \times X _{2k_1}, \C )$ and $\widetilde{K}_{a\mu\nu} ^i  \in C^{k_0} ( U\times X_{2k_1},   \Sigma _{2k_1}\cap P_c(\omega _0)L^2)$.
 Suppose also that the sums \eqref{eq:Krho} and \eqref{eq:tildeKrho} do not
 contain terms in normal form and that $\widetilde{K}( 0 ,b,B  )=0$. Then there exists
 a  neighborhood  $V\subseteq U$ of  0 in $\R ^4$  and a unique choice of functions
 $(b(\varrho  ), B(\varrho   )) \in C ^{k_0} (V,X _{2k_1} )$ such that for $\chi (\varrho  )
  =\chi (b(\varrho  ), B(\varrho  ))$,    $\widetilde{K} (\varrho  ) =\widetilde{K}  (\varrho ,b(\varrho  ),B(\varrho  ))$  we have

\begin{equation}
\label{eq:homologicalEq} \left\{ \chi (\varrho  ) ,H_2 (\varrho  )  \right \} = K   (\varrho)+\widetilde{K} (\varrho  ) + Z(\varrho    )
\end{equation}
where $Z(\varrho  )$ is in normal form and homogeneous of degree $M_0+1$ in $(z,\overline{z},f)$. If the coefficients of $K$ satisfy
the symmetries in \eqref{eq:ExpKcoeff2},
 the same is true for $\chi$.

\end{lemma}
\proof    Summing on repeated indexes,  we get \begin{equation} \label{eq:homologicalEq1}  \begin{aligned}&    \{
H_2, \chi  \} =    \im \lambda ( \varrho )\cdot (\mu -
\nu) z^{\mu} \overline{z}^{\nu} b_{\mu \nu } \\& +   \im \langle f |\sigma _1\sigma _3( \lambda ( \varrho )\cdot (\mu -
\nu)- \mathcal{H} )B _{\mu \nu }
\rangle   + \widehat{ K}  (\varrho   ,b,B)   ,
\end{aligned}
\end{equation}
where, with an abuse of notation,
  \begin{equation} \label{eq:homologicalEq2}  \begin{aligned}& \widehat{ K}  (\varrho ,  b,B) :=
\sum _{\mu \nu \mu ' \nu '}
 a_{\mu \nu}^{(\ell   )}( \varrho   ) \left ( b_{\mu ' \nu '}  +  \langle  \sigma _1\sigma _3  B_{\mu ' \nu '}
  | f \rangle\right )\{z^\mu
\overline{z}^\nu , z ^{\mu'}
\overline{z}^{ \nu '} \} ,
 \end{aligned}
\end{equation}
 with the sum only on $|\mu +\nu |=2$ with $(\mu , \nu )\neq (\delta _j , \delta _j )$ for all $j$.
 $\widehat{ K}  (\varrho  , b,B) $ is 0 for $  \varrho  =0$ and is a homogeneous polynomial of the same type of the above ones.   Denote by $\widehat{ Z}  (\varrho ,  b,B) $ the sum of its
monomials in normal form and set $\textbf{K } :=\widetilde  {K}  +\widehat{ K}  -\widehat{ Z}   $.
 We look at   \begin{equation} \label{eq:homologicalEq1}  \begin{aligned} &   \im \lambda (\varrho )\cdot (\mu -
\nu) z^{\mu} \overline{z}^{\nu} b_{\mu \nu }    +z^{\mu} \overline{z}^{\nu}  \im \langle f |\sigma _1\sigma _3( \lambda  (\varrho ) \cdot (\mu -
\nu)- \mathcal{H} )B _{\mu \nu }
\rangle    \\& +  \textbf{K }(\varrho   ,b ,B )  +K(\varrho  ) =0
\end{aligned}
\end{equation}
 that is  at

\begin{equation} \label{eq:NLhom11}\begin{aligned}&
k_{\mu \nu} (\varrho )  + \textbf{{k}}_{\mu \nu}(\varrho   ,b,B)  +\im b_{\mu \nu}  \lambda (\varrho )
\cdot (\mu - \nu)=0   \\&  {K}_{\mu \nu }(\varrho ) +  \textbf{{K}} _{\mu \nu}(\varrho   ,b,B)
  -  \im   (\mathcal{H}
-\lambda  (\varrho ) \cdot (\mu - \nu) )B_{\mu   \nu }    =0  ,
\end{aligned}
\end{equation}
with $\textbf{{k}}_{\mu \nu}$ and  $\textbf{{K}} _{\mu \nu}$ the coefficients of $\textbf{K}$.
Notice that by $ \textbf{k}_{\mu \nu}(0  ,b,B) =0$ and $\textbf{K}_{\mu \nu}(0 ,b,B) =0$,  for $\varrho = 0$ there is a unique solution $(b,B)\in X_{2k_1}$ given by
\begin{equation} \label{eq:NLhom12}\begin{aligned}&
  b_{\mu \nu}
 =\frac{ \im k_{\mu\nu}(0)}{
    \lambda  (0) \cdot(\mu-\nu)}   \, , \quad  B_{\mu   \nu } (0 )  =-\im R _{\mathcal{H}}      (\lambda (0) \cdot
 (\mu -\nu )  )  K_{\mu   \nu
}(0)  .
\end{aligned}
\end{equation}
Notice   that  for $i=0,1$ we have $R _{\mathcal{H}}      (\zeta  )P_c(\omega _0)\circ \partial _a^i\in B (\Sigma _{2l},\Sigma _{2l}) $ for any $l\in \Z$, $\zeta \not \in \sigma _e (\mathcal{H})$ and $a\in \{1,2,3\}$.
The case $i=0$ can be proved by induction over $l$ using material in Sect. \ref{sec:sigma}. The case $i=1$
holds for $\mathcal{H}_0:=\sigma _2(-\Delta +\omega )$. Finally, these facts and the resolvent identity
yield the case $i=1$ for $\mathcal{H}$.
Then Lemma \ref{lem:NLhom1} is a consequence of   the implicit function theorem.

\qed

Substituting $\varrho = \varrho (f)$ we obtain what follows:
 \begin{lemma}
\label{lem:NLhom2} Set  $ {K}_1 = {K}( \varrho  (f)   ) $,
$\widetilde{K}_1  = \widetilde{K} ( \varrho (f) , b(\varrho (f)) , B(\varrho (f)) ) $
    and  $\chi _1=\chi ( \varrho   (f) ) $. Then we have

\begin{equation}
\label{eq:homologicalEq1} \left\{ \chi _1, H_2  \right \} =  K_1+\widetilde{K}_1 +Z_1+L_1
\end{equation}
where $Z_1$ is in normal form and homogeneous of degree $M_0+1$ in $(z,\overline{z},f)$ and
 \begin{equation}  \label{eq:L0} \begin{aligned}  & L_1 =   \langle
  V_j(\varrho (f)) f|
  f   \rangle  \widetilde{ \chi}  _j  +   \langle
  T_j  f|
  f   \rangle  \widehat{ \chi}  _j ,
\end{aligned}
\end{equation}
where: $V_j(\varrho )\in C^{k_0-1}(U, \Sigma_{2k_1-1})$, $T_j\in B (\Sigma_{-2k_1 }, \Sigma_{l } )$ for all $l$; $ \widetilde{ \chi}  _j$ and $ \widehat{ \chi}  _j$ polynomials like $\chi _1$, with monomials of no smaller degree and with coefficients in $C^{k_0-1}$ in $\varrho$.
\end{lemma}
\proof By direct computation \eqref{eq:homologicalEq1} holds with
 \begin{equation}  \label{eq:L1} \begin{aligned}  & L_1 = \im \langle
 \beta ^\prime (\phi ^2  )\phi ^2 \sigma _3 f|
  f   \rangle  \partial _{Q(f)}\chi _1+\frac{1}{2}  \langle
  \sigma _3(\partial _a V_{\omega _0}) f|
  \sigma _1 f   \rangle \partial _{\Pi _a(f)}\chi _1+\\& z^\mu \overline{z}^\nu  \partial _{\varrho _i}a _{\mu  \nu}^{(\ell )}\left ( \partial _{\varrho _j}\chi  _{1}
  \{ \varrho _i(f),  \varrho _j(f) \} +z^{\mu '   } \overline{z}^{  \nu ' } \langle \sigma _1 \sigma _3 B _{\mu ' \nu ' }, \{ \varrho _i(f) ,f\} \rangle \right ) ,
\end{aligned}
\end{equation}
with $\varrho _0(f)=Q(f) $ and $\varrho _a(f)=\Pi _a(f) $ for $a=1,2,3.$
We have  \begin{equation}  \label{eq:L2}  \begin{aligned}
 &  \{ Q(f), \Pi _a(f) \} = \im \langle  \sigma _1f|P_c (\omega _0)\partial _af \rangle   =  \im \langle  \sigma _1f|  \partial _aP_d (\omega _0) f \rangle     ;
\\
 &
\langle  \sigma _1\sigma _3  B_{\mu   \nu
} |\{ Q(f), f \}  \rangle = \langle     B_{\mu   \nu
} |  P_c (\omega _0) \sigma _3f \rangle ;
\\
 &  \{ \Pi _b(f), \Pi _a(f) \} = \im  \langle  \sigma _1\sigma _3\partial _bf|P_c (\omega _0)\partial _af \rangle   =  - \im  \langle  \sigma _1\sigma _3\partial _bf|P_d (\omega _0)\partial _af \rangle;
\\
 &
\langle  \sigma _1\sigma _3  B_{\mu   \nu
} |\{ \Pi _a(f), f \}  \rangle = \langle     B_{\mu   \nu
} |  P_c (\omega _0) \partial _af \rangle  \\& = -\langle    \partial _a B_{\mu   \nu
} |    f \rangle - \langle     B_{\mu   \nu
} |  P_d (\omega _0) \partial _af \rangle .
\end{aligned}
\end{equation}
  \eqref{eq:L1}-- \eqref{eq:L2}   yield the properties of $L_1$.\qed

 \section{Canonical transformations}
\label{sec:LieTransf}
  We consider functions $\chi$

  \begin{equation}
\label{eq:chi1}\chi   =\sum _{|\mu +\nu |=M_0+1} b_{\mu\nu} (\varrho (f))  z^{\mu} \overline{z}^{\nu} + \sum _{|\mu +\nu |=M_0} z^{\mu}
\overline{z}^{\nu}
 \langle  \sigma _1\sigma _3  B_{\mu   \nu
}(\varrho (f))
  | f \rangle .
\end{equation}
 We assume   $ b_{\mu   \nu } \in C^{k_0 }(\mathbb{R}  ^4 , \mathbb{C})$ and $ B_{\mu   \nu }
\in C^{k_0 }(\mathbb{R} ^4,   \Sigma _{4k_1} (\mathbb{R}^3,
\mathbb{C}^2))$   satisfying the symmetries in
\eqref{eq:ExpKcoeff2}. Here $k_0$ and $k_1$ are fixed,   very large
and in $\N$ and $k_1\gg k_0$.
We will always assume $B_{\mu   \nu }
\in     L^2_c(\mathcal{H}_{\omega _0})$.
 We want to consider the  flow $\phi ^t$ associated to the
Hamiltonian vector field $X_{\chi}$ at time $t=1$ and use it to change coordinates.

Our first step consists in setting up the hamiltonian system associated to $\chi$.   It  is a quasilinear symmetric hyperbolic system.

\begin{lemma}\label{lem:qlin1chi}  Consider $\chi$ as in \eqref{eq:chi1}
satisfying the symmetries in
\eqref{eq:ExpKcoeff2}.
Summing on repeated indexes, the following holds.

\begin{equation}\label{eq:qlin1chi}
\begin{aligned} &
     \{f, \chi
\}  =\mathcal{ L} f  + \mathcal{D} \ ,   \quad    \{ z _j  , \chi  \}  =Z _j  \ ,    \quad \mathcal{ L} :=\mathcal{A }_a \partial _a  +  \mathcal{C}  \sigma _3 \end{aligned}
\end{equation}
where the coefficients are given by the following formulas:
\begin{equation}\label{eq:qlin1chi1}
\begin{aligned} & \mathcal{A }_a = \partial _{\Pi _{a} (f)}\chi
\quad , \quad  \mathcal{C }  =-\im \partial _{Q (f)}  \chi
\quad , \quad
Z_j= -\im \partial _{\overline{z}_j}\chi , \\& \mathcal{D}=-\im z^{\mu}
\overline{z}^{\nu}
    B_{\mu   \nu
}(\varrho (f) )-P_d(\omega _0)\mathcal{L} f.
      \end{aligned}
\end{equation}
The coefficients can be thought as dependent on $(z,f,\varrho (f))$.
If we substitute $\varrho (f)$ with an independent  variable $\varrho \in \R ^4$, then we have $ \mathcal{A }_a  \in C ^{k_0 -1} (\mathfrak{V}, \R)$,  $ \im \mathcal{B }_a , \im \mathcal{C} \in C ^{k_0-1 } (\mathfrak{V}, \R)$ and $   \mathcal{D }  \in C ^{k_0-1 } (\mathfrak{V}, \Sigma _{ 4k_1}) $,
with $\mathfrak{V}$ a neighborhood of the origin in  $\C ^m\times \Sigma _{8k_0-4k_1}\times \R^4$.

 The following inequalities hold:
\begin{equation}\label{eq:qlin1chi2}
\begin{aligned} & |\mathcal{A } |+|\mathcal{C } | \le C   |z|^{M_0 }(|z|+\| f \| _{\Sigma _{-4 k_1}})\\&   |Z |+\|\mathcal{D} \|_{\Sigma _{ 4k_1}} \le C  |z|^{M_0-1}(|z|+\| f \| _{\Sigma _{-4k_1}}) .
  \end{aligned}
\end{equation}

\end{lemma}
\proof  \eqref{eq:qlin1chi}--\eqref{eq:qlin1chi1} follow from a simple computation. The rest follows by Subsect. \ref{subsec:Ltrans} for $n_0=k_1$ and $n=k_0$.
\qed

 Lemma \ref{lem:qlin1chi} assures us that we are within the framework of
 Sec. \ref{sec:Quasilinear} and that  the Lie transform $\phi =\phi ^1$
 associated to the following system is well defined:

\begin{equation}\label{eq:systchi}
\begin{aligned} &
      \dot f= \{f, \chi
\}   \, ,\quad     \dot z   = \{ z  , \chi  \}    .  \end{aligned}
\end{equation}
In particular, we have:

\begin{lemma}\label{lem:systchi1} Suppose $k_1$ is sufficiently large.
Set     $(z^t,f^t)=\phi ^t (z,f)$.   Then we have
\begin{equation} \label{eq:systchi11}
\begin{aligned} &     z^t  =   z  +
 \mathcal{Z}_t  \,  , \quad
  & f^t  = e^{   \sigma _3     {\gamma}_t   }\tau_{\mathbf{A}_t}f
    +\mathcal{G}_t
\end{aligned}
\end{equation}
with $\tau_{\mathbf{A}_t}f
   (x)= f(x-\mathbf{A}_t)$, $\mathbf{A}_t=- \int  _0^t\mathcal{A}(\tau )d\tau ,$  $\mathcal{Z}_t=\int _0^tZ(\tau )d\tau  $,   $\gamma _t=\int _0^t  \mathcal{C}(s) ds$ and

\begin{eqnarray}\label{eq:G}
  &\mathcal{G}_t(x) =  \int _0^t e^{-\sigma _3  \int _s^t  \mathcal{C} (\tau )  d\tau} \mathcal{D}   (s, x+\int _s^t\mathcal{A }(\tau )d\tau   ) ds. \nonumber
\end{eqnarray}
We have  $\textbf{A}_t= \textbf{A}_t(z,f,\varrho (f)) $,
  and $\gamma _t= \gamma _t(z,f,\varrho (f))$,
with $C ^{k_0 -1} $ dependence in   $z\in \C ^m$, $f\in \Sigma _{8k_0-4k_1}$ and $\varrho \in \R ^4$. The same statement  holds for  $\mathcal{Z}_t=\mathcal{Z}_t(z,f,\varrho (f)) $ resp.
  $\mathcal{G}_t= \mathcal{G}_t(z,f,\varrho (f)) $
with values in  $\C^m$  resp. $\Sigma _{ 4k_1 -k_0 }$.

\noindent The $f^t$ has $C ^{k_0 -1} $ dependence in   $z\in \C ^m$, $f\in \Sigma _{8k_0-4k_1}$ and $\varrho \in \R ^4$
with values in   $\Sigma _{ 9k_0-4k_1 }$.

\noindent The $\varrho (f^t)$ has $C ^{k_0 -1} $ dependence in   $z\in \C ^m$, $f\in \Sigma _{8k_0-4k_1}$ and $\varrho  (f)\in \R ^4$
with values in   $\R^4 $.

\noindent There is a fixed constant  $C $ such that

\begin{eqnarray} \label{eq:systchi12}
  &    |\mathcal{Z}_t  |  +\norma{ \mathcal{G}_t  }
_{\Sigma _{4k_1}}  \leq C      |z | ^{M_0-1}
 ( |z|+ \norma{f }
_{\Sigma _{-4k_1}} ) ,\\& \label{eq:systchi120} |\mathbf{A}_t|  +|\gamma _t|    \leq C      |z | ^{M_0-1}
 ( |z|+ \norma{f }
_{\Sigma _{-4k_1}} )^2 .
\end{eqnarray}
\end{lemma}
\proof This is the analogue of   Lemmas \ref{lem:quasilin5} and \ref{lem:quasilin6} .
\qed

We will set $\phi =\phi ^1$, $(z',f')=\phi (z ,f )$ and we will drop the
subindex $t$ for $t=1$, that is  $\textbf{A} =\textbf{A}_t$ etc.

\begin{lemma}\label{lem:systchi2} In the above notation we have
\begin{eqnarray} \label{eq:systchi21}
  &     |Q( f') -Q( f ) |  \leq C
 |z|^{M_0-1}( |z|+ \norma{f }
_{\Sigma _{-4k_1} } )^{2},\\& \label{eq:systchi210} | \Pi (f')-\Pi (f)   |    \leq C
 |z|^{M_0-1}( |z|+ \norma{f }
_{\Sigma _{-4k_1+1}} )^{2}
\end{eqnarray}
for a fixed $C$ dependent on $\| (b _{\mu \nu},B _{\mu \nu}) \| _{C^1(\mathrm{U}_4)}$ with $\mathrm{U}_4\subset \R ^4$ a preassigned neighborhood of the origin.
\end{lemma}
\proof We have
\begin{equation}\label{eq:systchi22}  \begin{aligned}
 &   Q( f') =  Q( f )  + \int _0^1 \{ Q(f), \chi \} \circ \phi ^t dt  ,
\end{aligned}
\end{equation}
with
\begin{equation}  \begin{aligned}
 &  \{ Q(f), \chi \}= \{ Q(f), \Pi _a(f) \}  \partial   _{\Pi _a(f)}\chi + z^{\mu}
\overline{z}^{\nu}  \langle  \sigma _1\sigma _3  B_{\mu   \nu
}(\varrho (f))|\{ Q(f), f \}  \rangle .
\end{aligned}   \nonumber
\end{equation}
  We have  $|\{ Q(f), \chi \} |\le C (|b|+\| B\| _{\Sigma_{4k_1}})  |z|^{M_0-1}( |z|+ \norma{f }
_{\Sigma_{-4k_1}} )^{2} $ for a fixed $C$ dependent on $\| (b _{\mu \nu},B _{\mu \nu}) \| _{C^1}$
by the formulas in   \eqref{eq:L2}.
  The integral in \eqref{eq:systchi22} has the
same upper bound by Lemma \ref{lem:systchi1}, in particular by \eqref{eq:systchi11}
and inequalities \eqref{eq:systchi12}--\eqref{eq:systchi120}. This proves  \eqref{eq:systchi21}.
  For $ \Pi _a ( f')$ we have a similar argument  by  the following formulas and estimates: \begin{equation}   \begin{aligned}
 &  \{ \Pi _b(f), \Pi _a(f) \} = \im  \langle  \sigma _1\sigma _3\partial _bf|P_c (\omega _0)\partial _af \rangle   =  - \im  \langle  \sigma _1\sigma _3\partial _bf|P_d (\omega _0)\partial _af \rangle;
\end{aligned}  \nonumber
\end{equation}
$|\langle  \sigma _1\sigma _3\partial _bf|P_d (\omega _0)\partial _af \rangle|\le C_\ell \| f \| _{\Sigma _{-\ell}} ^2$ and  $|\langle     B_{\mu   \nu
} |  P_d (\omega _0) \partial _af \rangle |\le C _{\ell}\|B \| _{\Sigma _{-\ell}}   \| f \| _{\Sigma _{-\ell}}  $ for any $\ell$. We have
$|\langle    \partial _a B_{\mu   \nu
} |    f \rangle|\le C \|B \| _{\Sigma _{4k_1 }}    \| f \| _{\Sigma _{-4k_1 +1}}$ by  $\| \partial _aB \| _{\Sigma _{4k_1-1 }}\le c \|B \| _{\Sigma _{4k_1 }}$, see Lemma \ref{lem:sigma3}. Thanks to \eqref{eq:L2} and Lemma \ref{lem:systchi1} these inequalities yield \eqref{eq:systchi210}.
\qed

The information in Lemma \ref{lem:systchi1} is not sufficiently precise
for our purposes. Let us consider for any fixed $(z(0),f(0))$ the system

\begin{equation}\label{eq:systchi3}
\begin{aligned} &
      \dot g= -\im  \zeta ^{\mu} \zeta ^{\nu} B_{\mu \nu} (\varrho  (f(0)) \, , \\&     \dot \zeta _j   = -\im \nu _j \frac{ \zeta ^{\mu} \zeta ^{\nu}  }{\overline{\zeta} _j} \left ( b_{\mu \nu} (\varrho  (f(0))  +\langle  \sigma _1\sigma _3  B_{\mu   \nu
}(\varrho  (f(0))
  | f \rangle \right )   \\&  g(0) =f(0)  \, , \quad \zeta (0) =z(0) .  \end{aligned}
\end{equation}
Notice that well posedness, regularity of the flow, and smooth
dependence on the coefficients $ ( b  (\varrho  (f(0)) , B (\varrho  (f(0)))$
fall within the scope of the theory of ordinary equations. Denote $\phi _0^t$ the flow of \eqref{eq:systchi3}.
In particular  for $(\zeta ^t , g^t)=  \phi _0^t(\zeta  , g ) $ we have \begin{equation} \label{eq:comp}\begin{aligned} &\zeta ^t=\zeta  +  \mathbf{Z} _t( \zeta  ,g ,  b  (\varrho  (f ) ), B (\varrho  (f ) )) \\& g ^t=g + \mathbf{{G}}_t( \zeta  ,g ,  b  (\varrho  (f ) ), B (\varrho  (f )) ),
 \end{aligned}  \end{equation} with
$ \mathbf{{ Z}}_t( \zeta  ,g ,  b   , B   )$ (resp. $ \mathbf{{ G}}_t( \zeta  ,g ,  b   , B   )$)
with $C^{\infty}$ dependence on $t$,  $\zeta \in \C ^m$ and  $g\in \Sigma _{ -4k_1}$,
and     $( b    , B   ) $, with values in $\C ^m$ (resp.  $\Sigma _{4k_1}(\R ^3, \C ^2)$)
  and with
   \begin{equation}\label{eq:comp1}
    \mathbf{Z} _t( \zeta  ,g ,  0, 0 )\equiv 0, \, \quad \mathbf{G} _t( \zeta  ,g ,  0, 0 )\equiv 0.
   \end{equation}
  Furthermore $ \mathbf{{ Z}}_t  $ resp. $ \mathbf{{ G}}_t$ satisfy uniformly in $t$ the same bounds \eqref{eq:systchi12} of $\mathcal{ Z}   $ resp. $\mathcal{ G} $.

We   compare the solutions of \eqref{eq:systchi} with those of \eqref{eq:systchi3}. Denote   $(z',f')=\phi ^1(z,f)$ and $(\zeta ',g')=\phi _0 ^1(z,f)$.

\begin{lemma}\label{lem:systchi51}
  For a $C$ like in Lemma \ref{lem:systchi2} we have for any $j\le  k_1$
\begin{equation} \label{eq:systchi511}
\begin{aligned} &
      | z'- \zeta '|+\| f'- g '\|  _{\Sigma _{-4j-1}}\le C
 ( |z|+ \norma{f }
_{\Sigma _{ -4j +1 }} ) ^{M_0+1}  .
\end{aligned}
\end{equation}
  \end{lemma}
\proof
Set $\textbf{M}(t)=|z(t )|+  |\zeta(t )|+\norma{f(t ) }
_{\Sigma _{1-4j}}+\norma{g(t ) }
_{\Sigma _{1-4j }}.$
We have

\begin{equation}
\begin{aligned} & \dot f -\dot g= \im \zeta ^{\mu}
\overline{\zeta }^{\nu}
    B_{\mu   \nu
}(\varrho (f(0)) ) -\im z^{\mu}
\overline{z}^{\nu}
    B_{\mu   \nu
}(\varrho (f) ) -P_d(\omega _0)\mathcal{L} f + (\mathcal{A }_a \partial _a  +  \mathcal{C}  \sigma _3) f .
\end{aligned} \nonumber
\end{equation}
The rhs has  $ \Sigma _{-4j}$ norm bounded by
\begin{equation}
\begin{aligned} &   |\zeta ^{\mu}
\overline{\zeta }^{\nu} -z^{\mu}
\overline{z}^{\nu}| \, \|
    B_{\mu   \nu
}(\varrho (f(0)) )  \|  _{ \Sigma _{-4j} }
+|
z^{\mu}
\overline{z}^{\nu} | \|
    B_{\mu   \nu
}(\varrho (f ) ) - B_{\mu   \nu
}(\varrho (f(0)) ) \|  _{ \Sigma _{-4j} }
  \\&  + \| P_d(\omega _0)\mathcal{L} f \|  _{ \Sigma _{-4j} }+ \| (\mathcal{A }_a \partial _a  +  \mathcal{C}  \sigma _3) f\|  _{ \Sigma _{-4j} } .
\end{aligned} \nonumber
\end{equation}
Then \begin{equation}
\begin{aligned} &
      \| f(t)-g (t)\| _{\Sigma _{-4j }}\le  C \int _0^t
\textbf{M} ^{M_0-1}(t)   | z(\tau )-\zeta (\tau ) |  d\tau    + C \int _0^t
\textbf{M} ^{M_0+2}(t)     d\tau  .
\end{aligned}\nonumber
\end{equation}
Similarly

\begin{equation}
\begin{aligned} &
      | z(t)- \zeta (t)| \le  C \int _0^t
\textbf{M} ^{M_0-1}(t) (| z(\tau  )- \zeta (\tau )|  +\| f(\tau )-g (\tau )\| _{\Sigma _{-4j }} )  d\tau  \\& + C \int _0^t
\textbf{M} ^{M_0+1}(t)     d\tau  .
\end{aligned}\nonumber
\end{equation}
Then the statement follows from Gronwall inequality since $  \textbf{M}(t)
\le C   \textbf{M}(0)$. For instance, $|z(t )| +\norma{f(t ) }
_{\Sigma _{1-4j}}\le C (|z(0 )| +\norma{f(0 ) }
_{\Sigma _{1-4j}})$ follows  by formulas \eqref{eq:systchi11} and by inequalities
\eqref{eq:systchi12}--\eqref{eq:systchi120}.
\qed

\begin{lemma}\label{lem:systchi52} In the above notation of Lemma
\ref{lem:systchi52} we have

\begin{eqnarray} \label{eq:systchi521}
  &     |Q( f') -Q( g ') |  \leq C
  ( |z|+ \norma{f }
_{\Sigma _{-4k_1+1} } )^{M_0+2},\\& \label{eq:systchi522} | \Pi (f')-\Pi (g')   |    \leq C
  ( |z|+ \norma{f }
_{\Sigma _{-4k_1+3} } )^{M_0+2}.
\end{eqnarray}
\end{lemma}
\proof
We have

\begin{equation*} \begin{aligned} & \frac{d}{dt}( \Pi _b (f )-\Pi _b(g ))=
\im \zeta ^{\mu}
\overline{\zeta }^{\nu}
    \langle B_{\mu   \nu
}(\varrho (f(0)) ) | \sigma _3 \sigma _1\partial _bg\rangle \\& + \im z^{\mu}
\overline{z}^{\nu}
   \langle  B_{\mu   \nu
}(\varrho (f) )| \sigma _3 \sigma _1\partial _bf\rangle   -
\langle P_d(\omega _0)\mathcal{L} f|  \sigma _3 \sigma _1\partial _bf\rangle .
\end{aligned}
\end{equation*}
The right hand side can be bounded above by  the rhs of  \eqref{eq:systchi522}  computed at time $t$. Integrating, using the fact that  $(z,f)=(\zeta , g)$ at $t=0$ and by \eqref{eq:systchi11} and
\eqref{eq:systchi12}--\eqref{eq:systchi120} we get  \eqref{eq:systchi522}.
The proof of  \eqref{eq:systchi521} is similar.

\qed

 \section{Birkhoff normal forms}
\label{sec:Canonical} Our goal in this section is to prove the
following result.

\begin{theorem}
\label{th:main}   For any integer $2\le \ell \le  2N+1$   there are  a $\delta _0>0$ and   $M\gg N$ large
such that in the subset of $\Sigma  _{4M}$ defined by $|z|+\| f \| _{H^1}<\delta _0 $
is defined a canonical transformation $\Tr_r$ which is differentiable as a map
with values in  $\Sigma _1$  and
whose image contains a similar
subset of $\Sigma _{4M}$ defined by $|z|+\| f \| _{H^1}<\delta _0' $,  s.t.
\begin{equation}
\label{eq:bir1} H^{(\ell )}:= K\circ \Tr_\ell =  {\psi} (\varrho (f)) +H_2^{(\ell )}+Z^{(\ell )}+\resto^{(\ell )},
\end{equation}
     with ${\psi} (\varrho (f))$   the same of \eqref{eq:ExpH1} and where:
\begin{itemize}
\item[(i)] $H_2^{(\ell )}=H_2^{(2)}$ for $\ell \ge 2$, is of the
form \eqref{eq:ExpK2}
where  $a_{\mu \nu}^{(\ell )} $  satisfy
\eqref{eq:ExpHcoeff1}--\eqref{eq:ExpHcoeff2};

\item[(ii)]$Z^{(\ell )}$ is in normal form, with monomials of
degree $\le \ell$  whose coefficients satisfy \eqref{eq:ExpKcoeff2};
\item[(iii)] we have   $\Tr _\ell  =\phi _\ell \circ ...\circ \phi _1$,
 with each  $\phi _j$ a Lie   transformation  associated to a function
 \eqref{eq:chi1} with $M_0=j$;

\item[(iv)] we have $\resto^{(\ell )} = \sum _{d=0}^6\resto^{(\ell )}_d$
with the following properties (for $k_2(\ell )\ll k_3(\ell )\ll M$ pairs of appropriate large numbers  with $k_j(\ell +1) \ll k_j(\ell ) $ for all $\ell$ and for $j=2,3$):
\begin{itemize}
\item[(iv.0)]  we have  with $  |\partial _{  \varrho} ^l a_{\mu \nu }^{(\ell  )}(   \varrho ) |   \le C_l \text{ for
  $|l|\le  k_2(\ell )$},
$
\begin{equation} \resto^{(\ell )}_0=
\sum _{|\mu +\nu |  = \ell +1 } z^\mu \overline{z}^\nu  a_{\mu \nu }^{(\ell )}(   \varrho ( f) ):
\nonumber
\end{equation}

\item[(iv.1)]   we have with $  \|\partial _{  \varrho} ^l G_{\mu \nu }^{(\ell )}(   \varrho ) (\cdot )\| _{\Sigma _{ 4k_3(\ell )}(\R ^3, \C ^2)} \le C_l  \text{ for  $|l|\le k_2(\ell )$}$,
\begin{equation} \resto^{(\ell )}_1=
\sum _{|\mu +\nu |  = \ell   }z^\mu \overline{z}^\nu \langle  \sigma _1 \sigma _3G_{\mu \nu
}^{(\ell )}(   \varrho ( f) ) | f\rangle  \text{ with }   \nonumber
\end{equation}

\item[(iv.2--5)] in $\U  $  we have for $2\le d \le 5$ and for $\eta ^T=(\zeta , \overline{\zeta })$ with $\zeta \in \C$,

\begin{equation} \resto^{(\ell )}_d=
 \int
_{\mathbb{R}^3} F_d^{(\ell )}(x, z ,f, f(x),\varrho ( f)) f^d(x)   dx +  \widehat{\resto}^{(\ell )}_d,\nonumber
\end{equation}
 with   for
  $|l|\le  k_2(\ell )$

\begin{equation}\label{eq:coeff F}   \| \partial _{ z,\overline{z},
\zeta,\overline{\zeta},f,\varrho} ^l F_d^{(\ell )} (\cdot ,z,f,\eta ,
\varrho )\| _{\Sigma _{ 4k_3(\ell )} (\mathbb{R}^3,   B   (
 (\mathbb{C}^2)^{\otimes d},\mathbb{C} )} \le C_l  ,
\end{equation}
  with $F_2^{(\ell )}(x, 0 ,0,0,0)=0$
    and
  with $\widetilde{\resto}^{(\ell )}_d (z ,f, \varrho ( f))$ s.t.

\begin{equation}\label{eq:Rhat}\begin{aligned} &
\widehat{\resto}^{(\ell )}_d
 \in C^{ k_2(\ell )} ( \U   \times \R ,\R )\, , \quad  |
\widehat{\resto}^{(\ell )}_d (z ,f, \varrho )|\le C \| f \|
_{\Sigma _{-4k_3(\ell )}}^d , \\&    | \widehat{\resto}^{(\ell )}_2 (z ,f, \varrho
)|  \le C (|z|+|\varrho |+ \| f \| _{ \Sigma _{-4k_3(\ell )}}) \| f \| _{
\Sigma _{-4k_3(\ell )}}^2;
\end{aligned}\end{equation}

\item[(iv.6)] $ \resto^{(\ell )}_6=
 \int _{\mathbb{R}^3}  B ( | f (x)| ^2/2) dx$.

\end{itemize}
\end{itemize}
\end{theorem}

\subsection{Pullback of multilinear forms}
\label{subsec:pullback}

The method of Birkhoff normal forms is implemented using   the flows
of auxiliary hamiltonians
  $\chi$   like in \eqref{eq:chi1}. In particular, in we will assume for the moment that the degree is $M_0+1$ and that   $ b_{\mu   \nu } \in C^{k_0 }(\mathbb{R}  ^4 , \mathbb{C})$ and $ B_{\mu   \nu }
\in C^{k_0 }(\mathbb{R} ^4,   \Sigma _{4k_1} (\mathbb{R}^3,
\mathbb{C}^2))$. In the proof of Theorem \ref{th:main}, $\chi $ needs to solve
a \textit{homological equation}. In this section we look at pullbacks
of the various terms of the hamiltonian by means of the Lie transform
associated to $\chi$.  In general these terms are pulled back
 into other  reminder terms which are  less of regular. By this we mean
 both that their coefficients are less than $C^{k_0 }$ and   with values in some $\Sigma _{4k }$
 with $k<k_1$. In general this loss of regularity is harmless. However we have
 to make sure that the terms which enter in the homological equation of
  $\chi$, which is used to find  a useful $\chi$,
   have same regularity of $\chi$. It is at this stage that we use
the associated  simplified system \eqref{eq:systchi3}.
We will consider now a number of technical lemmas.

\begin{lemma}
  \label{lem:pull0} Let $F=F(  z,f,\varrho )$ be      $C ^{k_0} $  in  $z\in \C ^m$, $f\in \Sigma _{-4k_1}$ and $\varrho \in \R ^4$
with values in $\Sigma _{ 4k_1   }(\R ^3, B^2(\C , \C))$.
For $M_0=1$  we assume $F(0,0,0)=0$.
Then \begin{equation}  \label{eq:pull00}\begin{aligned} &  \langle F (   z ',f',\varrho (f') )|     \mathcal{G} ^2   \rangle =  \sum _{|\mu + \nu |=M_0 +1 }  {{k}}_{\mu \nu}  ( \varrho (f), b( \varrho (f)), B( \varrho (f)) )z^{\mu} \overline{z}^{\nu} \\&  +\sum _{|\mu + \nu  |=M_0 } z^{\mu} \overline{z}^{\nu}
 \langle  \sigma _1\sigma _3    {{K}}  _{\mu   \nu
}( \varrho (f), b( \varrho (f)), B( \varrho (f)) )
  | f \rangle  +    {\textbf{R}}  ,
\end{aligned}
\end{equation}
where the following holds.

\begin{itemize}

\item [(i)] ${k}_{\mu   \nu} (0,b,B)=0$ and   ${K}_{\mu   \nu} (0,b,B)=0$; ${k}_{\mu   \nu} (\varrho ,b,B) \in \C  $ and   ${K}_{\mu   \nu} (\varrho ,b,B)\in \Sigma _{4 k_1} $ are $C^{ k_0-1}$ in $\varrho$,
  in   $b _{\mu \nu}\in \C$
and   in $B _{\mu \nu}\in \Sigma _{4k_1}$.

\item [(ii)] ${\textbf{R}} $ is a sum of terms of the form $\resto  ^{(M_0+1)} $,   that is like  in the statement of Theorem \ref{th:main}, with $k_3= k_1-2k_0$ and $k_2= k_0-M_0-3$.
\end{itemize}
If $M_0>1$ formula \eqref{eq:pull00} holds with only ${\textbf{R}} $
in the rhs.

\end{lemma}

\proof  We have for $\varrho =\varrho (f ) $ and $\varrho '=\varrho (f ') $
and $\delta \varrho =\varrho '-\varrho   $

\begin{equation}  \label{eq:pull01}\begin{aligned} &  \langle F (   z ',f' ,\varrho' )|     \mathcal{G} ^2   \rangle =   \langle F (   z ',f' ,\varrho  )|     \mathcal{G} ^2   \rangle   +
\int _0^1\langle \partial_{\varrho  }F (   z ',f' ,\varrho   +t \delta \varrho )|     \mathcal{G} ^2   \rangle dt \cdot \delta \varrho.
\end{aligned}
\end{equation}
By Lemma \ref{lem:systchi1} the second term in the rhs is  $C ^{k_0-1} $  in  $z\in \C ^m$, $f\in \Sigma _{8k_0-4k_1}$ and $\varrho \in \R ^4$.
Furthermore $\delta \varrho$
satisfies $|\delta \varrho |\le C
 |z|^{M_0-1}( |z|+ \norma{f }
_{\Sigma _{-4k_1+1}} )^{2}$ by Lemma \ref{lem:systchi2}. By $k_0= k_2+M_0+3$  we can write
the second term in the rhs of  \eqref{eq:pull01} as  $\resto  _0^{(M_0+1)}+\resto  _1^{(M_0+1)}+\widehat{\resto}  _2^{(M_0+1)}$, just by performing an appropriate
and partial Taylor expansion.
If  $M_0>1$ the same result holds for the first  term in the rhs of  \eqref{eq:pull01}. Let now $M_0=1$.  By Lemmas  \ref{lem:systchi1}  and  \ref{lem:quasilin7} it is $C ^{k_0-1} $  in  $z\in \C ^m$, $f\in \Sigma _{8k_0-4k_1}$ and $\varrho \in \R ^4$ and can be expressed as
$\langle F (   0,0 ,\varrho  )|     \mathcal{G} ^2   \rangle $ plus a term
which is like the above ones and can be absorbed in ${\textbf{R}} $. Writing $F (    \varrho  )=F (   0,0 ,\varrho  )$, we have succinctly

\begin{equation}\label{eq:pull02} \begin{aligned}
 & \langle F (    \varrho  )|     \mathcal{G} ^2   \rangle =  \langle F (    \varrho  )|      \mathbf{{G}} _1^2   \rangle    -2 \langle F (    \varrho  )|      (\mathbf{{G}} _1 - \mathcal{G}) \mathbf{{G}} _1  \rangle + \langle F (    \varrho  )|      (\mathbf{{G}} _1 - \mathcal{G})^2   \rangle  ,\end{aligned}
\end{equation}
with $\mathbf{{G}} _1$ from the associated system, see \eqref{eq:comp}.
By Lemma \ref{lem:systchi51}   the rhs is  $C ^{k_0-1} $  in  $z\in \C ^m$, $f\in \Sigma _{8k_0-4k_1}$ and $\varrho \in \R ^4$.  We have
\begin{equation}\label{eq:pull03} \begin{aligned}
 &   \langle F (    \varrho  )|      (\mathbf{{G}} _1 - \mathcal{G}) \mathbf{{G}} _1  \rangle  = \langle F (    \varrho  )|      (g' - f') \mathbf{{G}} _1  \rangle  + \langle F (    \varrho  )|      (e^{\sigma _3\gamma}\tau_{\mathbf{A}}f - f ) \mathbf{{G}} _1  \rangle . \end{aligned} \nonumber
\end{equation}
Then
\begin{equation}\label{eq:pull04} \begin{aligned}
 &  | \langle F (    \varrho  )|      (g' - f') \mathbf{{G}} _1  \rangle |\le  \| F\| _{\Sigma _{ 4k_1 }} \| \mathbf{{G}} _1 \| _{\Sigma _{ 4k_1 }}  \| g' - f' \| _{\Sigma _{ 4k_1 }}  \\& \le C   \| F\| _{\Sigma _{ 4k_1 }} ( |z|+ \norma{f }
_{\Sigma _{1-4k_1 }} ) ^{3}. \end{aligned} \nonumber
\end{equation}
Similarly
\begin{equation}\label{eq:pull05} \begin{aligned}
 &  | \langle F (    \varrho  )|      (e^{\sigma _3\gamma}\tau_{\mathbf{A}}f - f ) \mathbf{{G}} _1   \rangle |\le  \| F\| _{\Sigma _{ 4k_1 }} \| \mathbf{{G}} _1 \| _{\Sigma _{ 4k_1 }}  \|  (e^{\sigma _3\gamma}\tau_{\mathbf{A}}f - f )  \| _{\Sigma _{ 4k_1 }}  \\& \le C   \| F\| _{\Sigma _{ 4k_1 }} ( |z|+ \norma{f }
_{\Sigma _{1-4k_1 }} ) ^{3}. \end{aligned} \nonumber
\end{equation}
Hence the second term in the rhs of \eqref{eq:pull02} can be absorbed in
$\mathbf{R}$. Similar reasoning applies to the third term in the rhs of \eqref{eq:pull02}.  We finally show that
the first term in the rhs of \eqref{eq:pull02} yields the first two terms
in the rhs of  \eqref{eq:pull00} plus a term which can be absorbed in
${\textbf{R}} $.  We know that $\mathbf{{G}} _1=\mathbf{{G}} _1(z,f,b_{\mu \nu} , B_{\mu \nu} )$ is  $C^{\infty}$  in  $z \in \C ^m$,  $f\in \Sigma _{ -4k_1}$
and     $( b_{\mu \nu} , B_{\mu \nu}   ) $, with values in    $\Sigma _{4k_1}(\R ^3, \C ^2)$ with   $\| \mathbf{{G}}_1\| _{\Sigma _{4k_1}} =O(|z|+\| f \|_{\Sigma _{ -4k_1}})$.  We can consider
\begin{equation}\label{eq:pull06} \begin{aligned}
 & \langle F (    \varrho  )|      \mathbf{{G}}_1 ^2   \rangle = \sum _{|\mu + \nu |=2}  h_{\mu \nu}  ( \varrho  , b , B  )z^{\mu} \overline{z}^{\nu}  \\&  +\sum _{|\mu + \nu  |=1 } z^{\mu} \overline{z}^{\nu}
 \langle  \sigma _1\sigma _3    H  _{\mu   \nu
}( \varrho  , b , B  )
  | f \rangle  +   \widetilde{ \resto} , \end{aligned}
\end{equation}
with  $ \widetilde{ \resto} = O(|z|+\| f \|_{\Sigma _{ -4k_1}})^3$ and $C^{k_0}$
in $z,f,\varrho$, and with
\begin{equation}\label{eq:pull07} \begin{aligned}
 & h_{\mu \nu}  ( \varrho  , b , B  ) := \frac{1}{\mu ! \nu !} \partial ^\mu _{z}  \partial ^\nu _{\overline{z}}\langle F (    \varrho  )|     \mathbf{{G}}_1 ^2   \rangle _{\mid (0,0,\varrho )} \ ,\\&  H_{\mu \nu}  ( \varrho  , b , B  ) :=  \sigma _3\sigma _1  \partial ^\mu _{z}  \partial ^\nu _{\overline{z}}\nabla _f\langle F (    \varrho  )|    \mathbf{{G}}_1 ^2   \rangle _{\mid (0,0,\varrho )}\  . \end{aligned} \nonumber
\end{equation}
$ \widetilde{ \resto}$ can be absorbed in $\mathbf{R}$.  The polynomial in
\eqref{eq:pull06} is like the one in the statement because of the hypothesis
$F( \varrho )=F(0,0,\varrho )=0$ when $\varrho =0$ if $M_0=1$.
\qed

\begin{lemma}
  \label{lem:pull1} Let $F=F(  z,f,\varrho )$ with the same properties  as in Lemma \ref{lem:pull0}.
   Then, for a  rhs which satisfies  the same properties
stated in  Lemma \ref{lem:pull0} but with $k_3= k_1-3k_0 $ and $k_2= k_0- M_0-4$, \begin{equation}  \label{eq:pull11}\begin{aligned} &  \langle F (   z ',f',\varrho (f') )|      \mathcal{G}  f '    \rangle =  \sum _{|\mu + \nu |=M_0 +1 }  {{k}}_{\mu \nu}  ( \varrho (f), b( \varrho (f)), B( \varrho (f)) )z^{\mu} \overline{z}^{\nu} \\&  +\sum _{|\mu + \nu  |=M_0 } z^{\mu} \overline{z}^{\nu}
 \langle  \sigma _1\sigma _3    {{K}}  _{\mu   \nu
}( \varrho (f), b( \varrho (f)), B( \varrho (f)) )
  | f \rangle  +    {\textbf{R}}  .
\end{aligned}
\end{equation}

 \end{lemma}
\proof We have
\begin{equation}  \label{eq:pull12}\begin{aligned} &  \langle F (   z ',f',\varrho (f') )|  \mathcal{G}    f '     \rangle =\langle F (   z ',f',\varrho (f') )|   \mathcal{G}    f  \rangle   +\\& \langle F (   z ',f',\varrho (f') )|   \mathcal{G} (e^{\sigma _3\gamma }\tau _{\mathcal{\mathbf{A}}} -1)  f  \rangle
 +\langle F (   z ',f',\varrho (f') )|   \mathcal{G} ^2  \rangle .
\end{aligned}
\end{equation}
The third  term in the rhs is like in Lemma \ref{lem:pull1}. The second
can be absorbed in ${\textbf{R}} $ by \eqref{eq:systchi120}. We focus on first  term in the rhs
of \eqref{eq:pull12}. By   \eqref{eq:G} we have

\begin{equation}  \label{eq:pull13}\begin{aligned} &
   \mathcal{G} (x) =  \int _0^1  \mathcal{D}   (s, x  ) ds \\& +
  \int _0^1 \left (e^{-\sigma _3  \int _s^t  \mathcal{C} (\tau )  d\tau} -1\right )\mathcal{D}   (s, x+\int _s^t\mathcal{A }(\tau )d\tau   ) ds \\&  + \int _0^1 \left ( ( \mathcal{D}   (s, x+\int _s^t\mathcal{A }(\tau )d\tau   ) - \mathcal{D}   (s, x  )\right  ) ds
  .
\end{aligned}
\end{equation}
 The last two lines are $C ^{k_0-1} $  in  $z\in \C ^m$, $f\in \Sigma _{8k_0-4k_1}$ and $\varrho \in \R ^4$ with values in  $  \Sigma _{  4k_1-4k_0}$  where they have
  norm  smaller than $ C
 ( |z|+ \norma{f }
_{\Sigma _{8k_0-4k_1 }} ) ^{3} $. This implies that  when we substitute \eqref{eq:pull13}
in the first  term in the rhs
of \eqref{eq:pull12},  the  last two lines of  \eqref{eq:pull13}   can be absorbed in $ {\textbf{R}} $.
When we substitute \eqref{eq:qlin1chi1}, the first term in the rhs of \eqref{eq:pull13}  is equal to what follows:

\begin{equation}  \label{eq:pull14}\begin{aligned} &
 -\im  \int _0^1   (z^{\mu}
\overline{z}^{\nu}
    B_{\mu   \nu
}(\varrho (f) ) ) \circ \phi ^t dt - \int _0^1  P_d(\omega _0)(\mathcal{L} f) \circ \phi ^t dt
  .
\end{aligned}
\end{equation}
By Lemma  \ref{lem:qlin1chi}  the second term in \eqref{eq:pull14}
is $C ^{k_0-1} $  in  $z\in \C ^m$, $f\in \Sigma _{8k_0-4k_1}$ and $\varrho \in \R ^4$ with values in  $  \Sigma _{ l}$ for any $l$ and
with norm  $ O
 ( |z|+ \norma{f }
_{\Sigma _{4k_0-4k_1 }} ) ^{M_0+2} $. The corresponding terms in \eqref{eq:pull12} can be then absorbed in $ {\textbf{R}} $.
The first term in \eqref{eq:pull14} can be written as
\begin{equation}  \label{eq:pull15}\begin{aligned} &
 -\im  B_{\mu   \nu
}(\varrho (f) )   \int _0^1   (z^{\mu}
\overline{z}^{\nu}
   ) \circ \phi ^t_0 dt
\end{aligned}
\end{equation}
plus an error term which can be then absorbed in $ {\textbf{R}} $ since it
 is $C ^{k_0-1} $  in  $z\in \C ^m$, $f\in \Sigma _{8k_0-4k_1}$ and $\varrho \in \R ^4$ with values in  $  \Sigma _{ 4k_1}$   and
with norm  smaller than $ C
 ( |z|+ \norma{f }
_{\Sigma _{8k_0-4k_1 }} ) ^{M_0+2} $.   By \eqref{eq:comp} we have that \eqref{eq:pull15} is
\begin{equation}  \label{eq:pull16}\begin{aligned} &
 -\im  B_{\mu   \nu
}(\varrho (f) )   z^{\mu}
\overline{z}^{\nu}
    -\im  B_{\mu   \nu
}(\varrho (f) ) F_{\mu   \nu
} (z,f, \varrho (f),b(\varrho (f)),B(\varrho (f)),
\end{aligned}
\end{equation}
with $F_{\mu   \nu
} (z,f, \varrho ,b ,B )\in \C$, $C^{\infty}$ in  $\zeta \in \C ^m$,  $f\in \Sigma _{ -4k_1}$, $\varrho \in \R ^4$
and in  $( b    , B   ) $, Furthermore $|F_{\mu   \nu
}|\le  C      |z | ^{M_0-1}
 ( |z|+ \norma{f }
_{\Sigma _{-4k_1}} ) .$ The contribution in \eqref{eq:pull12} is

\begin{equation}  \label{eq:pull16}\begin{aligned} &   ( z^{\mu}
\overline{z}^{\nu}
    +   F_{\mu   \nu
})\langle F (   z ',f',\varrho (f') )|   B_{\mu   \nu
}(\varrho (f) )    f  \rangle     .
\end{aligned}
\end{equation}
Then proceeding as  in Lemma \ref{lem:pull0}  we get a contribution like in the rhs of \eqref{eq:pull11}.

\qed

\begin{lemma}
  \label{lem:pulldbis} Let  $\widehat{\resto} _d=\widehat{\resto} _d(  z,f,\varrho )$ be      $C ^{k_0+2} $  in  $z\in \C ^m$, $f\in \Sigma _{ -4k_1}$ and $\varrho \in \R ^4$
with values in $\R$ for $d\ge 2$. Suppose that the following inequalities hold:

\begin{equation}\label{eq:Rhat1}\begin{aligned} & |
\widehat{\resto} _d (z ,f, \varrho )|\le C \| f \|
_{\Sigma _{-4k_1}}^d , \\&    | \widehat{\resto} _2 (z ,f, \varrho
)|  \le C (|z|+|\varrho |+ \| f \| _{ \Sigma _{-4k_1}}) \| f \| _{
\Sigma _{-4k_1 }}^2.
\end{aligned}\end{equation}
Then for any $d=2,...,5$ and for any pair $(k_2,k_3)$  there are $k_1(d)$ and $k_0(d)$ such that for $k_1\ge k_1(d)$, $k_0\ge k_0(d)$ and
$k_1\ge C k_0$ for some fixed large constant $C$, then the following occurs:
if
$d\ge 3$,
$\widehat{\resto} _d (   z ',f',\varrho (f') ) $ is of the form
$\resto ^{(M_0+1)} $;
if $d=2$ we have
\begin{equation}  \label{eq:pull21}\begin{aligned} &  \widehat{\resto} _2 (   z ',f',\varrho (f'))   =   \widehat{\resto} _2 (   z  ,f ,\varrho (f ) )
+\text{rhs like \eqref{eq:pull00}} .
\end{aligned}
\end{equation}
\end{lemma}
  \proof We will only sketch the case $d=2$, the others being similar.
   Schematically $\widehat{\resto} _2=G f^2$, with  $G(z  ,f ,\varrho   ) \in C ^{k_0 } $ with values in $B^2(\Sigma _{ -4k_1}, \C )$  and with  $\|   G(z  ,f ,\varrho   )  \|  _{B^2(\Sigma _{ -4k_1}, \C )} =O(|z|+|\varrho |+ \| f \| _{ \Sigma _{-4k_1}})$.
	 We have
   \begin{equation}   \label{eq:pull22}\begin{aligned} &   G (   z ',f',\varrho (f')  ) (f')^2=G (   z ',f',\varrho (f')  )  \mathcal{G}^2    +\\&   2  G (   z ',f',\varrho (f')  )   \mathcal{G}  e^{ \im  \sigma _3    {\gamma}   }\tau _{\mathcal{\mathbf{A}}} f     +G (   z ',f',\varrho (f')  )  \left ( e^{ \im  \sigma _3   {\gamma}   } \tau _{\mathcal{\mathbf{A}}}f\right  )^2   .
\end{aligned}
\end{equation}
The last term can be easily see to be of the  form $\widehat{\resto} ^{(M_0+1)}_2$. The
first two terms can be treated like in  Lemmas \ref{lem:pull0} and \ref{lem:pull1}.
For example, for a $\widetilde{\mathbf{R}}$ which can be absorbed in $ {\mathbf{R}}$, we have
\begin{equation*}
    G (   z ',f',\varrho (f')  )  \mathcal{G}^2  = G (   0,0,\varrho (f )  )    \mathbf{{G}}_1^2+ \widetilde{\mathbf{R}}.
\end{equation*}
This follows from the same argument used for \eqref{eq:pull02}
The first term in the rhs yields a term like the rhs of \eqref{eq:pull00} exactly
by the same argument used for \eqref{eq:pull06}, this time using the second inequality in \eqref{eq:Rhat1}

\qed

\begin{lemma}
  \label{lem:pullpol}
  Let
\begin{equation}  \label{eq:pullpol3}\begin{aligned} &    \psi  =   \sum _{|\mu +\nu |=M +1} d_{\mu\nu} (\varrho (f))  z^{\mu} \overline{z}^{\nu} + \sum _{|\mu +\nu |=M } z^{\mu}
\overline{z}^{\nu}
 \langle  \sigma _1\sigma _3  D_{\mu   \nu
}(\varrho (f))
  | f \rangle
\end{aligned}\nonumber
\end{equation}
  be another polynomial like $\chi$, in particular with $ d_{\mu   \nu } \in C^{k_0 }(\mathbb{R}  ^4 , \mathbb{C})$ and $ D_{\mu   \nu }
\in C^{k_0 }(\mathbb{R} ^4,   \Sigma _{4k_1} (\mathbb{R}^3,
\mathbb{C}^2))$ with $M_1\ge 2$.
Then,  for any pair $(k_2,k_3)$  there are $k _{10}  $ and $k _{00}$ such that for $k_1\ge k_{10}$, $k_0\ge k_{00}$ and
$k_1\ge C k_0$ for some fixed large constant $C$, then
  we have,  for an  $\textbf{{R}}  $ like in
  Lemma \ref{lem:pull0},

  \begin{equation}  \label{eq:pullpol1} \begin{aligned} & \psi \circ \phi  -  \psi =   {\textbf{R}} .
\end{aligned}
\end{equation}

  \end{lemma}

\proof
We have $   \psi  \circ
\phi  = \psi  +  \int _0^1  \{ \psi  , \chi \}  \circ \phi ^t dt .$
By elementary computation using \eqref{PoissBra2}  we get
\begin{equation}  \label{eq:pullpol4}\begin{aligned} &     \{ \psi  , \chi \} = \im
 \partial_{\overline{j}}\psi \partial_{j}\chi -\im
 \partial_{\overline{j}}\chi \partial_{j}\psi  +\im z^{\alpha +\mu }\overline{z}^{\beta  +\nu } \langle D_{\alpha   \beta
}  |\sigma _1\sigma _3B_{\mu   \nu
}  \rangle  + \\&    \partial _{\varrho _i}\psi      z^{ \mu }\overline{z}^{ \nu } \langle \sigma _1\sigma _3 B_{\mu   \nu
}  |\{ \varrho _i(f),f \} \rangle -  \partial _{\varrho _i}\chi      z^{ \alpha }\overline{z}^{ \beta  } \langle \sigma _1\sigma _3 D_{\alpha   \beta
}  |\{ \varrho _i(f),f \} \rangle
 \\& + (\partial _{\varrho _i}\psi \partial _{\varrho _j}\chi -
\partial _{\varrho _j}\psi \partial _{\varrho _i}\chi )   \{ \varrho _i(f)  , \varrho _j(f) \} .
\end{aligned}
\end{equation}
By the formulas \eqref{eq:L2}  all the terms   the last   line
are of the form $\widehat{\resto} _2$ like in Lemma \ref{lem:pulldbis}. Then by Lemma \ref{lem:pulldbis} we have
\begin{equation}  \label{eq:pullpol5}\begin{aligned} &  \widehat{\resto} _2 \circ \phi ^t   =   \widehat{\resto} _2
+\text{rhs like \eqref{eq:pull00}} _t,
\end{aligned}
\end{equation}
where the last term depends on $t$.  Integrating in $t$ we eliminate this dependence. Hence the terms from the last   line  of \eqref{eq:pullpol4} are absorbed
in $\mathbf{R} $. A similar conclusion holds for the terms from the second   line  of \eqref{eq:pullpol4}, this time using the last line of \eqref{PoissBra2}. Finally, the first line of \eqref{eq:pullpol4} is of the form

\begin{equation}  \label{eq:pullpol6}\begin{aligned} &   h=    \sum _{|\mu +\nu |=M +M_0 } h_{\mu\nu} (\varrho  )  z^{\mu} \overline{z}^{\nu} + \sum _{|\mu +\nu |=M +M_0 -1} z^{\mu}
\overline{z}^{\nu}
 \langle  \sigma _1\sigma _3  H_{\mu   \nu
}(\varrho  )
  | f \rangle
\end{aligned} \nonumber
\end{equation}
with coefficients  $ h_{\mu   \nu } \in C^{k_0 }(\mathbb{R}  ^4 , \mathbb{C})$ and $ H_{\mu   \nu }
\in C^{k_0 }(\mathbb{R} ^4,   \Sigma _{4k_1} (\mathbb{R}^3,
\mathbb{C}^2))$. Then $\int _0^1 h\circ \phi ^t dt$ can be absorbed in $\mathbf{R} $.
\qed

\begin{lemma}
  \label{lem:pullfin} Let $F=F(  z,f,\eta ,\varrho )$ be      $C ^{k_0} $  in  $z\in \C ^m$, $f\in \Sigma _{ -4k_1}$, $\eta \in \C^2$ and $\varrho \in \R ^4$
with values in $\Sigma _{ 4k_1   }(\R ^3, B^d(\C , \C))$. If $d=2$
 let $F(0,0, 0,0,0)=0$.
Set
\begin{equation*}
    R(z,f,\varrho )= \int
_{\mathbb{R}^3} F (x, z  ,f , f (x ),\varrho  )  f  ^d(x)   dx.
\end{equation*}
Then, for any $d=2,...,5$ and for any pair $(k_2,k_3)$  there are $k_1(d)$ and $k_0(d)$ such that for $k_1\ge k_1(d)$, $k_0\ge k_0(d)$,
$k_1\ge C k_0$ for some fixed large constant $C$  and for ${k}_{\mu   \nu}  $,   ${K}_{\mu   \nu}  $ and ${\textbf{R}} $ like in Lemma \ref{lem:pull0}, we have

\begin{equation}  \label{eq:fin0}\begin{aligned} &
R \circ \phi =  \sum _{|\mu + \nu |=M_0 +1 }  {{k}}_{\mu \nu}  ( \varrho (f), b( \varrho (f)), B( \varrho (f)) )z^{\mu} \overline{z}^{\nu} \\&  +\sum _{|\mu + \nu  |=M_0 } z^{\mu} \overline{z}^{\nu}
 \langle  \sigma _1\sigma _3    {{K}}  _{\mu   \nu
}( \varrho (f), b( \varrho (f)), B( \varrho (f)) )
  | f \rangle  +    {\textbf{R}}  .
\end{aligned}
\end{equation}
If $d>2$ formula \eqref{eq:fin0} holds with only ${\textbf{R}} $
in the rhs.
\end{lemma}
\proof
The rhs of \eqref{eq:fin0}  can
be written as a sum of terms of the form for $0\le i \le d$

\begin{equation} \label{eq:fin1} \begin{aligned} &  \int
_{\mathbb{R}^3} F (x  , z' ,f',e^{   \sigma _3     {\gamma}    } \tau _{ \mathbf{A}}f(x)
    +  \mathcal{G}(x ) ,\varrho ( f'))
   \mathcal{G} ^{d-i} (x)
    (e^{   \sigma _3     {\gamma}    } \tau_{\mathbf{A } } f )^i  dx  . \end{aligned}
\end{equation}
 The terms with $i\ge 2$ satisfy the statement if we change variable of integration. In particular for $i=2$
 we have $F=0$ at $(z,f)=(0,0)$, where we exploit $\mathcal{G}=0$. If $d\ge 3$, then the \eqref{eq:fin1} can be incorporated in $\mathbf{R}$. We then consider the case
 $d=2$ and $i=0,1$. For $i=0$,   \eqref{eq:fin1}  is a sum of the form

\begin{equation} \label{eq:fin2} \begin{aligned} &  \int
_{\mathbb{R}^3} F (x  , z' ,f',  \mathcal{G}(x ) ,\varrho ( f'))
   \mathcal{G} ^{2} (x)  dx\, +\\&  \int
_{\mathbb{R}^3} G(x +\mathbf{A} , z' ,f',e^{   \sigma _3     {\gamma}    }  f(x)
    , \tau _{ -\mathbf{A}} \mathcal{G}(x ) ,\varrho ( f'))
  \mathcal{G} ^{2} (x+\mathbf{A})
   e^{   \sigma _3     {\gamma}    } f(x) dx , \end{aligned}\nonumber
\end{equation}
where the second line can be incorporated in $\mathbf{R}$ and the first line
  is like Lemma \ref{lem:pull0}.
If in \eqref{eq:fin1}  we have $d=2$ and $i=1$ we have an expression of the form

\begin{equation} \label{eq:fin4} \begin{aligned} &  \int
_{\mathbb{R}^3} F (x  , z' ,f',  \mathcal{G}(x ) ,\varrho ( f'))
   \mathcal{G} (x)  e^{ \sigma _3 \gamma }\tau _{\mathbf{A}  }f(x)   dx\\&  +\int
_{\mathbb{R}^3} G(x +\mathbf{A} , z' ,f',e^{   \sigma _3     {\gamma}    }  f(x)
    ,  \mathcal{G}(x+\mathbf{A} ) ,\varrho ( f'))
 (e^{ \sigma _3 \gamma } f(x))^2  dx  . \end{aligned}\nonumber
\end{equation}
The second line is absorbed in $\mathbf{R}$.
The first line  is

\begin{equation} \label{eq:fin5} \begin{aligned} &  \int
_{\mathbb{R}^3} F (x +\mathbf{A}  , z' ,f',  \tau _{-\mathbf{A}  }\mathcal{G}(x ) ,\varrho ( f'))
   \mathcal{G} (x+\mathbf{A})  e^{ \sigma _3 \gamma } f(x)   dx \\& =    \int
_{\mathbb{R}^3} F (x  , z' ,f',  0 ,\varrho ( f'))
   \mathcal{G}(x)    f  (x) dx + \widetilde{\mathbf{R}}     \end{aligned}
\end{equation}
where $\widetilde{\mathbf{R}}  $ can be  absorbed in $\mathbf{R}$.
So we can apply Lemma \ref{lem:pull1}.

\qed

\subsection{Proof of Theorem \ref{th:main}: the step $\ell =2$}
\label{subsec:step1} At this stage  our goal is to obtain a
hamiltonian similar to $H$ but with  $\widetilde{\resto ^{(1)}} =
0$ in \eqref{eq:ExpH2resto}. In Lemma \ref{lem:ExpH} we can assume $M$
arbitrarily large. We
consider a polynomial  $\chi$, initially   unknown, like in \eqref{eq:chi} with $M_0=1$  and with $k_0$ and $k _1$ arbitrarily large with $1\ll k_0 \ll k_1 \ll M$.
We   choose $2N\ll k_2(2) \ll k_3(2) \ll  k_1-k_0$ with $k_2(2)$
as large as needed.
We write
\begin{equation} \label{eq:ExpH11} \begin{aligned} & H \circ
\phi=
 (\psi +H_2^{(1)} +\widetilde{\resto ^{(1)}}   +\widetilde{\resto ^{(2)}})
 \circ \phi ,
\end{aligned}
\end{equation}
for $\phi$ the Lie transform of $\chi$.   We have

\begin{equation*}
    H_2^{(1)}\circ
\phi  = H_2^{(1)}   +  \int _0^1 \{  H_2^{(1)}   , \chi \} \circ \phi ^t dt .
\end{equation*}
By the computations in Sect. \ref{sec:Normal form} we have schematically, for $\ell =1$,
\begin{equation}   \begin{aligned} &  \{  H_2^{(\ell)}   , \chi \} =   \im \sum _{|\mu +\nu |=\ell +1}\lambda ^{(\ell)}  ( \varrho )\cdot (\mu -
\nu) z^{\mu} \overline{z}^{\nu} b_{\mu \nu }   \\& + \im \sum _{|\mu +\nu |=\ell  }   z^{\mu} \overline{z}^{\nu} \langle f |\sigma _1\sigma _3( \lambda ^{(\ell)}  ( \varrho )\cdot (\mu -
\nu)- \mathcal{H} )B _{\mu \nu }
\rangle   \\& +\sum _{\substack{|\alpha  +\beta |=2\\ (\alpha , \beta )\neq (\delta _j , \delta _j )\, \forall \,j  }} a_{\alpha \beta}^{(\ell   )}( \varrho    )    \sum _{\mu \nu}\left ( b_{\mu   \nu  }       +
   \langle  \sigma _1\sigma _3  B_{\mu   \nu  }
  |   f   \rangle  \right ) \{z^\alpha
\overline{z}^\beta , z ^{\mu }
\overline{z}^{ \nu  }  \}  +L   ,
\end{aligned}\nonumber
\end{equation}
with $L$ like \eqref{eq:L1} with $\chi _1$ replaced by $\chi$.
Then, in the notation of  \eqref{eq:L0} and for  $\phi _0 ^t$   the flow of the     simplified system \eqref{eq:systchi3}, we have
\begin{equation}  \label{eq:ExpH13}\begin{aligned} & H_2^{(1)}\circ
\phi    - H_2^{(1)}   = \\& \sum _{|\mu +\nu |=2}\lambda ^{(1)}  ( \varrho (f))\cdot (\mu -
\nu)   b_{\mu \nu } (\rho (f))   \int _0^1   (z^{\mu}
\overline{z}^{\nu}) \circ \phi ^t_0 dt  \\& +  \sum _{|\mu +\nu |=1}
 \langle  \sigma _1\sigma _3B _{\mu \nu }(\varrho (f))
  ,      ( \lambda ( \varrho (f) )\cdot (\mu -
\nu)+ \mathcal{H} )\int _0^1(z^{\mu}
\overline{z}^{\nu}f )   \circ \phi ^t_0 dt\rangle  \\& +\sum _{\substack{|\alpha  +\beta |=2\\ (\alpha , \beta )\neq (\delta _j , \delta _j )\, \forall \,j  }} a_{\alpha \beta}^{(1   )}( \varrho  (f) ) \big [  \sum _{|\mu +\nu |=2}  b_{\mu   \nu  } ( \varrho  (f) )   \int _0^1\{z^\alpha
\overline{z}^\beta , z ^{\mu }
\overline{z}^{ \nu  } \}   \circ \phi ^t_0 dt \\& +
 \sum _{|\mu +\nu |=1}
  \langle  \sigma _1\sigma _3  B_{\mu   \nu  }( \varrho  (f) )
  |  \int _0^1 (f \{z^\alpha
\overline{z}^\beta , z ^{\mu }
\overline{z}^{ \nu  }  \} )\circ \phi ^t_0 dt \rangle \big ]
         \\&+    \widehat{{\textbf{R}}}_1   +\int _0^1  ( \langle
  V_j(\varrho (f)) f|
  f   \rangle  \widetilde{ \chi}  _j) \circ \phi ^t dt  ,
\end{aligned}
\end{equation}
 where  by Lemmas   \ref{lem:systchi2}  and  \ref{lem:systchi51} and by $\widehat{\textbf{R}}_1 =O(|z|+\| f\| _{\Sigma _{4k_0+4-4k_1 }})^{3}$ s.t.
 $\widehat{\textbf{R}}_1$ is $C^{ k_0-1 }$ in $z\in \C ^m$, $f\in \Sigma _{4k_0+4-4k_1}$ and $\varrho (f)$.
 In particular we have used estimates like
 \begin{equation}  \label{eq:ExpH133}\begin{aligned} & |\langle  \sigma _1\sigma _3  B_{\mu   \nu
} (\varrho (f))
|\mathcal{H}    f \circ \phi  ^t - \mathcal{H}    f \circ \phi _0^t   \rangle |
\le  \| B_{\mu   \nu
}(\varrho (f))\|  _{\Sigma _{ 4k_1}} \\& \times  \| \mathcal{H}    f \circ \phi  ^t - \mathcal{H}    f \circ \phi _0^t\|  _{\Sigma _{-4k_1}}\le C  \|      f \circ \phi  ^t - f \circ \phi _0^t\|  _{\Sigma _{2-4k_1 }}\\& \le C' ( |z|+ \norma{f }
_{\Sigma _{ 3-4k_1  }} ) ^{2} ,
 \end{aligned}
\end{equation}
with the latter a consequence of Lemma \ref{lem:systchi51}. In \eqref{eq:ExpH13} the last term
 is like in Lemma \ref{lem:NLhom2}.  It can be treated by
Lemma \ref{lem:pulldbis}.
By our choice of $k_2(2)$ and $k_3(2)$,  if we denote by
$\widetilde{\textbf{{R}}}_1$  the last line of  \eqref{eq:ExpH13},
we conclude that
$\widetilde{\textbf{{R}}}_1$
can be absorbed in $\resto _0^{(2)}+\resto _1^{(2)}
 +\widehat{\resto} _2^{(2)}$.
 We then
 obtain for $\ell =1$
\begin{equation}  \label{eq:ExpH14}\begin{aligned} & H_2^{(\ell )}\circ
\phi  = H_2^{(\ell )}   +\im \sum_{|\mu +\nu
|=\ell +1}b_{\mu \nu}^{(\ell )}( \varrho (f)) \lambda ( \varrho (f))\cdot (\mu - \nu) z^{\mu} \overline{z}^{\nu}
 \\& - \im \sum_{|\mu +\nu |=\ell } z^{\mu} \overline{z}^{\nu} \langle f
|\sigma _1\sigma _3 (\mathcal{H} -\lambda   \cdot (\mu - \nu)
)B_{\mu   \nu }^{(\ell )}( \varrho (f))  \rangle \\& +\sum _{|\mu + \nu |=\ell +1 }  {{k}}_{\mu \nu}^{(\ell )} ( \varrho (f), b( \varrho (f)), B( \varrho (f)) )z^{\mu} \overline{z}^{\nu} \\&  +\sum _{|\mu + \nu  |=\ell } z^{\mu} \overline{z}^{\nu}
 \langle  \sigma _1\sigma _3    {{K}} ^{(\ell )}_{\mu   \nu
}( \varrho (f), b( \varrho (f)), B( \varrho (f)) )
  | f \rangle  +    {\textbf{R}}_\ell\end{aligned}
\end{equation}
where: $\ell =1$;  $ {\textbf{R}}_1$ is like $ \widetilde{{\textbf{R}}}_1$;
${k}_{\mu   \nu} (0,b,B)=0$ and   ${K}_{\mu   \nu} (0,b,B)=0$ (follows by
\eqref{eq:ExpHcoeff1}); ${k}_{\mu   \nu} (\varrho ,b,B) \in \C  $ and   ${K}_{\mu   \nu} (\varrho ,b,B)\in \Sigma _{4 k_1} $ are $C^{ k_0 }$ in $\varrho$,
  in   $b _{\mu \nu}\in \C$
and   in $B _{\mu \nu}\in \Sigma _{4k_1}$. Notice that to get  $C^{ k_0 }$
regularity it is crucial the use   of $\phi _0^t$ and   its properties
stated  under \eqref{eq:comp}. The  $C^{ k_0 }$ regularity is key   for the homological equation.

 By Lemma \ref{lem:pullpol} we have
\begin{equation}  \label{eq:ExpH16}\begin{aligned} &   \widetilde{\resto ^{(1)}}  \circ
\phi  -   \widetilde{\resto ^{(1)}}  =  \sum _{|\mu + \nu |=2 } \widetilde{{k}}_{\mu   \nu}  ( \varrho (f), b( \varrho (f)), B( \varrho (f)) )z^{\mu} \overline{z}^{\nu}  \\& +\sum _{|\mu + \nu |=1} z^{\mu} \overline{z}^{\nu}
 \langle  \sigma _1\sigma _3   \widetilde{{K}} _{\mu   \nu
}( \varrho (f), b( \varrho (f)), B( \varrho (f)) )
  | f \rangle   +\widetilde {\mathbf{S}}_1 ,
\end{aligned}
\end{equation}
with:     $ \widetilde{\textbf{S}}_1 $ like
 $ \resto  ^{(2)}$;    $\widetilde{{k}}_{\mu   \nu}   $ (resp. $\widetilde{K}_{\mu   \nu}   $) is  like  $ {{k}}_{\mu   \nu}  $ (resp. $ {K}_{\mu   \nu}  $).

 By  Lemma  \ref{lem:systchi52} we have
  \begin{equation*}
  \psi (\varrho (f))  \circ
\phi  =  \psi (\varrho (f) ) \circ
\phi _0  +\mathbf{T}  _\ell
  \end{equation*}
 for $\ell =1$ where   $\mathbf{T}  _1$ is like  $ \resto  ^{(2)}$.
 By Lemma  \ref{lem:quasilin7} applied to $\phi _0^t$, exploiting the fact that    $\psi (0) =0$, we have that
 \begin{equation}  \label{eq:pullpsi}\begin{aligned} &   \psi (\varrho (f))  \circ
\phi _0 -\psi (\varrho (f))  = \sum _{|\mu + \nu |=\ell +1 }  {{k}}_{\mu \nu}^{(\ell ) \prime } ( \varrho (f), b( \varrho (f)), B( \varrho (f)) )z^{\mu} \overline{z}^{\nu} +\\&    \sum _{a=1}^3 \sum _{i=0}^1\sum _{|\mu + \nu  |=\ell } z^{\mu} \overline{z}^{\nu}
 \langle  \sigma _1\sigma _3   \partial _a ^i{{K}} ^{(\ell )ia}_{\mu   \nu
}( \varrho (f), b( \varrho (f)), B( \varrho (f)) )
  | f \rangle  +    {\textbf{R}}_\ell '
\end{aligned}
\end{equation}
  with $  {\textbf{R}}_\ell ' $ like $  {\textbf{R}}_\ell .$

By Lemma \ref{lem:pullfin}
\begin{equation}  \label{eq:ExpH18}\begin{aligned} &  \langle B_{2 } (   z ',\varrho (f') )|   ( f ')^2   \rangle     =   \widehat{\mathbf{S}}_1\\& +    \sum _{|\mu + \nu |=2 }  \upsilon_{\mu   \nu} ( \varrho (f), b( \varrho (f)), B( \varrho (f)))z^{\mu} \overline{z}^{\nu} \\& +\sum _{|\mu + \nu |=1} z^{\mu} \overline{z}^{\nu}
 \langle  \sigma _1\sigma _3   \Upsilon _{\mu   \nu
}( \varrho (f), b( \varrho (f)), B( \varrho (f)))
  | f \rangle ,
\end{aligned}
\end{equation}
 where  $\upsilon _{\mu   \nu}$ and  $\Upsilon_{\mu   \nu}$
have the same properties of $ {{k}}_{\mu   \nu}$,  $ {{K}}_{\mu   \nu}$, $\widehat{\mathbf{S}}_1 $ is
like $ \resto  ^{(2)}$.
Notice that the fact that $\upsilon _{\mu   \nu}$ and  $\Upsilon_{\mu   \nu}$
are $C^{k_0}$ is key here for the homological equation.

By Lemma \ref{lem:pulldbis}, where we are using  $k_0\ll k_1 \ll M$, which is much more than needed,
we have for $\ell =1$
\begin{equation}  \label{eq:h192}  \begin{aligned} &   \widehat{\resto}^{(\ell )}_2 \circ \phi =   \widehat{\resto}^{(\ell )}_2   +\sum _{|\mu + \nu |=\ell +1 }  {{\kappa}}_{\mu   \nu} ^{(\ell  )} ( \varrho (f), b( \varrho (f)), B( \varrho (f)) )z^{\mu} \overline{z}^{\nu}  \\& +\sum _{|\mu + \nu |=\ell } z^{\mu} \overline{z}^{\nu}
 \langle  \sigma _1\sigma _3    {{\mathcal{K}}} ^{(\ell  )} _{\mu   \nu
}( \varrho (f), b( \varrho (f)), B( \varrho (f)) )
  | f \rangle   +\widehat {\mathbf{S}}_\ell ,
\end{aligned}
\end{equation}
where $ {{\kappa}}_{\mu   \nu}^{(1  )} $, ${{\mathcal{K}}}  _{\mu   \nu
}^{(1  )}$ and $\widehat {\mathbf{S}}_1  ^{(1 )} $ are like   $ {{k}}_{\mu   \nu}  $, $ {K}_{\mu   \nu}   $     and $ \resto  ^{(2)}$.
Consider now $\mathbf{K}:=\widetilde{\resto ^{(1)}}$   and  $\widetilde{\mathbf{K}}  $  the polynomial of the form \eqref{eq:tildeKrho} with
\begin{equation} \begin{aligned} & \widetilde{\mathbf{k}}_{\mu \nu}  (\varrho ,b,B):= (k _{\mu   \nu}  +k _{\mu   \nu}^{(1)\prime } +\widetilde{k} _{\mu   \nu} +\upsilon _{\mu   \nu}  +\kappa ^{(1  )} _{\mu   \nu} ) (\varrho ,b,B) \, ,\\& \widetilde{\mathbf{K}}_{\mu \nu}  (\varrho ,b,B):=( K _{\mu   \nu} +\widetilde{K} _{\mu   \nu} +\sum _{i=0}^1
\sum _{a=1}^3\partial _a ^i  {K} _{\mu   \nu}^{(1  )ia} +\Upsilon _{\mu   \nu} +\mathcal{K} _{\mu   \nu}^{(1  )}) (\varrho ,b,B).
\end{aligned}\nonumber \end{equation}
  Then $\mathbf{K} (\varrho )$ and  $\widetilde{\mathbf{K}}  (\varrho ,b,B)$
  are like in Lemmas \ref{lem:NLhom1}-\ref{lem:NLhom2}. This means that we can choose $\chi$ in \eqref{eq:chi1} so that \begin{equation} \begin{aligned} & \im \lambda (\varrho )\cdot (\mu -
\nu) z^{\mu} \overline{z}^{\nu} b_{\mu \nu }    +z^{\mu} \overline{z}^{\nu}  \im \langle f |\sigma _1\sigma _3( \lambda  (\varrho ) ,\cdot (\mu -
\nu)- \mathcal{H} )B _{\mu \nu }
\rangle\\&  + \widetilde{\resto ^{(1)}}+\widetilde{\mathbf{K}}  (\varrho (f), b( \varrho (f)), B( \varrho (f)))= Z_1 (\varrho (f))  ,
\end{aligned}\nonumber \end{equation}
where  $Z_1(\varrho)$ is in normal form and  homogeneous of degree $2$ in $(z,\overline{z},f)$ and with $Z_1(0)=0$.

 For $d>2$, the terms $  \langle B_{d } (   z ',\varrho (f') )|   ( f ')^d   \rangle
$ can
can be incorporated  in  ${\resto} ^{(2)}$ by Lemma \ref{lem:pullfin}.
This is true also for the $d=5$ term and for $E_P(f')$.
We set
\begin{equation} \begin{aligned}    \widetilde{\resto} ^{(2)}:&=\widehat{\resto} ^{(1)}_2 +  \mathbf{R}_1+   \widehat{\mathbf{S}}_1+\widetilde{\mathbf{S}}_1+\textbf{T}_1 +  \mathbf{R}_1'\\& + \widetilde{\resto ^{(2)}}
 \circ \phi -  \langle B_{2 } (   z' ,\varrho (f') )|    (f ') ^{   2} \rangle -\widehat{\resto} ^{(1)}_2\circ \phi .
\end{aligned}\nonumber \end{equation}
All the terms in the rhs  have the properties required
to   ${\resto} ^{(2)}$. Hence, if we also set ${\resto} ^{(2)}:=\widetilde{\resto} ^{(2)}$, we conclude  the proof of case $\ell =2$ in Theorem
\ref{th:main}.

\subsection{Proof of Theorem \ref{th:main}:  the step $\ell >2$}
\label{subsec:step2}

  Case $\ell=2$  has been  treated in Subsection \ref{subsec:step1}.
 We proceed by
induction to complete the proof of Theorem \ref{th:main}. From the
argument below one can see that $H_2^{(\ell )} =H_2^{(2)} $ for all
$\ell\ge 2$.  Suppose that the statement of Theorem \ref{th:main} holds
for an $\ell\ge 2$.

Set $k_0=k_2(\ell )-\underline{k}$ and $k_1=k_3(\ell )-\underline{k}$
for a fixed and appropriately large $\underline{k}$.
We will choose
$ 2N \ll k_2(\ell +1 )\ll k_3(\ell +1 ) \ll k_1-k_0   $ with $k_2(\ell +1 )$ as large as needed.

  Since
$H^{(\ell )}=H\circ \Tr _\ell $ is real valued (because $H$ is real valued),
    $ {a}_{\mu \nu }^{(\ell )}$ and $ {G}_{\mu \nu
}^{(\ell )}$ satisfy \eqref{eq:ExpKcoeff2}. We seek an appropriate  polynomial $\chi$ as in
\eqref{eq:chi1} with $M_0=\ell $. For any such polynomial, we consider its
Lie transform $\phi =\phi ^1$.   Proceeding like in the previous step of the proof,
  we   obtain formula \eqref{eq:ExpH14}
with:    $ {\textbf{R}}_\ell =O(|z|+\| f\| _{\Sigma _{-4k_1+2}})^{\ell  +2}$, with $ {\textbf{{R}}}_\ell \in C ^{ k_0-1}$
 in $z\in \C^m$, $f\in \Sigma _{4k_0-4k_1 }$ and $\rho (f)$;  $ {k}^{(\ell )}_{\mu   \nu} ( 0,b,B  )=0$ resp. $ {K}^{(\ell )}_{\mu   \nu} ( 0,b,B )=0$;  $\widehat{k}_{\mu   \nu} ^{(\ell )}( \varrho ,b,B  ) \in \C $ and $\widehat{K}^{(\ell )}_{\mu   \nu} ( \varrho,b,B ) \in \Sigma _{4k_1} $ are $C^{ k_0 }$ in $\varrho \in \R^4$, $b_{\mu   \nu}\in \C$ and $B_{\mu   \nu} \in \Sigma _{4k_1} $. By our choice of $k_2(\ell +1)$ and $k_3(\ell +1)$,  $ {\textbf{R}}_\ell $ can be absorbed in $\resto ^{(\ell +1 )}_0+\resto ^{(\ell +1 )}_1+\widehat{\resto} ^{(\ell +1 )}_2$.

By Lemma \ref{lem:pullfin} we have

\begin{equation}  \label{eq:main15}\begin{aligned} &  \langle F_{2 }^{(\ell )} (   z ',f',  f'(\cdot ), \varrho (f') )|   (f')^2  \rangle  = \langle F_{2 }^{(\ell )} (   z  ,f ,  f  (\cdot ), \varrho (f ) )|   f^2     \rangle \\&   + \widehat{\mathbf{S}}_\ell+
     \sum _{|\mu + \nu |=\ell +1 }  \upsilon_{\mu   \nu}^{(\ell )} ( \varrho (f), b( \varrho (f)), B( \varrho (f)) )z^{\mu} \overline{z}^{\nu}   \\& +\sum _{|\mu + \nu |=\ell} z^{\mu} \overline{z}^{\nu}
 \langle  \sigma _1\sigma _3   \Upsilon ^{(\ell )}_{\mu   \nu
}( \varrho (f) , b( \varrho (f)), B( \varrho (f)) )
  | f \rangle  ,
\end{aligned}
\end{equation}
where:    $\upsilon ^{(\ell )}_{\mu   \nu} $ (resp.  $\Upsilon ^{(\ell )}_{\mu   \nu}  )$  is like  $k_{\mu   \nu}^{(\ell )} $ (resp.  $K_{\mu   \nu} ^{(\ell )} )$;
  $ \widehat{\mathbf{S}}_\ell $ is like  $ \resto ^{(\ell +1)}$.
    Proceeding like in \eqref{eq:main15},  by Lemma \ref{lem:pulldbis} we get \begin{equation}  \label{eq:h193}  \begin{aligned} &   \widehat{\resto}^{(\ell )}_2 \circ \phi -   \widehat{\resto}^{(\ell )}_2   =\sum _{|\mu + \nu |=\ell +1 }  {{\kappa}}_{\mu   \nu} ^{(\ell  )} ( \varrho (f), b( \varrho (f)), B( \varrho (f)) )z^{\mu} \overline{z}^{\nu}  \\& +\sum _{|\mu + \nu |=\ell } z^{\mu} \overline{z}^{\nu}
 \langle  \sigma _1\sigma _3    {{\mathcal{K}}} ^{(\ell  )} _{\mu   \nu
}( \varrho (f), b( \varrho (f)), B( \varrho (f)) )
  | f \rangle   +\widehat {\mathbf{S}}_\ell ,
\end{aligned}
\end{equation}
with  $ \widetilde{\textbf{S}}_\ell $   like $ \resto ^{(\ell +1)}$,  $\widetilde{{k}}_{\mu   \nu}^{(\ell )}  $ resp. $\widetilde{K}_{\mu   \nu}^{(\ell )}  $ with the properties of  $ {{k}}_{\mu   \nu}^{(\ell )}  $ resp. $ {K}_{\mu   \nu}^{(\ell )}  $. Proceeding as for the $\ell =2$ case,
we have that $\psi (\varrho (f))\circ \phi -\psi (\varrho (f))$ is
like the right hand side of \eqref{eq:pullpsi}.

\noindent Set $\mathbf{K}^{(\ell )}(\varrho (f)):=\resto^{(\ell )}_{0 }+ \resto^{(\ell )}_{1} .$
Consider the polynomial $\widetilde{\mathbf{K}} $   of the form \eqref{eq:tildeKrho} with coefficients
\begin{equation} \begin{aligned} & \widetilde{\mathbf{k}} ^{(\ell )}_{\mu \nu} (\varrho ,b,B):= (k _{\mu   \nu}^{(\ell )} +k _{\mu   \nu}^{(\ell)\prime } +\widetilde{k} _{\mu   \nu}^{(\ell )}  + \upsilon ^{(\ell )}_{\mu   \nu}) (\varrho ,b,B)\, , \\& \widetilde{\mathbf{K}}^{(\ell )} _{\mu \nu} (\varrho ,b,B):= (K _{\mu   \nu} ^{(\ell )}  +\sum _{i=0}^1
\sum _{a=1}^3\partial _a ^i  {K} _{\mu   \nu}^{(\ell  )ia} +{\mathcal K} _{\mu   \nu} ^{(\ell )} +\Upsilon _{\mu   \nu} ^{(\ell )}) (\varrho ,b,B).
\end{aligned}\nonumber \end{equation}

 Then $\mathbf{K}  ^{(\ell )}(\varrho )$ and  $\widetilde{\mathbf{K}} ^{(\ell )}  (\varrho ,b,B)$
  are like in Lemma  \ref{lem:NLhom1}.
  This means that we can choose $\chi$ in \eqref{eq:chi1} so that \begin{equation} \begin{aligned} & \im \lambda (\varrho )\cdot (\mu -
\nu) z^{\mu} \overline{z}^{\nu} b_{\mu \nu }    +z^{\mu} \overline{z}^{\nu}  \im \langle f |\sigma _1\sigma _3( \lambda  (\varrho ) ,\cdot (\mu -
\nu)- \mathcal{H} )B _{\mu \nu }
\rangle\\&  + \mathbf{K}  ^{(\ell )}(\varrho )+\widetilde{\mathbf{K}} ^{(\ell )}  (\varrho ,b,B) = Z_\ell  (\varrho  )  ,
\end{aligned}\nonumber \end{equation}
where  $Z_{ \ell +1 }(\varrho)$ is in normal form and  homogeneous of degree $\ell +1$ in $(z,\overline{z},f)$.  Set $Z ^{(\ell +1)}:=Z ^{(\ell )}+Z _{ \ell +1 }$ and

 \begin{equation}
\label{Hamr1}  \begin{aligned} & H^{(\ell +1)}:=H^{(\ell )}\circ\phi  =H_2^{(\ell )}+Z^{(\ell +1)} +\widetilde{\resto }^{(\ell +1)},
  \end{aligned}
\end{equation}
where
\begin{equation}
\label{Hamr2}  \begin{aligned}   \widetilde{\resto }^{(\ell +1)} :& =Z^{(\ell +1)} \circ\phi -Z^{(\ell +1)}
\\& +  \sum _{d=3}^{6}  {\resto}^{(\ell )}_d\circ\phi + \widetilde{\resto} _{2  }^{(\ell)} + \mathbf{R}_\ell+\widetilde{\mathbf{S}}_\ell+ \widehat{\mathbf{S}}_\ell  + \mathbf{R}_\ell ' .
\end{aligned}\nonumber \end{equation}
  The terms in the last   line  are of the form requested for terms
  of $ {\resto }^{(\ell +1)}$.  By Lemma \ref{lem:pullpol} also $Z^{(\ell +1)} \circ\phi -Z^{(\ell +1)} $ is of the same type. Then we set $ {\resto }^{(\ell +1)}:=\widetilde{{\resto }}^{(\ell +1)}$ and the proof is finished.

\qed

 \section{Dispersion}
 \label{sec:dispersion}
We apply Theorem \ref{th:main} for $\ell =2N+1 $  (recall $N=N_1$ where
$N_j\lambda _j<\omega _0 <(N_j+1)\lambda _j).$ In the rest of the
paper we work with the hamiltonian $H^{(2N+1 )}$. We will drop the upper
index. So we will set $H=H^{(2N+1 )}$, $H_2=H_2^{(2N+1 )}$, $\lambda
_j=\lambda _j^{(2N+1)}$,    $\lambda  =\lambda  ^{(2N+1 )}$,
$Z_i=Z_i^{(2N+1 )}$ for $i=0,1$ and $\resto =\resto ^{(2N+1 )}$. In
particular we will denote by $G_{\mu \nu}$ the coefficients $G_{\mu
\nu}^{(2N+1)}$ of $Z_1^{(2N+1)}$. We will show:
\begin{theorem}\label{proposition:mainbounds} Consider the constant   $0<\epsilon <\varepsilon  $  of Theorem \ref{theorem-1.2}.
There is a fixed
$C >0$ such that for    $\varepsilon   $  sufficiently small and for any
$\epsilon \in (0, \varepsilon  )$ we have
\begin{eqnarray}
&   \|  f \| _{L^p_t( [0,\infty ),W^{ 1 ,q}_x)}\le
  C \epsilon \text{ for all admissible pairs $(p,q)$,}
  \label{Strichartzradiation}
\\& \| z ^\mu \| _{L^2_t([0,\infty ))}\le
  C \epsilon \text{ for all multi indexes $\mu$
  with  $\lambda\cdot \mu >\omega _0 $,} \label{L^2discrete}\\& \| z _j  \|
  _{W ^{1,\infty} _t  ([0,\infty )  )}\le
  C \epsilon \text{ for all   $j\in \{ 1, \dots , m\}$ }
  \label{L^inftydiscrete} .
\end{eqnarray}
\end{theorem}
Notice that by the time reversibility of the NLS, the above estimates
imply the ones with $\R$ replacing $ [0,\infty )$,
doubling constants in \eqref{Strichartzradiation}--\eqref{L^2discrete}.

\eqref{L^inftydiscrete} is a consequence of the already known orbital stability, so we do not need to prove it.
To obtain
 Theorem
\ref{proposition:mainbounds}  it is enough to show that   there are fixed constants $C_1$, $C_2$ (large) and $\varepsilon$ (small) such that if for $\epsilon \in (0, \varepsilon  )$ (where  $  \epsilon$ and $\varepsilon  $  are those of Theorem \ref{theorem-1.2})

\begin{eqnarray}
&   \|  f \| _{L^p_t([0,T],W^{ 1 ,q}_x)}\le
  C _1\epsilon \text{ for all admissible pairs $(p,q)$,} \label{4.4a}
\\& \| z ^\mu \| _{L^2_t([0,T])}\le
 C_2 \epsilon \text{ for all multi indexes $\mu$
  with  $\omega \cdot \mu >\omega _0 $} ,\label{4.4}
\end{eqnarray}
 then in fact   \eqref{4.4a} and
\eqref{4.4}  hold  but with   $C_1$, $C_2$ replaced
by $C_1/2$, $C_2/2$. In fact we conclude that these estimates hold
for all $T$ and so \eqref{Strichartzradiation}--\eqref{L^2discrete} hold.
The proof consists in three main steps.
\begin{itemize}
\item[(i)] Estimate $f$ in terms of $z$.
\item[(ii)] Substitute the variable $f$  with a
new "smaller" variable $g$ and find smoothing estimates for $g$.
\item[(iii)] Reduce the system for $z$ to a closed system involving
only the $z$ variables, by insulating the  part of $f$  which
interacts with $z$, and by decoupling the rest (this reminder is
$g$). Then clarify the nonlinear Fermi golden rule.
\end{itemize}

\subsection{Proof of Theorem \ref{proposition:mainbounds}:   step (i)}
\label{subsec:stepi}

Step (i) is encapsulated by the following proposition:

\begin{proposition}\label{prop:conditional4.2} Assume \eqref{4.4a}--\eqref{4.4}. Then there exist constants $C=C(C_1,C_2)$ and $ K_1=K_1(C_1)$,  such that, if
  $C(C_1,C_2) \epsilon  $  is sufficiently small,   then we have
\begin{eqnarray}
&   \|  f \| _{L^p_t([0,T],W^{ 1 ,q}_x)}\le
  K_1  \epsilon \text{ for all admissible pairs $(p,q)$}\ .
  \label{4.5}
\end{eqnarray}
\end{proposition}
\proof
Consider $Z_1$ of the form \eqref{e.12a}. Set:
 \begin{equation}\label{eq:G^0}      G_{\mu
\nu}^0=G_{\mu \nu}( \varrho (0) ) ;  \quad    \lambda
^0_j=\lambda _j(\omega _0).\end{equation}
Then we have (with
finite sums)
\begin{equation}\label{eq:f variable} \begin{aligned}  &\im \dot f -
\mathcal{H}f -      (\partial _{ Q(f)} H) P_c(\omega _0)
\sigma _3  f -   \im    (\partial _{ \Pi _a(f)} H) P_c(\omega _0)
\partial   _{x_a}  f\\& = \sum _{\substack{|\lambda  ^0 \cdot(\nu-\mu)|> \omega _0,\\
|\mu+\nu|\leq 2N_1+1}}
z^\mu \overline{z}^\nu
  G_{\mu \nu}^0  \\&  +
  \sum _{\substack{|\lambda ^0 \cdot(\nu-\mu)|>m-\omega _0,\\
|\mu+\nu|\leq 2N_1+1}} z^\mu
\overline{z}^\nu  (G_{\mu \nu} -  G_{\mu \nu}^0)    +  \sigma
_3 \sigma _1\widehat{\nabla} _f \resto  ,
\end{aligned}\end{equation}
with $\widehat{\nabla} _f \resto (z,f,\rho )$ the gradient in $f$,
with no differentiation in $\varrho (f)$.
In order to obtain bounds on $f$, we need bounds on the right hand term of the
equation especially the last two terms. They are provided by the following
lemma.

\begin{lemma}\label{lem:bound remainder}
Assume \eqref{4.4a}--\eqref{4.4}.
Then
there is a   constant $ C(C_1,C_2
)$ independent of $\epsilon$
such that the following is true: we have $\sigma
_3 \sigma _1\widehat{\nabla} _f \resto =R_1+ R_2$ with
\begin{equation}
\label{bound1:z1} \| R_1 \| _{L^1_t([0,T],H^{ 1 }_x)}+\|
   R_2 \| _{L^{2
 }_t([0,T],W^{ 1
,
 \frac{6}{5}}_x)}\le C( C_1,C_2
 ) \epsilon^2.
\end{equation}

\end{lemma}
\proof  The proof is standard,  a combination of \cite{bambusicuccagna} and \cite{cuccagnamizumachi}.\qed

 \begin{lemma}\label{lem:conditional4.21} Consider $\im \dot \psi -
\mathcal{H}\psi - \varphi  (t)
 \sigma_3P_c \psi -\im  A _a (t)
 P_c \partial _{x_a}\psi=F$  where: $P_c=P_c(\omega _0)$,  $\psi= P_c\psi$, $\varphi $ and each $A_a$ are real valued.   Then there exist
 $c _0>0$ and
   $C>0$ such that if  $\| (\varphi , A )  \| _{L^\infty _t [0,T]}<c _0$   then
   for  $(p,q) $ as in Theorem \ref{proposition:mainbounds}   we have
\begin{equation} \label{eq:421}  \begin{aligned}
 &     \|  \psi \| _{L^p_t( [0,T ],W^{  1,q } )  }\le C\| \psi(0) \| _{H^{  1 }}+ C \|  F \|
_{L^1_t([0,T],H^{ 1 }_x) + L^2_t([0,T],W^{ 1  ,  \frac{6}{5}}_x)}  \end{aligned}
\end{equation}
\end{lemma}
  \proof This result is due to Beceanu, see for example Theorem 3.8 \cite{beceanu}.
\qed

{\it Continuation of the proof of Proposition \ref{prop:conditional4.2}.}
By \eqref{eq:f variable}
we can apply to $f$ Lemma  \ref{lem:conditional4.21}   by taking $\varphi (t)
= \partial _{Q(f)} H $,   $ {A}_a (t)
= \partial _{\Pi _a(f)} H $
and  $F=\text{rhs\eqref{eq:f
variable}}-\varphi   [\sigma _3,P_d]f $. Then, for fixed constants
\begin{equation} \label{eq:4212}  \begin{aligned}
 &     \|  f \| _{L^p_t( [0,T ],W^{1,q   }  ) }\le C_1\|
f(0)
\| _{H^{  1}}+ C_1 \|  F \| _{L^1_t([0,T],H^{ 1}_x) +  L^2_t([0,T],W^{ 1
,
 \frac{6}{5}}_x)}  \\& \le C_1\|
f(0)
\| _{H^{  1}}+
 C  \sum _{\lambda\cdot \mu >m-\omega _0 }\| z ^\mu \| _{L^2_t( 0,T ) }^2\\&  +C\|
R_1
\| _{L^1_t([0,T], H^{ 1}_x )  }+\|
   R_2 \| _{L^{2
 }_t([0,T],W^{ 1
,
 \frac{6}{5}}_x)} +C
\epsilon
  \|  f \| _{  L^2_t ([0,T], L^{ 6  }_x) }. \end{aligned}
\end{equation}
 For $\epsilon $ small this yields  Proposition \ref{prop:conditional4.2}
 by Lemma   \ref{lem:conditional4.21} and by \eqref{4.4}.\qed

\begin{lemma}
\label{lem:asymptotically f} Assume the conclusions of Theorem
\ref{proposition:mainbounds}. Then there exists a fixed $C>0$ and $f_+ \in
H^{1 } $ with $\| f_+ \| _{H^{1}}<C \epsilon $ such that   we have
\begin{equation}\label{scattering111}  \lim_{t\to +\infty}
\left \|  \tau _{X(t)}e^{\im \chi (t) \sigma _3}f (t) -
 e^{  \im t \Delta  \sigma _3  }  {f}_+   \right
\|_{H^{1}}=0
 \end{equation} for  $ \chi(t) :=    t\omega _0+\int _0^t \partial _{Q(f)} H (t')dt'  $ and $X(t):= \int _0^t  \partial _{\Pi  (f)} H ((t')dt'.$
\end{lemma} \proof For $\psi (t)=f(t)$,   $F=\text{rhs\eqref{eq:f
variable}}-\varphi (t) [\sigma _3,P_d]f$, $\varphi (t)
= \partial _{Q(f)} H $,   $\mathcal{A}_a (t)
= \partial _{\Pi _a(f)} H $, $\mathcal{U} (t ) = e^{  \int _0^t \left (\im \sigma _3\varphi (\tau )+ \mathcal{A}(\tau )\cdot \nabla   \right )d\tau}$
and for $t_1<t_2$, we have
 \begin{equation} \label{eq:423scat} \begin{aligned} & \| e^{ \im \mathcal{H}_0 t_2}
\mathcal{U} (t_2) f (t _2)-e^{ \im \mathcal{H}_0 t_1}
\mathcal{U} (t_1) f (t _1) \| _{H^{1}} \le\\&
\| \int   _{t_1}^{t_2}
 e^{ \im \mathcal{H}_0  t' } \mathcal{U} (t')\left [ F(t')+ Vf (t ')  -
\varphi  (t')
  \sigma_3P_d f (t ') -\im \mathcal{A}_aP_d\partial _af(t')  \right ] dt' \| _{H^{1}}   \\& \le
 C  (\sum _{|\lambda  ^0 \cdot \mu |>m-\omega _0}
\| z^\mu \| _{L^2(t_1,t_2)} +\| R_1 \| _{L^1_t([t_1,t_2],H^{ 1 }_x)}\\& +\|
   R_2 \| _{L^{2
 }_t([t_1,t_2],W^{ 1, \frac{6}{5}}_x)}+ \|
   f \| _{L^{2
 }_t([t_1,t_2],W^{ 1, 6}_x)}) .\end{aligned}\nonumber
 \end{equation}
Since the rhs has limit 0 as $t_1\to +\infty$,
there exists $ {f}_+ \in H^{ 1 } $ such that
 \begin{equation}   \lim_{t\to +\infty}
\left \| \mathcal{U} (t )f (t) -
 e^{ -\im \mathcal{H}_0 t }{f}_+   \right \|_{H^{1}}=0 . \nonumber
 \end{equation}
   This yields Lemma \ref{lem:asymptotically f}.\qed

\begin{lemma}
\label{lem:part-1.2}  Assume the conclusions of Theorem
\ref{proposition:mainbounds} and the notation of   Theorem \ref{theorem-1.2}.
 Then  the conclusions of  Theorem \ref{theorem-1.2} hold
 with
 the $f_+$ of
     \eqref{scattering111}
     and with
   \begin{equation}
\label{eq:part-1.2}\begin{aligned} &  \widehat{\vartheta} =\chi +\im
\sum _{\ell =1}^{2N+1} \gamma _\ell \, , \quad  \widehat{D} =X -
\sum _{\ell =1}^{2N+1}  \mathbf{A} _ \ell ,
      \end{aligned}
\end{equation}
 with  $ \mathbf{A} _\ell$ and $ \gamma _\ell$
the terms in \eqref{eq:systchi11} corresponding to the lie transforms $\phi _\ell$ of Theorem \ref{th:main}.
\end{lemma}
\proof  This follows immediately from    \eqref{eq:systchi11}. Indeed, schematically \begin{equation*}
   e^{\im \chi-\sum _{\ell =1}^{2N+1} \gamma _\ell}\,  \tau _{X -
\sum _{\ell =1}^{2N+1}  \mathbf{A} _ \ell}f_{\text{\ref{theorem-1.2}}} \approx  e^{\im \chi }\,  \tau _{ X  }f_{\text{\ref{proposition:mainbounds}}} \approx e^{  \im t \Delta \sigma _3  }  {f}_+,
 \end{equation*}
 where $f_{\text{\ref{theorem-1.2}}}$ (resp. $f_{\text{\ref{proposition:mainbounds}}}$)   is the coordinate in Theorem \ref{theorem-1.2} (resp. Theorem\ref{proposition:mainbounds}).
 \qed

\begin{lemma}
\label{lem:vomega as} Denote by $(\omega ,v,z',f')$ the coordinates
\eqref{eq:decomp2} of the solution
in the initial system of coordinates (we omit $\vartheta ,D$).
Then $\lim _{t\to   \infty } z'(t)=0$,
there are      functions   $\theta
\in C^1(\R , \R )$ and
  $y\in C^1(\R , \R ^3 )$ s.t. \begin{equation}\label{scattering12}\lim _{t\to + \infty }
\left \|  \tau _{y  (t)}e^{\im \sigma _3\left (v (t)\cdot \frac{ x }{2} +  \theta   (t)\right )  }f '(t) -
 e^{  \im t \Delta \sigma _3  }  {f}_+   \right
\|_{H^{1}}=0
,
\end{equation}
with the $f_+$  of Lemma \ref{lem:asymptotically f},  and
$\lim _{t\nearrow \infty}\omega (t)=\omega _+$ and $\lim _{t\nearrow \infty}v (t)=v _+$.
\end{lemma}
\proof If we denote by $(  z', f')$ the initial coordinates  and by $(  z _{\text{\ref{theorem-1.2}}}, f_{\text{\ref{theorem-1.2}}} )$ the coordinates in \eqref{eq:SystK2} considered in Theorem \ref{theorem-1.2},
we have $z'=z+O(|z_{\text{\ref{theorem-1.2}}} |+ \| f _{\text{\ref{theorem-1.2}}} \| _{L _x ^{2,-2}} )$. So the asymptotic behavior of $  z' $ and of $  z _{\text{\ref{theorem-1.2}}}  $ is the same. By \eqref{eq:quasilin511} we get

\begin{equation}  \label{eq:quasilin611}  \begin{aligned} & \tau _{\widehat{D} }e^{\im \widehat{\vartheta}   \sigma _3}f_{\text{\ref{theorem-1.2}}}    = \tau _{\widehat{D}-\mathbf{A}} e^{  \sigma _3  (    \frac \im {2}{v \cdot x }  - \widetilde{\gamma}+\im \widehat{\vartheta} ) } f'
    - \tau _{\widehat{D}-\mathbf{A}} e^{  \sigma _3  (    \frac \im {2}{v \cdot x }  - \widetilde{\gamma}+\im \widehat{\vartheta} ) }\mathcal{G}.
 \end{aligned}
\end{equation}
By \eqref{eq:fl131} the second term on the rhs converges to 0 in $H^1$ as $ t\nearrow \infty$.
Hence, for $\theta := \widehat{\vartheta}+\im \widetilde{\gamma}
 $  and for $y  :=   \widehat{D}-\mathbf{A}  $, we obtain
 \eqref{scattering12}.

We   have
$
q\left ( \omega   (t) \right ) = q\left (\omega _0 \right )
-\frac{\| f'(t) \|_2^2 }{2}    + O( |z'(t)|+\| f'(t) \| _{L^{2,-2 }_x} )
$ by   $ q\left (\omega _0 \right )=q\left (\omega   \right ) +Q(R)$.
Then \eqref{scattering12} and  $|z'(t)|+\| f'(t) \| _{L^{2,-2 }_x}\to 0$
  imply
\begin{equation}   \begin{aligned}
 &   \lim _{t \to + \infty }
q\left ( \omega   (t) \right ) =
      q\left (\omega _0 \right ) - \lim _{t \to + \infty  }\frac{\|  e^{ {\im }t
\sigma _3   \Delta } f_+ \|_2^2 }{2}  = q\left (\omega _0 \right )
-\frac{\|   f_+  \|_2^2 }{2}=q(\omega _+ ),
\end{aligned} \nonumber
\end{equation}
 where  $\omega _+ $ is the unique element near $\omega _0$ for which the last  equality holds. So $ \lim _{t
\to + \infty } \omega (t) = \omega _+  $.
By  $  v =2   (\Pi  (U_0)- \Pi  (  R)  ) Q^{-1}(U_0)$ we obtain   $$  v =2(\Pi  (U_0)- \Pi  (  f')  ) Q^{-1}(U_0) +  O( |z'(t)|+\| f'(t) \| _{L^{2,-2 }_x} )$$
which implies $
     \lim _{t \nearrow \infty }
 v \left ( t\right ) =  2   (\Pi  (U_0)- \Pi  (  f_+ )  ) Q^{-1}(U_0) =:v_+
$.\qed

\begin{lemma}
\label{lem:scattpar} For  $(\theta ,y )$ the functions of   \eqref{scattering12}
and  $( {\vartheta}, {D})$ the  coordinates of \eqref{eq:anzatz},
there are   $\vartheta _0\in \R$, $y_0\in \R ^3$  and   $o (1)\to 0$ as $t\to +\infty$ s.t.

\begin{equation} \label{eq:scattpar4} \begin{aligned}&
 {\theta}  (t)  = \vartheta (t)  +\vartheta _0+o(1) \, , \quad y  (t)  = D (t)   + y_0  .
          \end{aligned}
\end{equation}

\end{lemma}
\proof  Consider the representation  $U=\tau _D e^{ \im \sigma _3 (\frac{v \cdot x}{2} +\vartheta  ) }  (\Phi _ \omega +R)$ of the solution of $\im \dot U = \sigma _3 \sigma _1 \nabla E(U)$. We have the identity
\begin{equation*} \begin{aligned}&\im \dot  U= -\sigma _3
(\dot \vartheta  - \frac{v\cdot \dot D}{2}) \tau _D e^{\im \sigma _3\Theta }
 (\Phi _\omega  + R)  -\im \dot D\cdot \tau _D e^{\im \sigma _3\Theta }
 \nabla _x(\Phi _\omega  + R)
 \\& -\frac{\dot v}{2}\cdot  \tau _D e^{\im \sigma _3\Theta }\sigma _3x
 (\Phi _\omega  + R) + \im \dot \omega\tau _D e^{\im \sigma _3\Theta }
  \partial _\omega \Phi _\omega +  \im   \tau _D e^{\im \sigma _3\Theta }
  \dot R.\end{aligned}
\end{equation*}
By    Lemma \ref{lem:gauge}  we have $ \nabla E(U)=$
\begin{equation*} \begin{aligned}& = \tau _D e^{ -\im \sigma _3 (\frac{v \cdot x}{2} +\vartheta  ) }  \left ( \nabla E( \Phi _ \omega +R) -v_a    \nabla \Pi _a( \Phi _ \omega +R)+ \frac{v^2}{4} \nabla Q( \Phi _ \omega +R) \right ) .
   \end{aligned}
\end{equation*}
Then, using also
\eqref{eq:charge1},
$\im \dot U = \sigma _3 \sigma _1 \nabla E(U)$ can be expanded for $\varpi =\omega _0$ into

\begin{equation}\label{eq:asympt} \begin{aligned}&  -\sigma _3
(\dot \vartheta  - \frac{v\cdot \dot D}{2} +\frac{v^2}{4} -\varpi )
 (\Phi _\omega  + R)  -\im (\dot D -v)\cdot
 \nabla _x(\Phi _\omega  + R)+  \im
  \dot R
 \\&  -\frac{\dot v\cdot x}{2}    \sigma _3
 (\Phi _\omega  + R) + \im \dot \omega
  \partial _\omega \Phi _\omega =\sigma _3 \sigma _1 \left (  \nabla E( \Phi _ \omega +R)
   + \varpi Q ( \Phi _ \omega +R)\right ).
          \end{aligned}
\end{equation}
Using the first system of coordinates \eqref{eq:coordinate}, and denoting
 the $f$--coordinate by $f'$,
we have for $\| G_1 \| _{L^\infty _t ([0,+\infty ), L^1_x)}\le C \epsilon ^2$

\begin{equation} \label{eq:scattpar1} \begin{aligned}&  -\sigma _3
(\dot \vartheta  - \frac{v\cdot \dot D}{2} +\frac{v^2}{4} -\omega _0)
 f '-\im (\dot D -v)\cdot
 \nabla _xf'
 -\frac{\dot v\cdot x}{2}    \sigma _3
 f '+  \im
  \dot f ' \\& =\sigma _3   (-\Delta +\omega _0) f' + G_1.
          \end{aligned}
\end{equation}
  Now we substitute  in \eqref{eq:scattpar1} the variables of the last coordinate system.
In particular $f'$ and $f$ are related by a formula like \eqref{eq:quasilin511}:
\begin{equation} \label{eq:scattpar2}
\begin{aligned} &    f'(x) =  e^{  \im\sigma _3  (  - \frac{1}{2} v \cdot (x- \mathbb{A} )   -\im     \widehat{\gamma}  ) } f(x- \mathbb{{A}} )
    +\mathcal{G} (x) \ , \\&  \widehat{\gamma}  =\widetilde{\gamma} +\sum _{\ell =1}^{2N+1} \gamma _\ell  \, , \quad   \mathbb{{A}}=\mathbf{A}+
\sum _{\ell =1}^{2N+1}  \mathbf{A} _ \ell  \, ,
\end{aligned}
\end{equation}
with the $(\widetilde{\gamma} ,\mathbf{A})$ of \eqref{eq:quasilin611} and the $(  \gamma _\ell, \mathbf{A} _ \ell)$
of Lemma \ref{lem:part-1.2}.   Substituting  \eqref{eq:scattpar2} in \eqref{eq:scattpar1}
we get after various cancelations, for $\| G_2\| _{L^\infty _t ([0,+\infty ), L^1_x\cap H^1_x)}\le C \epsilon  $,

\begin{equation} \label{eq:scattpar3} \begin{aligned} - &  \sigma _3
 (\dot \vartheta  -\omega _0)
 f -\im (\dot D  +\dot {\mathbb{A}})\cdot
 \nabla  f  +  \im
  \dot f     =\sigma _3   (-\Delta +\omega _0) f  + G_2.
          \end{aligned}
\end{equation}
  We claim that equation
\eqref{eq:scattpar3} is equivalent to equation \eqref{eq:f variable}.  Indeed,  taking their difference
we have
\begin{equation*}  a  _0(t)\sigma _3  f +\im  a_j(t)  \partial _{x_j} f= \textbf{G}
\end{equation*}
with  $\textbf{G} $  (resp .  $a_j $ )   a continuous  functional with values $  L^{\infty}(\R,  L^1(\R ^3)
\cap H^1_x) $
 (resp .  $  L^{\infty}(\R ) $ )  bounded in the space of solutions we are considering.  Then $ a  _0(t) \int  f (t,x)dx =\sigma _3\int   \textbf{G}  dx $.  If  $a  _0(t_0)\neq 0$ for a given   solution, we can find solutions for which
$f_n (t,\cdot ) \in {\mathcal S}(\R ^3)$,  $ f_n (t_0,\cdot )\to f (t_0,\cdot )$ in $H^1(\R ^3)$,
$\|  f_n (t_0  ) \|  _{L^1(\R ^3)}\nearrow \infty$,
$\textbf{G}  _{ n}(t_0)\to \textbf{G}  (t_0)$
and $a  _{0n}(t_0)\to a  _{0}(t_0)$.  This yields a contradiction. So    $a  _0(t)\equiv 0$.  By similar reasons
  $   a_j(t) \equiv 0$. This implies    $   \textbf{G}   \equiv 0$.
Equivalence of  \eqref{eq:scattpar3} and \eqref{eq:f variable} yields
  \begin{equation} \label{eq:scattpar31} \begin{aligned}&  \dot \vartheta   + \frac{1}{2}\frac{d}{dt} (v\cdot \mathbb{A})-\im \dot {\widehat{\gamma}} = \dot {\chi}   \, , \quad  \dot D _a +\dot {\mathbb{A}}_a= \partial _{ \Pi _a(f)} H \ ,
          \end{aligned}
\end{equation}
with the first one a consequence of $ \dot \vartheta  -\omega _0 =  \partial _{ Q(f)} H $.
 Using the notation of Lemmas \ref{lem:asymptotically f}--\ref{lem:vomega as},  \eqref{eq:scattpar31} yields what follows:
 \begin{equation}  \label{eq:scattpar32} \begin{aligned}&
 \dot \theta =\dot \vartheta   + \frac{1}{2}\frac{d}{dt} (v\cdot \mathbb{A})\, ;
  \\& \dot D    +\dot {\mathbb{A}} = \dot X =\dot {\widehat{D}} +
\sum _{\ell =1}^{2N+1} \dot { \mathbf{A} } _ \ell =\dot y   +\dot {\mathbb{A}}
\text{ and so }  \dot y   =\dot D    .
          \end{aligned}
\end{equation}
Notice that by Lemma \ref{lem:quasilin5} and inequality  \eqref{eq:systchi120}
there exists $\lim _{t\infty}\mathbb{A} (t)$. By Lemma \ref{lem:vomega as} we have $\lim _{t\nearrow \infty}v (t)=v _+$. This proves the existence of $\vartheta _0$
in \eqref{eq:scattpar4}, which is then proved.
    \qed

\begin{lemma}
\label{lem:Conv der}    The functions  $(\vartheta  ,D)$ in Lemma \ref{lem:scattpar} satisfy $\dot D =v+o(1)$  and $\dot \vartheta =\omega + \frac{v^2}{4}+o(1)$, with  $\lim _{t\to \infty}o(1)=0$.
\end{lemma}
\proof Consider \eqref{eq:asympt} with $\varpi =\omega$. Applying to it the linear operator $  | x_a \Phi _\omega \rangle $ and using the $(z,f)$ of the initial coordinate system \eqref{eq:coordinate}
we get $|\dot D_a-v_a|\le C (|z|+\| f \| _{L ^{2,-S}}) $ for arbitrary $S$ and fixed $C$. So the rhs is o(1) and we conclude  $\dot D =v+o(1)$.
 Applying     $  |  \sigma _3\partial _\omega \Phi _\omega \rangle $   to \eqref{eq:asympt} we get similarly
$|\dot \vartheta  - \frac{v\cdot \dot D}{2} +\frac{v^2}{4} -\omega|\le C (|z|+\| f \| _{L ^{2,-S}})$. Using $\dot D =v+o(1)$ we conclude $\dot \vartheta =\omega + \frac{v^2}{4}+o(1)$.

\qed

\subsection{Steps (ii) and (iii): the Fermi golden rule}
\label{subsec:FGR}

Step (ii) in the proof of Theorem \ref{proposition:mainbounds}
consists in introducing the variable

\begin{equation}
  \label{eq:g variable}
g=f+ \sum _{|\lambda ^0 \cdot(\mu-\nu)|>\omega _0} z^\mu
\overline{z}^\nu
   R ^{+}_{\mathcal{H}}  (\lambda  ^0 \cdot(\mu-\nu) )
    G_{\mu \nu}^0 .
\end{equation}
Substituting the new variable  $g$ in \eqref{eq:f variable}, the
first line on the rhs of  \eqref{eq:f variable} cancels out.  The following
result has been proved in a variety of places in the absence of translation, see for example \cite{boussaidcuccagna} .
Thanks to Sect. 3.3 \cite{beceanu}, essentially the same proof holds here. We skip the proof.

\begin{lemma}\label{lemma:bound g} For $\epsilon$ in Theorem \ref{theorem-1.2} sufficiently small, for $C_0=C_0(\mathcal{H})$  a fixed constant, we have $\| g
\| _{L^2_tL^{2,-S}_x}\le C_0 \epsilon + O(\epsilon
^2)$
for a fixed $S>1$.
\end{lemma}
  We have arrived at  step (iii) of the proof of Theorem \ref{proposition:mainbounds}: the Fermi Golden Rule.

   \begin{proposition}\label{prop:FGR} There is a  new set of variables
   $\zeta =z+O(z^2)$ such that     for a fixed $C$
   we have

\begin{equation}  \label{equation:FGR3} \begin{aligned}   & \| \zeta  -
 z  \| _{L^2_t}
\le CC_2\epsilon ^2\, , \quad  \| \zeta  -
 z \| _{L^\infty _t} \le C \epsilon ^3
\end{aligned}
\end{equation}
  and we have
\begin{equation} \label{eq:FGR5} \begin{aligned}
 &\partial _t \sum _{j=1}^m \lambda  _j ^0
 | \zeta _j|^2  =  2  \sum _{j=1}^m \lambda _j^0\Im \left (
\mathcal{D}_j \overline{\zeta} _j \right ) -\\&    -2 \sum _{
\substack{ \lambda ^0 \cdot  \alpha =\lambda ^0\cdot   \nu   >
\omega _0
   \\ \lambda
\cdot  \alpha -\lambda _k   < \omega _0   \, \forall \, k \, \text{
s.t. } \alpha _k\neq 0\\ \lambda \cdot  \nu -\lambda _k  <\omega _0
\, \forall \, k \, \text{ s.t. } \nu _k\neq 0}} \lambda ^0\cdot \nu
\Im \left ( \zeta ^{ \alpha } \overline{\zeta }^ { \nu  } \langle
R_{ \alpha 0}^+  G_{ \alpha 0}^0|  \sigma _1\sigma _3  G ^0 _{0\nu
}\rangle \right )
\end{aligned}
\end{equation}
     where $ \sum _j\|\mathcal{D}_j \overline{\zeta} _j\|_{
L^1[0,T]}\le (1+C_2)c_0 \epsilon ^{2}$ for a
fixed constant $c_0$.

   \end{proposition}
   \proof See \cite{Cu1} and  \cite{Cu4}.\qed

For the sum in the second line of \eqref{eq:FGR5}
we get finite sums

\begin{equation} \label{eq:FGR8} \begin{aligned} & 2\sum _{\Lambda>\omega
_0 } \Lambda
    \Im \left   \langle R_{
\mathcal{H}}^+ (\Lambda )\sum _{  \lambda ^0\cdot \alpha =\Lambda }\zeta ^{
\alpha } G _{  \alpha 0}^0| \sigma _1 \sigma _3\sum _{  \lambda
^0\cdot \nu =\Lambda}\overline{\zeta} ^{ \nu } G ^0_{0\nu
 }  \right \rangle      =\\&  2\sum _{\Lambda>\omega
_0 } \Lambda     \Im \left    \langle R_{ \mathcal{H}}^+ (\Lambda )\sum _{
\lambda ^0\cdot \alpha =\Lambda }\zeta ^{ \alpha } G _{  \alpha 0}^0|
\sigma _3\overline{\sum _{  \lambda ^0\cdot \alpha =\Lambda}\zeta ^{
\alpha } G ^0_{  \alpha 0} }\right \rangle    ,
\end{aligned}
\end{equation}
where we have used $G_{\mu \nu }^0=-\sigma _1 \overline{G^0}_{  \nu
\mu} $. Notice that the existence of $R_{ \mathcal{H}}^+(\Lambda )$  for $\Lambda \in \sigma _e(\mathcal{H})$ is proved in
\cite{CPV}.

 We have:

\begin{lemma}[Semipositivity]
\label{lemma:FGR8} We have rhs\eqref{eq:FGR8}$\ge 0.$
\end{lemma}
\proof   See \cite{Cu1}.
\qed

Now we will assume the following hypothesis.

\begin{itemize}
\item[(H10)]  We assume
that for some fixed constants for any vector $\zeta \in
\mathbb{C}^n$ we have:
\begin{equation*} \label{eq:FGR} \begin{aligned} &  \sum _{ \substack{
(\alpha , \nu ) \text{ as} \\ \text{in   \eqref{eq:FGR5}}}}
\lambda ^0\cdot \nu
\Im \left ( \zeta ^{ \alpha } \overline{\zeta }^ { \nu  } \langle
R_{ \alpha 0}^+  G_{ \alpha 0}^0|  \sigma _1\sigma _3  G ^0 _{0\nu
}\rangle \right )
 \approx \sum _{ \substack{ \lambda ^0\cdot  \alpha
> \omega _0
\\
   \lambda ^0
\cdot  \alpha -\lambda _k  ^0  < \omega _0 \\ \forall \, k \, \text{
s.t. } \alpha _k\neq 0}}  | \zeta ^\alpha  | ^2 .
\end{aligned}
\end{equation*}

\end{itemize}

 By (H10) we have

\begin{equation*} \label{eq:FGR10} \begin{aligned} &
2\sum _j \lambda _j^0\Im \left ( \mathcal{D}_j   \overline{\zeta} _j
\right )\gtrsim \partial _t \sum _j \lambda _j^0| \zeta _j|^2  +
      \sum _{ \substack{ \text{$\alpha$ as in (H10)}}}
     | \zeta ^\alpha  | ^2.
\end{aligned}
\end{equation*}
Then, for $t\in [0,T]$,  by the  last line in Proposition \ref{prop:FGR} we have

\begin{eqnarray}& \sum _j \lambda _j ^0 | \zeta
_j(t)|^2 +\sum _{ \substack{ \text{$\alpha$ as in (H10)}}}  \| \zeta ^\alpha \| _{L^2(0,t)}^2\lesssim
\epsilon ^2+ C_2\epsilon ^2.\nonumber
\end{eqnarray}
By \eqref{equation:FGR3} this implies $\| z ^\alpha \|
_{L^2(0,t)}^2\lesssim \epsilon ^2+ C_2\epsilon ^2$ for all the above
multi indexes. So, from  $\| z ^\alpha \| _{L^2(0,t)}^2\le
  C_2^2\epsilon ^2$ we conclude $\| z ^\alpha \|
_{L^2(0,t)}^2\le \kappa  C_2\epsilon ^2$ for a fixed $\kappa$, which is an improvement if $C_2$ is sufficiently large.
 So if we take $C_1>2K_1(C_2)$  and  $C_1$-- $C_2$ sufficiently large,
 in particular so that $ \sqrt{\kappa C_2}<C_2/2$, we conclude as desired that
 \eqref{4.4a} and
\eqref{4.4} imply the same estimate but with $C_1$, $C_2$ replaced
by $C_1/2$, $C_2/2$.   This yields Theorem \ref{proposition:mainbounds}.

\begin{remark}
\label{rem:H11} Suppose for simplicity that $ \lambda
_1(\omega )=   ...= \lambda _m(\omega )=:\lambda (\omega )$ and let us see the meaning of (H10). We have  $ G ^0_{  \alpha 0}  \in \mathcal{S}(\R ^3, \C ^2)$ for all $\alpha$.  For
   $W(\omega ) =\lim_{t\to +\infty}e^{-\im t \mathcal{H}_\omega }e^{\im t\sigma_3
(-\Delta+\omega  )}$, there exist  $ F_{  \alpha }  \in W^{k,p}(\R ^3, \C ^2)$
 for all $k\in \R$ and $p\ge 1$ with
 $ G ^0_{  \alpha 0} =W(\omega _0)F_\alpha $, \cite{Cu6}. Let $^t{ {F} _\alpha }=( {F} ^{(1)} _\alpha , {F}^{(2)}  _\alpha )  $.
Then the left hand side  of (H10) can be expressed as

\begin{equation} \label{eq:FGR81} \begin{aligned} &  \int _{|\xi | =\sqrt{(N+1) \lambda _0-\omega _0 } } \big | \sum _{|\alpha |=N+1} \zeta ^{\alpha}\widehat{ {F}}_\alpha ^{(1)} (\xi )\big |^2 dS(\xi )   ,
\end{aligned}
\end{equation}
where we are taking the   standard Fourier transform and $dS(\xi ) $
is the standard measure on a sphere. \eqref{eq:FGR81} is
equivalent to the linear independence of the finite family of functions $\{ \widehat{ {F}}_\alpha ^{(1)} \} _{\alpha }$ on the sphere of radius $\sqrt{(N+1) \lambda _0-\omega _0} $. This independence is in general
  expected to be true. This point is discussed in
\cite{zhouweinstein1} for a special situation involving small ground states and  $N=1$.
\end{remark}

Department of Mathematics and Computer Sciences,  University of Trieste, via Valerio  12/1  Trieste, 34127  Italy

{\it E-mail Address}: {\tt scuccagna@units.it}

\end{document}